\documentclass[11pt]{article}
\usepackage{amsmath,amsfonts,amssymb,amsthm,bbm,mathptm,cases}
\usepackage{graphicx,graphics,epsfig,subfigure}
\usepackage{graphicx}
\usepackage{booktabs}
\usepackage{threeparttable}
\usepackage{appendix}
\usepackage{mathrsfs}
\usepackage{graphicx}
\usepackage{subfigure}

\usepackage{xcolor,color,soul,array,cite,float,version,times}
\numberwithin{equation}{section}
\numberwithin{figure}{section}
\numberwithin{table}{section}
\numberwithin{footnote}{section}

\theoremstyle{definition}

\newtheorem{sch}{Scheme}[section]
\theoremstyle{plain}
\newtheorem{thm}{Theorem}[section]

\theoremstyle{remark}

\newcommand{\ben}{\begin{eqnarray}}
\newcommand{\een}{\end{eqnarray}}
\newcommand{\bea}{\begin{array}}
\newcommand{\eea}{\end{array}}
\newcommand{\bes}{\begin{subequations}}
\newcommand{\ees}{\end{subequations}}
\newcommand{\bef}{\begin{figure}[H]}
\newcommand{\eef}{\end{figure}}
\newcommand{\bet}{\begin{tikzpicture}}
\newcommand{\eet}{\end{tikzpicture}}
\newcommand{\beq}{\begin{equation}}
\newcommand{\eeq}{\end{equation}}
\def\bena#1\eena{\begin{eqnarray}\begin{array}{l}#1\end{array}\end{eqnarray}}
\def\besl#1\eesl{\begin{subequations}\begin{align}#1\end{align}\end{subequations}}

\newcommand{\parl}[2]{\ensuremath{\frac{\partial #1}{\partial #2}}}

\newcommand{\vparl}[2]{\ensuremath{\frac{\delta #1}{\delta #2}}}

\newcommand{\ohs}[1]{\ensuremath{\overline{#1}^{n+1/2}}}
\newcommand{\hs}[1]{\ensuremath{{#1}^{n+1/2}}}

\def\inc(#1){\includegraphics[height=3 cm]{pics/#1}}

\def\br{\mathbf{r}}

\def\bI{\mathbf{I}}

\begin{document}
%\begin{CJK*}{GBK}{song}
%ÖÐÎIJâÊÔ
%\end{CJK*}
\title{Second Order Linear Energy Stable Schemes for Allen-Cahn Equations with Nonlocal Constraints}
\author{
{Xiaobo Jing}\footnote{Beijing Computational Science Research Center, Beijing 100193, P. R. China. },
{Jun Li}\footnote{nkjunli@foxmail.com, School of Mathematical Sciences, Tianjin Normal University, Tianjin 300387, P. R. China; Beijing Computational Science Research Center, Beijing 100193, P. R. China. },
{Xueping Zhao} \footnote{xzhao@math.sc.edu, Department of Mathematics, University of South Carolina, Columbia, SC 29208, USA. } and
{Qi Wang}\footnote{qwang@math.sc.edu, Department of Mathematics,
 University of South Carolina, Columbia, SC 29028, USA;  Beijing Computational Science Research Center, Beijing  100193, P. R. China.}
%%{Xiaogang Yang}\footnote{xgyang@csrc.ac.cn, Beijing Computational Science Research Center, Beijing  100193, P. R. China.}
%%{M. Gregory Forest}\footnote{forest@unc.edu, Departments of Mathematics and Biomedical Engineering, University of North Carolina at Chapel Hill,
%%Chapel Hill, NC 27599, USA.}
}
\date{\today}
\maketitle
\begin{abstract}
%The Allen-Cahn equation with a nonlocal constraint can be used as an alternative to the Cahn-Hilliard equation in describing phase behavior in immiscible material mixtures or on its own for relaxation dynamics with a nonlocal free energy. In this paper,
We present a set of linear, second order, unconditionally energy stable schemes for the Allen-Cahn equation with nonlocal constraints that preserves the total volume of each phase in a binary material system. The energy quadratization strategy is employed to derive the energy stable semi-discrete numerical algorithms in time. Solvability conditions are then established for the linear systems resulting from the semi-discrete, linear schemes. The fully discrete schemes are obtained afterwards by applying second order finite difference methods on cell-centered grids in space. The performance of the schemes are assessed against two benchmark numerical examples, in which dynamics obtained using the volume-preserving Allen-Cahn equations with nonlocal constraints is compared with those obtained  using the classical Allen-Cahn as well as the Cahn-Hilliard model, respectively, demonstrating slower dynamics when volume constraints are imposed as well as their usefulness as alternatives to the Cahn-Hilliard equation in describing phase evolutionary dynamics for immiscible material systems while preserving the phase volumes.
Some performance enhancing, practical implementation methods for the  linear energy stable schemes are discussed in the end.
\end{abstract}

\section{Introduction}

\noindent \indent Thermodynamically consistent models for material systems represent a class of models that yield the energy dissipation law. One particular class of the models is known as the gradient flow, in which the time evolution of thermodynamical variables is proportional to the variation of the system free energy. When the thermodynamic variable is a phase variable, it's known as the  Allen-Cahn equation. This class of models describes relaxation dynamics of the thermodynamical system to equilibrium. There are many applications of such gradient flow equations, particularly in the materials science, life science and fluid dynamics \cite{allen1979microscopic, cahn1958free, gurtin1996two, doi1988theory, gurtin1996two, Leslie1979Theory, Gong2018Fully, Pethrick1988The, Gong2018Linear, Yang2017Linear, Zhao2016Numerical, Zhao2017Numerical, wise2009energy, yang2017hydrodynamic, cheng2015fast, guo2016local}. However, in the case of a phase field description,  when the phase variable represents the volume fraction of a material component, this model does not warrant the conservation of the volume of that component. In order to preserve the volume, the free energy functional has to be augmented by a volume preserving mechanism with a penalizing potential term or a Lagrange multiplier \cite{li2018unconditionally}\cite{Rubinstein1992Nonlocal }. This often results in a nonlocal term in the modified Allen-Cahn equation. In this paper, we call these nonlocal Allen-Cahn equations or Allen-Cahn equations with nonlocal constraints.

The Cahn-Hilliard equation is an alternative  model for the gradient flow. One feature of the Cahn-Hilliard equation is its volume preserving property. Rubinstein and Sternberg studied the Allen-Cahn model with the volume constraint analytically and compared it with the Cahn-Hilliard model \cite{Rubinstein1992Nonlocal}.  Their result favors using the Allen-Cahn model with a volume constraint in place of the Cahn-Hilliard model when studying interfacial dynamics of incompressible, immiscible multi-component material systems.  For the Allen-Cahn equation as well as the Cahn-Hilliard equation, there have been several popular numerical approaches to construct energy stable schemes for the equations, including the convex splitting approach \cite{elliott1993global, fan2017componentwise, eyre1998unconditionally, wise2009energy, shen2012second, wang2010unconditionally, han2015second } , the stabilizing approach \cite{shen2010numerical, du2018stabilized, chen2017uniquely}, the energy quadratization (EQ) approach \cite {zhao2018general, yang2017numerical, Gong2018Fully, Gong2018Linear} and the scalar auxiliary variable approach, which is a special form of EQ method, \cite {shen2018scalar, dong2018family, wang2017efficient}.  Recently, the energy quadratization (EQ) and its reincarnation in the scalar auxiliary variable (SAV) method have been applied to a host of thermodynamical and hydrodynamic models owing to their simplicity, ease of implementation, computational efficiency, linearity, and most importantly their energy stability property \cite{Zhaoetal2018-2, yang2017numerical, Gong2018Fully, Gong2018Linear, Yang2017Linear, Zhao2016Numerical, Zhao2017Numerical, chen2018regularized, yang2018efficient, zhao2016energy, zhao2017novel, zhao2016decoupled, dong2018family, wang2017efficient}. We have shown that these strategies are general enough to be useful for developing energy stable numerical approximations to any thermodynamically consistent models, i.e., the models derived using the second law of thermodynamics or the Onsager principle \cite{onsager1931reciprocal1, onsager1931reciprocal2, zhao2018general}.

In this paper,  we develop a set of linear, second order, unconditionally energy stable schemes using the energy quadratization (EQ)  and scalar auxiliary variable (SAV) approach to solve the Allen-Cahn equation, nonlocal Allen-Cahn equation and Cahn-Hilliard equations numerically. The numerical schemes for the Allen-Cahn and the Cahn-Hilliard model are not new. They are presented here for comparison purposes.  However, the schemes for the nonlocal Allen-Cahn models are new. In some of these schemes, both EQ and SAV approaches are combined to yield linear, energy stable schemes. We note that when a nonlocal Allen-Cahn model is discretized, it is inevitable to yield an SAV scheme when the EQ strategy is applied. When multiple integrals are identified as SAVs in the free energy functional, new iterative steps are proposed to solve the subproblems in which elliptic equations are solved efficiently. All these schemes are linear and second order accurate in time. The linear system resulting from the schemes are all solvable uniquely so that the solution existence and uniqueness in the semi-discrete system is warranted. When EQ is coupled with the discretized integrals, the Sherman-Morrison formula can lead us to an efficient numerical scheme. This can be equivalently be dealt with using the SAV method as well. The numerical schemes developed for the Allen-Cahn equation with  nonlocal volume-preserving constraints also preserves the volume at the discrete level in addition to preserving the energy dissipation rate. In the end, we conduct two numerical experiments to assess the performance of the schemes. The results  based on EQ and those based on SAV perform equally well in preserving the volume and energy dissipation rate.  In addition, the computational efficiency of the schemes is comparatively studied in one of the benchmark examples. Some performance enhancing, practical implementation methods are discussed to improve the accuracy of the schemes at large time step sizes. To simplify our presentation, we present the temporal discretization of the models using EQ and SAV approaches in detail. Then, we briefly discuss our strategy to obtain fully discrete schemes by discretizing the semi-discrete schemes in space, and refer readers to our early publications in \cite{Gong2018Fully,Gong2018Linear} for more details.

The rest of paper is organized as follows. In $\S 2$, we present the mathematical models for the Allen-Cahn, the Allen-Cahn with nonlocal constraints, and the Cahn-Hilliard equation. In $\S 3$, we study their near equilibrium dynamics. In $\S 4$, we present a set of second order, linear, energy stable numerical schemes for the models. In $\S 5$, we conduct mesh refinement tests on all the schemes and carry out two simulations on drop merging as well as phase coarsening experiments using the models. Finally, we give the concluding remark in section $\S 6$.

\section {Phase Field Models for Binary Materials}

\noindent \indent We briefly review two simple phase field models for a binary material system: the Allen-Cahn and the Cahn-Hilliard model, in which the free energy density of the binary material system is given by a functional of phase variable $\phi \in [0,1]$ and its gradients. For instance, to study drops of one fluid within the matrix of the other immiscible fluid while ignoring hydrodynamic effects, the free energy is customarily chosen as the following double-well potential:
\bena
F=\int_{\Omega} \gamma[\frac{\epsilon}{2}{(\nabla\phi)}^2+\frac{1}{\epsilon}\phi^2(1-\phi)^2]\mathrm{d{\bf r}},
\eena
where $\Omega$ is the material domain, $\epsilon$ is a parameter describing the width of the interface and $\gamma$ is proportional to the surface tension.
In general, the generic form of some commonly used free energies are given as follows
\bena
F=\int_{\Omega}[\gamma_1{(\nabla\phi)}^2+f(\phi)]\mathrm{d{\bf r}},
\eena
where $\gamma_1$ parametrizes the conformational entropy and  $f(\phi)$ is the bulk potential.

Dynamics of the binary material system is customarily governed by a time dependent partial differential equation model resulting from the Onsager's linear response theory\cite{onsager1931reciprocal1,onsager1931reciprocal2},
\bena
\parl{\phi}{t}= -M \mu, in \; \Omega,\\
\mu=\vparl{F}{\phi}=-2\gamma_1\nabla^2\phi+f'(\phi),
\eena
subject to appropriate boundary and initial conditions, where $M$ is the mobility matrix consisting of differential operators of even order and  $\mu$ is the chemical potential.

The time rate of change of the free energy, known as the energy dissipation functional, is given by
\ben
\frac{dF}{dt}=-\int_{\Omega} \mu M \mu \mathrm{d{\bf r}}+\int_{\partial \Omega} {\bf n}\cdot \frac{\partial f}{\partial \nabla \phi}\phi_t ds.\label{Energy-dissipation}
\een
The no-flux boundary condition
\ben
{\bf n}\cdot \frac{\partial f}{\partial \nabla \phi}=0
\een
annihilates the energy dissipation at the boundary, where $\bf n$ is the unit external normal of the boundary. This is the commonly used boundary condition. Another less-used boundary condition is
\ben
\phi_t=-\beta {\bf n}\cdot \frac{\partial f}{\partial \nabla \phi}.
\een
This specifies a decay rate of the phase variable at the boundary for $\beta \geq 0$.

 The system is dissipative with respect to both boundary conditions provided $M$ is a positive definite operator.
 The commonly used phase field models such as Allen-Cahn and Cahn-Hilliard are two special cases,
corresponding to
\ben
M=\left \{
\bea{ll}
M_0, & \hbox{Allen-Cahn},\\
-\nabla \cdot M_0 \nabla, & \hbox{Cahn-Hilliard},
\eea\right.
\een
respectively, where $M_0$ is a prescribed mobility coefficient matrix, which is a function of $\phi$. However, in some case, a constant mobility is used as an approximation instead. The Allen-Cahn equation defined this way does not conserve the total volume
$\int_\Omega \phi \mathrm{d\br}$ if $\phi$ is the volume-fraction while the Cahn-Hilliard equation does. However, these two models predict  similar near equilibrium dynamics revealed in their linear stability analyses below. On the other hand, the Allen-Cahn equation is an equation of lower spatial derivatives, and presumably costs less when solved numerically. Thus, one opted to  impose the volume conservation as a constraint to the Allen-Cahn equation for it to be able to describe dynamics in which the volume is conserved in some cases. Next, we will briefly recall several ways to impose  volume conservation to dynamics described by the Allen-Cahn equation, then discuss how to design efficient and energy stable numerical algorithms to compute their numerical solutions.

\subsection{Allen-Cahn model}

\noindent \indent The  Allen-Cahn equation with the no-flux Neumann boundary condition and initial conditions is given as follows:
\bena
\parl{\phi}{t}= -M \mu,  in \; \Omega \\ \label{BAC}
\parl{ \phi}{n}=0, in \;\partial\Omega, \\
\phi|_{t=0}=\phi(\br, 0), \br \; in \; \Omega.
\eena
%\begin{align*}
%\begin{cases}
%\parl{\phi}{t}= -M \mu,  & \text{in} \;\Omega \\
%\parl{ \phi}{n}=0, & \text{in} \;\partial\Omega \\
%\phi|_{t=0}=\phi(0),
%\end{cases}
%\end{align*}

The energy dissipation rate of the Allen-Cahn equation is given by \eqref{Energy-dissipation}
%\bena
%\frac{d F}{dt}=\int_{\Omega} \vparl{F}{\phi}\phi_t \mathrm{d{\bf r}}
%=-\int_{\Omega} \mu (M \mu) \mathrm{d{\bf r}} \leq 0,
%\eena
%provided $M$ is nonnegative.
We denote the volume of $\Omega$ by
\ben
V(t)=\int_{\Omega} \phi \mathrm{d{\bf r}}.
\een
 Then,
\ben
\frac{dV}{dt}=-\int_{\Omega} M \mu \mathrm{d{\bf r}}.
\een
It is normally nonzero provided $M\neq 0$, which implies that $V(t)$ is not conserved. One simple fix is to enforce the volume constraint $V(t)=V(0)$ by coupling it to the Allen-Cahn dynamical equation. The resulting equation is termed the Allen-Cahn equation with nonlocal constraints.

\subsection {Allen-Cahn models with nonlocal constraints}

\noindent \indent In addition to the approximate volume defined above,  $V(t)=\int_\Omega \phi(t) \mathrm{d{\bf r}}$,  we can introduce a more general definition using a function $h(\phi)$,
which is (i) monotonically increasing for $\phi \in [0,1]$, and (ii) $h(0)=0, h(1)=1$ as follows,
\ben
V(t)=\int_{\Omega} h(\phi) d\br.
\een

 We next discuss methods to enforce volume conservation for the Allen-Cahn model. First, we consider the method that  minimizes $(V(t)-V(0))^2$ by penalizing it in the free energy functional.

\subsubsection{Allen-Cahn model with a penalizing potential}

\noindent \indent Here, we  augment the free energy with  a penalizing potential term as follows
\bena
F=\int_{\Omega} [\gamma_1{(\nabla\phi)}^2+f(\phi)]\mathrm{d{\bf r}}+\frac{\eta}{2}(V(t)-V(0))^2, \label{PPM}
\eena
where $\eta$ is the penalizing parameter, a large positive constant.
The transport equation for $\phi$ is given by the Allen-Cahn equation with the modified chemical potential
%\bena
%\parl{\phi}{t}= -M \tilde \mu,\quad  in \; \Omega \\ \label{BPAC}
%\parl{ \phi}{n}=0, in \; \partial \Omega \\
%\phi|_{t=0}=\phi(0),
%\eena
%where $M$ is the mobility coefficient and $\tilde \mu$ is the chemical potential given by
\bena
\tilde \mu=\vparl{F}{\phi}
=\mu+\sqrt{\eta}\zeta h'(\phi),\quad
\zeta=\sqrt{\eta}(V(t)-V(0)).
\eena
The energy dissipation rate is given by
\bena
\frac{d F}{dt}=\int_{\Omega} \vparl{F}{\phi}\phi_t \mathrm{d{\bf r}}
=-\int_{\Omega} \tilde \mu (M \tilde \mu) \mathrm{d{\bf r}}. \label{EDISS}
\eena
It is negative if  $M\geq 0$. The modified Allen-Cahn equation is approximately volume-conserving depending on the size of $\eta$. In principle, an appropriate $\eta>0$ can make  $V(t)$  close to $V(0)$. The choice of $\eta$ is, however, up to  the user.
%{\bf If $h'(\phi)=\frac{(m+1)(2m+1)}{m}[\phi(1-\phi)]^m, m \hbox{ is a  positive integer}$, Since $\int_\Omega h(\phi) \mathrm{d{\bf r}} $ is nonlinear, which will leads to a nonlinear EQ or SAV scheme. The energy dissipation law seems like not inconsistent}

\subsubsection{Allen-Cahn models with a Lagrange multiplier}

\noindent \indent To enforce the volume conservation in the Allen-Cahn model exactly, we use a Lagrange multiplier $L$ for  constraint $V(t)-V(0)=0$.
This can be accomplished by augmenting a penalty term with a Lagrange multiplier $L$ in the free energy functional as follows:
\ben
\tilde F=F-L (V(t)-V_0). \label{FM}
\een
The transport equation for $\phi$ is given by the Allen-Cahn equation
%\bena
%\parl{\phi}{t}= -M \tilde \mu, in \; \Omega \\ \label{BLAC}
%\parl{\phi}{n}=0 \quad in \quad \partial \Omega, \\
%\phi|_{t=0}=\phi(0),
%\eena
%where $M$ is the mobility coefficient and
with a modified chemical potential $\tilde \mu$ given by
\bena
\tilde \mu=\vparl{F}{\phi}
=\mu-L h'(\phi),
\eena
where $h'(\phi)=\vparl{h(\phi)}{\phi}$ and $V(t)-V(0)=0$ is enforced.
%where
%\bena
%\mu=-2\gamma_1\nabla^2\phi+2\gamma_2\phi(2\phi-1)(\phi-1)\\
%\eena
Since the volume is conserved,
 \ben
\int_{\Omega}[h'(\phi) \frac{d\phi}{dt}]\mathrm{d{\bf r}}=\int_{\Omega}[ h'(\phi)M \tilde \mu]\mathrm{d{\bf r}}=0.
 \een
We arrive at
\bena
L=\frac{1}{\int_{\Omega} h'(\phi) M h'(\phi) \mathrm{d{\bf r}}}\int_{\Omega}[h'(\phi) M \mu]\mathrm{d{\bf r}}.
\eena
%We calculate the energy dissipation rate as follows
%\bena
%\frac{d F}{dt}=\int_{\Omega} \vparl{F}{\phi}\phi_t \mathrm{d{\bf r}}
%=\int_{\Omega} \tilde \mu (-M \tilde \mu) \mathrm{d{\bf r}} \leq 0,
%\eena
%provided $M$ is nonnegative definite.
The volume conserved Allen-Cahn equation is nonlocal because an integral is in the chemical potential as well as in the equation.
The choice of $h(\phi)$ in this paper includes the following two families.
\ben
\left \{
\bea{l}
h(\phi)=\phi,\\
h'(\phi)=\frac{(m+1)(2m+1)}{m}[\phi(1-\phi)]^m, m \hbox{ is a  positive integer}.
\eea\right.
\een
The energy dissipation is again given by  (\ref{EDISS}).
%The energy dissipation rate can be expressed in the original chemical potential as follows
%\bena
%\frac{dF}{dt}=\int_V \vparl{F}{\phi}\phi_t \mathrm{d{\bf r}}
%%=\int_V \mu (-M\mu +Mh'(\phi)L) \mathrm{d{\bf r}}
%%=\int_V(\mu-L+L) (-M\mu +ML) \mathrm{d{\bf r}}
%\\
%%=\int_V[(\mu-L)M (-\mu +L)+L (-M\mu +ML)] \mathrm{d{\bf r}}\\i
%%=-\int_V(\mu-L)M (\mu-L) \mathrm{d{\bf r}}+L\int_VM(-\mu +L)\mathrm{d{\bf r}}\\
%=-\int_V(\mu-h'(\phi)L){M}(\mu-h'(\phi)L) \mathrm{d{\bf r}}.
%\eena
%The nonlocal Allen-Cahn equation is not only volume conserving but also energy dissipative provided $M \geq 0$. %We next discuss how we approximate the  Allen-Cahn equations numerically.
%The two types of models are named as the Lagrangian model 1 and Lagrangian model 2 in this paper, respectively.

\subsection{Cahn-Hilliard model}

\noindent \indent In the Cahn-Hilliard model, the transport equation for $\phi$ is given by
\bena
\parl{\phi}{t}= \nabla \cdot (M \nabla \mu),   in \; \Omega \\ \label{BCH}
\parl{ \phi}{n}=0, \parl{\mu}{n}=0, in \; \partial \Omega \\
\phi|_{t=0}=\phi(\br, 0).
\eena
%where $M$ is the mobility coefficient and $\mu$ is the chemical potential.
% given by
%\bena
%\mu=-2\gamma_1\nabla^2\phi+2\gamma_2\phi(2\phi-1)(\phi-1)\\
%\eena
The energy dissipation rate of the Cahn-Hilliard equation is given by
\bena
\frac{d F}{dt}=\int_{\Omega} \vparl{F}{\phi}\phi_t \mathrm{d{\bf r}}
=-\int_{\Omega} (\nabla\mu) M (\nabla\mu) \mathrm{d{\bf r}} \leq 0,
\eena
provided  $M$ is nonnegative definite. We next compare the near equilibrium dynamics of  the modified Allen-Cahn equations and the Cahn-Hilliard equation.

\section{Near equilibrium dynamics}

\noindent \indent We conduct a linear stability analysis about a constant steady state in a rectangular domain to investigate near equilibrium dynamics of the phase field models mentioned above. Specifically, we perturb the steady state $\phi^{ss}$ of the equations by a small disturbance $\delta v$,
\bena
\phi=\phi^{ss}+\delta v. \label{BLS}
\eena
{For the Allen-Cahn model,}
substituting equation \eqref{BLS} into equation \eqref{BAC}, we get the linearized equation
\bena
\parl{\delta v}{t}= -M [-2\gamma_1\nabla ^2\delta v+f''(\phi^{ss})\delta v].
\eena
We solve it using the Fourier series method in domain $\Omega=[-\pi,\pi]^2$ with $\delta v= \sum_{k,l=0}^{\infty}a_{kl}(t)cos(kx)cos(ly)$ that satisfies the boundary condition. Then we obtain the following ordinary differential equation system for each single mode
\bena
%\sum_{k=0}^{\infty}\dot{a_k(t)}cos(kx)cos(ky)=-M \sum_{k=0}^{\infty}a_k(t)[4\gamma_1k^2+2\gamma_2(6{(\phi^{ss})}^2-6\phi^{ss}+1)]cos(kx)cos(ky),\\
\dot{a}_{kl}(t)=-M a_{kl}(t)[2\gamma_1(k^2+l^2)+f''(\phi^{ss})], k, l=0, \cdots, \infty.
\eena
In this system, instability may occur only if $f''(\phi^{ss})<0$ and $2\gamma_1(k^2+l^2)+f''(\phi^{ss})<0$ for some small $k,l$.

For the Allen-Cahn model with a penalizing potential, substituting equation (\ref{BLS}) into the transport equation, we get the linearized equation
\bena
\parl{\delta v}{t}= -M (-2\gamma_1\nabla ^2\delta v+f''(\phi^{ss})\delta v+\eta [h''(\phi^{ss})\delta v (\int_\Omega h(\phi^{ss}) \mathrm{d{\bf r}}-V(0))+h'(\phi^{ss})\int_\Omega h'(\phi^{ss})\delta v \mathrm{d{\bf r}}]).
\eena
We solve it analogously and obtain governing system of equations for the Fourier coefficients:
\bena
%\sum_{k=0}^{\infty}\dot{a_k(t)}cos(kx)cos(ky)=-M \sum_{k=0}^{\infty}a_k(t)[(4\gamma_1k^2+2\gamma_2(6{(\phi^{ss})}^2-6\phi^{ss}+1))cos(kx)cos(ky)\\
%-\eta\int_\Omega a_k(t)cos(kx)cos(ky)\mathrm{d{\bf r}}]=-M \sum_{k=0}^{\infty}a_k(t)cos(kx)cos(ky)[(4\gamma_1k^2+2\gamma_2(6{(\phi^{ss})}^2-6\phi^{ss}+1))].\\
\dot{a}_{kl}(t)=-M a_{kl}(t)[2\gamma_1(k^2+l^2)+f''(\phi^{ss})+\eta h''(\phi^{ss}) (\int_{\Omega} h(\phi^{ss}) \mathrm{d{\bf r}}-V(0)) +\eta h'(\phi^{ss})\int_\Omega h'(\phi^{ss}) \mathrm{d{\bf r}}\delta_{k0}\delta_{l0} ].
\eena
If $2\gamma_1(k^2+l^2)+f''(\phi^{ss})+\eta h''(\phi^{ss}) (\int_{\Omega} h(\phi^{ss}) \mathrm{d{\bf r}}-V(0)) +\eta  \int_\Omega (h'(\phi^{ss}))^2 \mathrm{d{\bf r}}\delta_{k0}\delta_{l0} <0$ for some small $k, l$, instability may occur. Compared with the Allen-Cahn model, this modified Allen-Cahn introduces a stabilizing  mechanism at the zero wave number $k=0$ and $l=0$ due to the average at the zero-wave number and
a potentially  destabilizing  mechanism for $h''(\phi^{ss}) < 0$ and stablizing  mechanism for $h''(\phi^{ss}) >
0$. For the simple  case $h(\phi) = \phi$, the only contribution is the stabilizing  mechanism at $k = 0$ and
$l = 0.$

For the Allen-Cahn model with a Lagrangian multiplier,
we denote $g(\phi)=h'(\phi)M(-2\gamma_1\nabla^2\phi+f'(\phi^{ss}))$.
Substituting equation (\ref{BLS}) into the transport equation, we get the linearized equation
\bena
\parl{\delta v}{t}= -M [-2\gamma_1\nabla ^2\delta v+f''(\phi^{ss})\delta v-\frac{h''(\phi^{ss})\delta v\int_\Omega g(\phi^{ss})\mathrm{d{\bf r}}}{\int_\Omega [h'(\phi^{ss})]M [h'(\phi^{ss})]\mathrm{d{\bf r}}}-\frac{h'(\phi^{ss})\int_\Omega g'(\phi^{ss})\delta v\mathrm{d{\bf r}}}{\int_\Omega [h'(\phi^{ss})]M [h'(\phi^{ss})]\mathrm{d{\bf r}}}+\\
2\frac{h'(\phi^{ss})\int_\Omega g(\phi^{ss})\mathrm{d{\bf r}}}{(\int_\Omega [h'(\phi^{ss})]M [h'(\phi^{ss})]\mathrm{d{\bf r}})^2}\int_\Omega [h''(\phi^{ss})]M [h'(\phi^{ss})]\delta v\mathrm{d{\bf r}}].
\eena
%Now we discuss the different linear stabilities between the two choices of $h'(\phi)$. \\
%{\bf Choice 1: $h'(\phi)=1$}
%\bena
%\parl{\delta v}{t}= -M (-2\gamma_1\nabla ^2\delta v+f''(\phi^{ss})\delta v-\frac{1}{\int_\Omega M \mathrm{d{\bf r}}}\int_\Omega %[Mf''(\phi^{ss})]\mathrm{d{\bf r}}. \\
%\eena
%We assume $\delta v= \sum_{k=0}^{\infty}a_kcos(kx)cos(ky)$, the domain $\Omega=[-\pi,\pi]^2$. Then we have
The dynamic equations for the Fourier coefficients are
\bena
\dot a_{kl}(t)=-M a_{kl}(t)[2\gamma_1(k^2+l^2)+f''(\phi^{ss})-\frac{h''(\phi^{ss})\int_\Omega g(\phi^{ss})\mathrm{d{\bf r}}}{\int_\Omega [h'(\phi^{ss})]M [h'(\phi^{ss})]\mathrm{d{\bf r}}}-(\frac{h'(\phi^{ss})\int_\Omega g'(\phi^{ss})\mathrm{d{\bf r}}}{\int_\Omega [h'(\phi^{ss})]M [h'(\phi^{ss})]\mathrm{d{\bf r}}}-\\
2\frac{h'(\phi^{ss})\int_\Omega g(\phi^{ss})\mathrm{d{\bf r}}}{(\int_\Omega [h'(\phi^{ss})]M [h'(\phi^{ss})]\mathrm{d{\bf r}})^2}\int_\Omega [h''(\phi^{ss})]M [h'(\phi^{ss})]\mathrm{d{\bf r}})\delta_{k0}\delta_{l0}].
\eena
For the simple case $h(\phi)=\phi$,
\bena
\dot a_{kl}(t)=-M a_{kl}(t)[2\gamma_1(k^2+l^2)+f''(\phi^{ss})-\frac{ \int_\Omega Mf''(\phi^{ss})\mathrm{d{\bf r}}}{\int_\Omega M \mathrm{d{\bf r}}}\delta_{k0}\delta_{l0}].
\eena
If $2\gamma_1(k^2+l^2)+f''(\phi^{ss})-\frac{ \int_\Omega Mf''(\phi^{ss})\mathrm{d{\bf r}}}{\int_\Omega M \mathrm{d{\bf r}}}\delta_{k0}\delta_{l0}<0$, for small $k, l$, instability may occur.  Compared with the Allen-Cahn dynamics, the only contribution of the Lagrange multiplier is introducing a correction at the zero wave number limit which may be stabilizing when $f''(\phi^{ss})<0$ and destabilizing otherwise.

For the Cahn-Hilliard model,
substituting equation (\ref{BLS}) into equation (\ref{BCH}), we get the linearized equation
\bena
\parl{\delta v}{t}= M [-2\gamma_1\nabla ^4\delta v+f''(\phi^{ss})\nabla ^2\delta v)].
\eena
Repeating the analysis, we have
\bena
%\sum_{k=0}^{\infty}\dot{a_k(t)}cos(kx)cos(ky)=-M \sum_{k=0}^{\infty}a_k(t)[4\gamma_1k^4+2\gamma_2\nabla ^2(6{(\phi^{ss})}^2-6\phi^{ss}+1)cos(kx)cos(ky))],\\
\dot{a}_{kl}(t)=-M a_{kl}(t)[2\gamma_1(k^2+l^2)+f''(\phi^{ss})](k^2+l^2).
\eena
If $kl[2\gamma_1(k^2+l^2)+f''(\phi^{ss})]<0$, instability may ensure. We note that the window of instability in the Cahn-Hilliard model is identical to that in the Allen-Cahn model. However, the rate of growth is different. These linear stability results dictate the initial transition of the solution towards or away from the steady state. For long time transient behavior of the solution, we have to resort to numerical computations.

We next discuss how to numerically approximate the model equations efficiently with  high order, linear, energy stable schemes.

\section{Numerical Approximations of the Phase Field  Models }

\noindent \indent We design  numerical schemes to solve the above nonlocal equations  to ensure that the energy dissipation property as well as the total volume conservation are respected.   We do it by  employing the energy quadratization strategy (EQ) and the scalar auxiliary variable approach (SAV) developed recently \cite{Zhao2016Numerical, zhao2016energy, shen2018scalar, yang2018efficient}. Both methods depend on a reformulation of the model into one with a quadratic energy, and provide effective ways to design linear  numerical schemes. SAV is a specialized EQ strategy suitable for thermodynamical systems. For a full review on EQ methods on thermodynamical systems, readers are referred to a recent review article \cite{zhao2018general}. All schemes presented below are firstly semi-discretized by the  Crank-Nilcolson method in time and then centrally discretized in space later.
In fact, we have shown recently that BDF and Runge-Kutta methods can be used to design energy stable schemes for thermodynamical systems up to arbitrarily high order in time \cite{Zhaoetal2018-2}. For simplicity, we present the schemes in their semi-discrete forms in time. For comparison purposes, we also present analogous schemes for the classical Allen-Cahn and the Cahn Hilliard model as well.
% and efficient computation of the linear system of equations in the schemes.
%For completeness, we also recall the methods for Allen-Cahn model.

\subsection{Temporal discretization}

\subsubsection{Numerical methods for the Allen-Cahn model using EQ}

\noindent \indent We reformulate the free energy density by introducing an intermediate variable $q$ and a constant $C_0>0$ such that $q$ is a real variable,
\bena
q=\sqrt{f(\phi)-\gamma_2 \phi^2+C_0}.\label{Scheme1}
\eena
Then, the free energy is recast into a quadratic form:
\bena
F
=\int_{\Omega} [\gamma_1(\nabla\phi)^2+\gamma_2\phi^2+q^2-C_0]\mathrm{d{\bf r}}.
\eena
The chemical potential is give by
\bena
\mu=\vparl{F}{\phi}=-2\gamma_1\nabla^2\phi+2\gamma_2\phi+2qq', \quad q'=\vparl{q}{\phi}.
\eena

We rewrite the Allen-Cahn equation given in (\ref{BAC}) using the new variables as follows
\bena
\parl{\phi}{t}= -M \mu,\quad \quad
\parl{q}{t}= q' \phi_t.
\eena
The initial condition of $q$ must be calculated from that of $\phi$.
We denote
\bena
\hs{(\cdot)}=\frac{{(\cdot)}^{n+1}+{(\cdot)}^{n}}{2},\quad
\ohs{(\cdot)}=\frac{3{(\cdot)}^{n+1}-{(\cdot)}^{n}}{2}.
\eena
A second order in time numerical algorithm is given below.

\begin{sch}
 Given initial conditions $\phi^0$ and $q^0$ (calculated from $\phi^0$), we compute $\phi^1, q^1$ by a first order scheme. Having computed $\phi^{n-1},q^{n-1}$, and $\phi^n,q^n$, we compute $\phi^{n+1},q^{n+1}$ as follows.
\bena
\phi^{n+1}-\phi^n= -\Delta t \ohs{M}\hs \mu,\\ \label{EQAC}
\hs \mu=-2\gamma_1 \hs{\nabla ^2\phi}+2 \gamma_2 \hs{\phi} +2\hs{q} \ohs{q'},\\
q^{n+1}-q^n= \ohs{q'}(\phi^{n+1}-\phi^n).
\eena
\end{sch}
We define the discrete energy as follows
\bena
F^n=\int_{\Omega} [\gamma_1(\nabla\phi^{n})^2+\gamma_2{(\phi^{n})}^2+({q}^{n})^2-C_0]\mathrm{d{\bf r}}.
\eena
Then,  the  energy dissipation rate preserving property is guaranteed by the following theorem.
\begin{thm}
The semi-discrete  system obeys the following energy dissipation law
\bena
F^{n+1}-F^n=-\Delta t\int_{\Omega} [\hs{\mu}\ohs{M}\hs{\mu}]\mathrm{d{\bf r}}.\label{energy-dissipation-dis}
\eena
So, the scheme is unconditional stable.
\end{thm}
\noindent {\bf Proof}: the proof is a limiting case of that  of
 theorem 4.3. We thus omit it.

The numerical implementation is done as follows
\bena
(\bI+\Delta t \ohs M [-\gamma_1\nabla^2+\gamma_2 \bI+{(\ohs{q'})}^2])\phi^{n+1}=b^n,\\
b^n=(\phi^{n}-\Delta t \ohs M [-\gamma_1\nabla^2\phi^{n}+\gamma_2\phi^{n}+2q^n\ohs{q'}-{(\ohs{q'})}^2\phi^{n}]),\\
q^{n+1}=q^n+\ohs{q'}(\phi^{n+1}-\phi^n).
\eena
I.e., the  equation of $\phi^{n+1}$ decouples from the equation of $q^{n+1}$ so that $\phi^{n+1}$ is solved independently, and then $q^{n+1}$ is updated. So, the scheme is a sequentially decoupled scheme.

\subsubsection{Numerical method for the Allen-Cahn model using SAV}

\noindent \indent The scalar auxiliary variable (SAV) method provides another way to arrive at linear  numerical schemes. This in fact is the EQ method in another form.
We rewrite the energy functional as follows
\bena
F=\int_{\Omega} [\gamma_1(\nabla\phi)^2+\gamma_2\phi^2]\mathrm{d{\bf r}}+\int_{\Omega} [f(\phi)-\gamma_2\phi^2]\mathrm{d{\bf r}}.
\eena
We define $E_1(\phi)=\int_{\Omega} [f(\phi)-\gamma_2\phi^2]\mathrm{d{\bf r}}$ and choose a constant $C_0$ such that $E_1(\phi)\geq -C_0$. Setting $U=\vparl{E_1}{\phi}$, and introducing $r=\sqrt{E_1+C_0}$ as the scalar auxiliary variable, we arrive at a reformulated Allen-Cahn equation
\bena
\parl{\phi}{t}= -M \mu, \quad \mu=-2\gamma_1 \nabla^2 \phi+2 \gamma_2 \phi +r g, \\
\parl{r}{t}= \int_\Omega \frac{g}{2} \phi_t \mathrm{d{\bf r}}, \quad g=2\vparl{r}{\phi}=\frac{U}{\sqrt{E_1+C_0}},
\eena
where the free energy is given by
\ben
F=\int_{\Omega} [\gamma_1(\nabla\phi)^2+\gamma_2\phi^2]\mathrm{d{\bf r}}+r^2-C_0.
\een

Using the newly formulated  equation system, we design a new scheme as follows.
\begin{sch} Given initial conditions $\phi^0$ and $r^0$, we compute $\phi^1$ and $r^1$ by a first order scheme. Having computed $\phi^{n-1}$, $r^{n-1}$, $\phi^n$ and $r^n$, we compute $\phi^{n+1}$ as follows.
\bena
\phi^{n+1}-\phi^n= -\Delta t \overline{M}^{n+1/2} \hs{\mu},\\ \label{SAV_AC}
\hs \mu=\hs{(-2\gamma_1 \nabla ^2\phi+2\gamma_2\phi)}+\hs r \ohs g,\\
r^{n+1}-r^n= \int_{\Omega} \frac{\ohs g}{2}(\phi^{n+1}-\phi^n)\mathrm{d{\bf r}},\\
\eena
where $\ohs g=\ohs{(\frac{ U( \phi)}{\sqrt {E_1 ( \phi )+C_0}})}$.
\end{sch}
We define the discrete energy as follows
\bena
F^n=\int_{\Omega} [\gamma_1(\nabla\phi^{n})^2+\gamma_2(\phi^n)^2]\mathrm{d{\bf r}}+({r}^{n})^2-C_0.
\eena
Then, energy stability follows from the  following theorem.
\begin{thm}
Scheme (\ref{SAV_AC}) obeys the following energy dissipation law
\bena
F^{n+1}-F^n=-\Delta t\int_{\Omega}[\hs \mu\ohs{M}\hs \mu]\mathrm{d{\bf r}}.
\eena
So, it is unconditionally stable.
\end{thm}
\noindent {\bf Proof:}The proof is similar to that of theorem 4.3 and is thus omitted.

We next discuss how to practically implement the scheme.
The SAV scheme at the nth level can be written into the following form
\bena
A\phi^{n+1}+(c,\phi^{n+1})d=b^n, \label{Example}
\eena
where $A=\bI-\Delta t \ohs{M} (\Delta \gamma_1 -\gamma_2 \bI)$, $c=\ohs g$, $d=\frac{\Delta t \ohs M}{4}\ohs g$ and $b^n=\phi^n-\Delta t \ohs M (-\Delta \gamma_1\phi^n+\gamma_2 \phi^n+r^n c)+d(c,\phi^n)$.
This system can be solved efficiently using the following technique. We multiply the inverse of $A$ and then take the inner product of $c$ to obtain
\bena
(c,\phi^{n+1})+(c,\phi^{n+1}) (c, A^{-1}d^n)=(c,A^{-1} b^n).
\eena
We solve this in the following steps,
\bena
(c,\phi^{n+1})=\frac{(c,A^{-1} b^n)}{1+(c, A^{-1}d^n)},\\ \label{Example1}
\phi^{n+1}=-(c,\phi^{n+1}) A^{-1}d+A^{-1}b^n,\\
r^{n+1}-r^n= \int_{\Omega} \frac{\ohs g}{2}(\phi^{n+1}-\phi^n)\mathrm{d{\bf r}}.
\eena
So in each time step, we only need to solve two Elliptic equations as follows:
\ben
A[x,y]=[d^n,b^n].
\een
This can be done very efficiently.

\subsubsection{Numerical method for the Allen-Cahn model with a penalizing potential  using EQ}

\noindent \indent In the  Allen-Cahn model with a penalizing potential, we reformulate the free energy density by introducing two intermediate variables
\bena
q=\sqrt{f(\phi)-\gamma_2\phi^2+C_0},\quad
\zeta=\sqrt{\eta}(\int_\Omega \phi(t) \mathrm{d{\bf r}}-V_0).
\eena
Then, the free energy is recast into
\bena
F=\int_{\Omega} [\gamma_1(\nabla\phi)^2+\gamma_2{\phi}^2+q^2-C_0]\mathrm{d{\bf r}}+\frac{\zeta^2}{2}.
\eena

We rewrite the nonlocal Allen-Cahn equation as follows:
\bena
\parl{\phi}{t}= -M \tilde \mu,\\
\parl{\zeta}{t}= \sqrt{\eta}\int_\Omega \parl\phi{t} \mathrm{d{\bf r}},\\
\parl{q}{t}= q' \phi_t,
\eena
where
\bena
\tilde \mu=-2\gamma_1 \nabla^2 \phi+2\gamma_2 \phi+2qq'+\zeta \zeta',\quad
q'=\parl{q}{\phi}, \quad \zeta'=\parl{\zeta}{\phi}=\sqrt{\eta}.
\eena
We now discretize it using the linear Crank-Nicolson method in time to arrive at a new scheme as follows.

\begin{sch} Given initial conditions $\phi^0,q^0$, we compute $\phi^1, q^1$ by a first order scheme. Having computed $\phi^{n-1},q^{n-1}$, and $\phi^n,q^n$, we compute $\phi^{n+1},q^{n+1}$ as follows.
\bena
\phi^{n+1}-\phi^n= -\Delta t\ohs{M}\hs {\tilde \mu},\label{EQPAC}\\
q^{n+1}-q^n= \ohs{q'}(\phi^{n+1}-\phi^n),\\
\zeta^{n+1}-\zeta^{n}= \sqrt{\eta}\int_\Omega (\phi^{n+1}-\phi^n) \mathrm{d{\bf r}}.\\
\eena
where
\bena
\hs {\tilde \mu}=-2\gamma_1 \hs{\nabla ^2\phi}+2 \gamma_2 \hs{\phi} +2\hs{q} \ohs{q'}+\sqrt{\eta}\hs{\zeta}, \\
\eena
\end{sch}
We define the discrete energy as follows
\bena
F^n=\int_{\Omega} [\gamma_1(\nabla\phi^{n})^2+\gamma_2(\phi^{n})^2+({q}^{n})^2-C_0]\mathrm{d{\bf r}}+\frac{(\zeta^n)^2}{2}.
\eena
\begin{thm}
Scheme (\ref{EQPAC}) obeys the following energy dissipation law
\ben
F^{n+1}-F^n=-\Delta t\int_{\Omega} [{\tilde \mu}^{n+1/2} \overline{M}^{n+1/2}{\tilde \mu}^{n+1/2}] \mathrm{d{\bf r}}.
\een
Hence, it is unconditionally stable.
\end{thm}
\noindent {\bf Proof:}
Taking the $L^2$ inner product of $\frac{\phi^{n+1}-\phi^n}{\Delta t}$ with $-\hs {\tilde \mu}$, we obtain
\bena
-(\frac{\phi^{n+1}-\phi^n}{\Delta t},\hs {\tilde \mu})
=(\ohs{M}[\hs{\mu}+\sqrt \eta \hs{\zeta}],\hs{\mu}+\sqrt \eta \hs{\zeta})\\
=\lVert\sqrt{\ohs{M}}(\hs{\mu}+\sqrt \eta \hs{\zeta})\rVert^2,\\
\eena
Taking the $L^2$ inner product of $\hs {\tilde \mu}$ with $ \frac{\phi^{n+1}-\phi^n}{\Delta t} $, we obtain
\bena
((\hs{-2\gamma_1 \nabla ^2\phi+2\gamma_2\phi)}+2\hs{q} \ohs{q'}+\sqrt \eta\hs{\zeta},\frac{\phi^{n+1}-\phi^n}{\Delta t})\\
%=((\hs{-\gamma_1 \nabla ^2\phi)},\frac{\phi^{n+1}-\phi^n}{\Delta t})+\frac{\gamma_2}{ \Delta t}(\lVert\phi^{n+1}\rVert^2-\lVert\phi^n\rVert^2)+(2\hs{q} \ohs{q'},\frac{\phi^{n+1}-\phi^n}{\Delta t})+(\sqrt \eta\hs{\zeta},\frac{\phi^{n+1}-\phi^n}{\Delta t})\\
=\frac{\gamma_1}{\Delta t}(\lVert\nabla \phi^{n+1}\rVert^2-\lVert\nabla \phi^n\rVert^2)+\frac{\gamma_2}{ \Delta t}(\lVert\phi^{n+1}\rVert^2-\lVert\phi^n\rVert^2)+(2\hs{q} \ohs{q'},\frac{\phi^{n+1}-\phi^n}{\Delta t})+(\sqrt \eta\hs{\zeta},\frac{\phi^{n+1}-\phi^n}{\Delta t}).
\eena
Taking the $L^2$ inner product of $q^{n+1}-q^n$ with $ \frac{q^{n+1}+q^n}{\Delta t} $, we obtain
\bena
\frac{1}{\Delta t}(\lVert q^{n+1}\rVert^2-\lVert q^{n}\rVert^2)=\frac{1}{\Delta t}(\ohs{q'}(\phi^{n+1}-\phi^n),q^{n+1}+q^{n}).
\eena
Taking the $L^2$ inner product of $\zeta^{n+1}-\zeta^n$ with $ \frac{\zeta^{n+1}+\zeta^n}{\Delta t} $, we obtain
\bena
\frac{1}{\Delta t}(\lVert \zeta^{n+1}\rVert^2-\lVert \zeta^{n}\rVert^2)=\frac{1}{\Delta t}(\sqrt{\eta}\int_\Omega (\phi^{n+1}-\phi^n) \mathrm{d{\bf r}},\zeta^{n+1}+\zeta^{n}).
\eena
Combining the above equations, we have
\bena
\frac{\gamma_1}{\Delta t}(\lVert\nabla \phi^{n+1}\rVert^2-\lVert\nabla \phi^n\rVert^2)+\frac{\gamma_2}{ \Delta t}(\lVert\phi^{n+1}\rVert^2-\lVert\phi^n\rVert^2)+\frac{1}{\Delta t}(\lVert q^{n+1}\rVert^2-\lVert q^{n}\rVert^2)+\frac{1}{2 \Delta t}(\lVert \zeta^{n+1}\rVert^2-\lVert \zeta^{n}\rVert^2)\\
=-\lVert\sqrt{\ohs{M}}(\hs{\mu}+\sqrt \eta \hs \zeta)\rVert^2.
\eena
This leads to the energy stability equality.

The new scheme can be recast into
\ben
A\phi^{n+1}+(\phi^{n+1},c)d=b^n,\label{Sherman-M}
\een
where
\bena
A=\bI+\Delta t \ohs{M}(-\gamma_1 \nabla^2+\gamma_2 \bI+(\ohs{q'})^2), \\
c=1, d=\frac{\Delta t \ohs{M} \eta}{2},\\
\eena
and $b^n=\phi^{n}-\Delta t \ohs M [-\gamma_1\nabla^2\phi^{n}+\gamma_2\phi^{n}+2q^n\ohs{q'}-{(\ohs{q'})}^2\phi^{n}]-\Delta t \ohs M\sqrt\eta\zeta^n+\frac{\Delta t \ohs M\eta}{2}\int_V\phi^{n}\mathrm{d{\bf r}}$.\\
This can be solved  using the  Sherman-Morrison formula to calculate it to avoid calculating the function with a full rank coefficient matrix (See Appendix).
Or equivalently, we can solve it effectively in the following steps:
\bena
A[x,y]=[d,b^n],\\
(\phi^{n+1},c)=\frac{(y,c)}{1+ (x,c)},\\
\phi^{n+1}=-(\phi^{n+1},c) x+y,\\ \label{Sherman-M1}
q^{n+1}=q^n+\ohs{q'}(\phi^{n+1}-\phi^n),\\
\zeta^{n+1}=\zeta^n+\sqrt{\eta}((\phi^{n+1},1)-(\phi^n,1)).
\eena

\subsubsection{Numerical method for the Allen-Cahn model with a penalizing potential  using SAV}

\noindent \indent We can solve the Allen-Cahn model with a penalizing potential utilizing the SAV approach. We rewrite the energy functional as follows
\bena
F=\int_{\Omega} [\gamma_1(\nabla\phi)^2+\gamma_2\phi^2]\mathrm{d{\bf r}}+\int_{\Omega} [f(\phi)-\gamma_2\phi^2]\mathrm{d{\bf r}}+\frac{\eta}{2}(\int_\Omega \phi(t) \mathrm{d{\bf r}}-V_0)^2.
\eena
We define $E_1(\phi)=\int_{\Omega} [f(\phi)-\gamma_2\phi^2]\mathrm{d{\bf r}}$ and choose a constant $C_0$ such that $E_1(\phi)\geq -C_0$. Setting $U=\vparl{E_1}{\phi}$ and introducing $r=\sqrt{E_1+C_0}$ as the scalar auxiliary variable,   we arrive at a new scheme as follows.

\begin{sch}
 Given initial conditions $\phi^0$ and $r^0$, we compute $\phi^1$ and $r^1$ by a first order scheme. Having computed $\phi^{n-1}$, $r^{n-1}$, $\phi^n$ and $r^n$, we compute $\phi^{n+1}$ as follows.
\bena
\phi^{n+1}-\phi^n= -\Delta t \overline{M}^{n+1/2} \hs{\mu},\\ \label{SAV_PP}
\hs \mu=\hs{(-2\gamma_1 \nabla ^2\phi+2\gamma_2\phi)}+\hs r \ohs g+\sqrt{\eta}\hs{\zeta},\\
\zeta^{n+1}-\zeta^{n}= \sqrt{\eta}\int_\Omega (\phi^{n+1}-\phi^n) \mathrm{d{\bf r}},\\
r^{n+1}-r^n= \int_{\Omega} \frac{\ohs g}{2}(\phi^{n+1}-\phi^n)\mathrm{d{\bf r}},\\
\eena
where $\ohs g=\ohs{(\frac{ U( \phi)}{\sqrt {E_1 ( \phi )+C_0}})}$.
\end{sch}
We define the discrete energy as follows
\bena
F^n=\int_{\Omega}[\gamma_1(\nabla\phi^{n})^2+\gamma_2(\phi^{n})^2]\mathrm{d{\bf r}}+\frac{(\zeta^n)^2}{2}+({r}^{n})^2-C_0.
\eena
Then, the following theorem guarantees that the scheme is unconditionally, energy stable. .
\begin{thm}
Scheme (\ref{SAV_PP}) obeys the following energy dissipation law
\bena
F^{n+1}-F^n=-\Delta t\int_{\Omega}[\hs \mu\ohs{M}\hs \mu]\mathrm{d{\bf r}}.
\eena
So, it is unconditionally stable.
\end{thm}
\noindent {\bf Proof:} The proof is similar to that of theorem 4.3 and is thus omitted.

The scheme can be recast into
\ben
A\phi^{n+1}+(\phi^{n+1}, c_1)d_1+(\phi^{n+1}, c_2)d_2=b^n,
\een
where
\bena
A=\bI+\Delta t \ohs M(-\gamma_1 \Delta +\gamma_2 \bI),c_1=\ohs g, d_1=\frac{\Delta t \ohs M}{4}\ohs g, d_2=\frac{\Delta t \ohs M\eta}{2}, c_2=1, \\
b^n=\phi^n-\Delta t\ohs M(-\gamma_1 \phi^n +\gamma_2 \phi^n+r^n\ohs g+\sqrt{\eta}\zeta^n)+(\phi^n,c_1)d_1+(\phi^n,c_2)d_2.
\eena
So we have
\bena
(\phi^{n+1},c_1)+ (\phi^{n+1},c_1) (A^{-1}d_1, c_1)+ (\phi^{n+1},c_2) (A^{-1}d_2, c_1)=(A^{-1} b^n, c_1),\\
(\phi^{n+1},c_2)+ (\phi^{n+1},c_1) (A^{-1}d_1, c_2)+ (\phi^{n+1},c_2) (A^{-1}d_2, c_2)=(A^{-1} b^n, c_2).
\eena
We solve for $(\phi^{n+1},c_1)$ and $(\phi^{n+1},c_2)$ from the above equation after we obtain $(x,y,z)$ from
\ben
A[x,y,z]=[d_1, d_2,b^n].
\een
Finally,
\bena
\phi^{n+1}=z-[ (\phi^{n+1}, c_1) x+(\phi^{n+1}, c_2) y],\\
r^{n+1}=r^n+(\phi^{n+1}-\phi^n,\frac{\ohs g}{2}),\\
\zeta^{n+1}=\zeta^n+\sqrt{\eta}((\phi^{n+1},1)-(\phi^{n},1)).
\eena

\subsubsection{Numerical method for the Allen-Cahn model with a Lagrange multiplier using EQ}

\noindent \indent We use
$
q=\sqrt{f(\phi)-\gamma_2\phi^2+C_0}
$
to recast the free energy into
\bena
F=\int_{\Omega} [\gamma_1(\nabla\phi)^2+\gamma_2{\phi}^2+q^2-C_0]\mathrm{d{\bf r}}.
%-L(\int_{\Omega}h(\phi(t))\mathrm{d{\bf r}}-\int_{\Omega}h(\phi(0)){\bf d r}).
\eena
We rewrite (\ref{BAC}) as
\bena
\parl{\phi}{t}= -M \tilde \mu,\quad
\parl{q}{t}= q' \phi_t,
\eena
where
\bena
\tilde u=u-h'(\phi)L,\quad
L=\frac{\int_{\Omega}{h'(\phi)}M \mu \mathrm{d{\bf r}}}{\int_{\Omega}{h'(\phi)}M{h'(\phi)}\mathrm{d{\bf r}}},\\
\mu=-2\gamma_1\nabla ^2\phi+2 \gamma_2 \phi +2q q',\quad
q'=\parl{q}{\phi}.
\eena

Discretizing it using the linear Crank-Nicolson method in time, we obtain the following scheme.
\begin{sch}
 Given initial conditions $\phi^0,q^0$, we compute $\phi^1, q^1$ by a first order scheme. Having computed $\phi^{n-1},q^{n-1}$, and $\phi^n,q^n$, we compute $\phi^{n+1},q^{n+1}$ as follows.
\bena
\phi^{n+1}-\phi^n= -\Delta t \ohs M \hs {\tilde \mu},\label{EQLAC} \\
q^{n+1}-q^n= \ohs{q'}(\phi^{n+1}-\phi^n).\label {model24}\\
\eena
where \bena
\hs {\tilde \mu}=\hs \mu-\ohs {h'(\phi)}\hs L,\\
\hs \mu=(\hs{-2\gamma_1 \nabla ^2\phi+2\gamma_2\phi)} +2\hs{q} \ohs{q'},\\
\hs L=\frac{\int_{\Omega}\ohs {h'(\phi)} \ohs M \hs \mu \mathrm{d{\bf r}}}{\int_{\Omega}\ohs {h'(\phi)} \ohs M \ohs {h'(\phi)} \mathrm{d{\bf r}}}.\\
\eena
\end{sch}
Note that $\hs L \neq \frac{L^n+L^{n+1}}{2}$.
\begin{thm}
The scheme preserves the volume, namely,
\bena
\int_{\Omega} h(\phi^{n+1}) \mathrm{d{\bf r}}=\int_{\Omega} h(\phi^n) \mathrm{d{\bf r}}.
\eena
\end{thm}
\noindent {\bf Proof:} Substituting  $\hs L$ into the   equation below, we have
\bena
\int_\Omega \frac{h(\phi^{n+1})-h(\phi^{n})}{\Delta t}\mathrm{d{\bf r}}=\int_\Omega \ohs {h'(\phi)}\frac{\phi^{n+1}-\phi^{n}}{\Delta t}\mathrm{d{\bf r}}\\
=\int_\Omega -\ohs {h'(\phi)}\ohs M(\hs \mu-\ohs {h'(\phi)}\hs L)\mathrm{d{\bf r}}=0.
\eena
It implies that it preserves the volume.

We define the discrete energy under the volume constraint as follows
\bena
F^n=\int_{\Omega} [\gamma_1(\nabla\phi^{n})^2+\gamma_2{(\phi^{n})}^2+({q}^{n})^2-C_0]\mathrm{d{\bf r}}.
\eena
\begin{thm}
The semi-discrete  scheme obeys the following energy dissipation law
\bena
F^{n+1}-F^n=-\Delta t\int_{\Omega}[\hs {\tilde \mu}\ohs M\hs {\tilde \mu}]\mathrm{d{\bf r}}
\eena
So, it is unconditionally stable.
\end{thm}
\noindent {\bf Proof:} The proof is similar to that of theorem 4.3 and is thus omitted.

We can  solve the linear system effectively in the following steps:
\bena
A[x,y,z]=[c, d, b^n],\\
(\phi^{n+1},c)=\frac{(z,c)}{1+(y, c)},\\
\phi^{n+1}=z-(\phi^{n+1}, c) x,\\
q^{n+1}=q^n+\ohs{q'}(\phi^{n+1}-\phi^n),
\eena
where
\bena
A=\bI+\Delta t \ohs M(-\gamma_1 \Delta +\gamma_2 \bI+(\ohs{q'})^2),\\
c=\ohs{h'}\ohs{M}(-\gamma_1 \Delta +\gamma_2\bI+(\ohs {q'})^2),\\
d=-\frac{\Delta t \ohs{M}\ohs{h'}}{(\ohs{h'},\ohs M \ohs {h'})},\\
b^n=\phi^n-\Delta t \ohs M (-\gamma_1\Delta\phi^n+\gamma_2\phi^n+2q^n\ohs{q'}
-(\ohs {q'})^2\phi^n)\\-(\ohs{h'}\ohs M,-\gamma_1\Delta \phi^n+\gamma_2\phi^n+2q^n\ohs{q'}-{(\ohs{q'})}^2\phi^n)d.
\eena

\subsubsection{Numerical method for the Allen-Cahn model with a Lagrange multiplier using SAV}

\noindent \indent We rewrite the energy functional as follows
\bena
F=\int_{\Omega} [\gamma_1(\nabla\phi)^2+\gamma_2\phi^2]\mathrm{d{\bf r}}+\int_{\Omega} [f(\phi)-\gamma_2\phi^2]\mathrm{d{\bf r}}-L(\int_{\Omega}h(\phi(t))\mathrm{d{\bf r}}-\int_{\Omega}h(\phi(0))\mathrm{d{\bf r})}.
\eena
We define $E_1(\phi)=\int_{\Omega} [f(\phi)-\gamma_2\phi^2]\mathrm{d{\bf r}}$ and choose a constant $C_0$ such that $E_1(\phi)\geq -C_0$. Setting $U=\vparl{E_1}{\phi}$ and introducing $r=\sqrt{E_1+C_0}$ as the scalar auxiliary variable, then we arrive at a new scheme as follows.

\begin{sch} Given initial conditions $\phi^0$ and $r^0$, we  compute $\phi^1$ and $r^1$ by a first order scheme. Having computed $\phi^{n-1}$, $r^{n-1}$, $\phi^n$ and $r^n$, we compute $\phi^{n+1}$ as follows.\\
\bena
\phi^{n+1}-\phi^n= -\Delta t \overline{M}^{n+1/2} \hs{\tilde \mu},\\ \label{SAVLAC}
\hs {\tilde \mu}=\hs \mu-\ohs {h'(\phi)}\hs L,\\
\hs \mu=\hs{(-2\gamma_1 \nabla ^2\phi+2\gamma_2\phi)}+\hs r \ohs g,\\
\hs L=\frac{\int_{\Omega}\ohs {h'(\phi)} \ohs M \hs \mu \mathrm{d{\bf r}}}{\int_{\Omega}\ohs {h'(\phi)} \ohs M \ohs {h'(\phi)} \mathrm{d{\bf r}}},\\
r^{n+1}-r^n= \int_{\Omega} \frac{\ohs g}{2}(\phi^{n+1}-\phi^n)\mathrm{d{\bf r}},\\
\eena
where $\ohs g=\ohs{(\frac{ U( \phi)}{\sqrt {E_1 ( \phi )+C_0}})}$.
\end{sch}
%Remark: $\hs L \neq \frac{L^n+L^{n+1}}{2}$.\\
\begin{thm}
The volume of each phase is conserved, i.e.,
\bena
\int_{\Omega} h(\phi^{n+1}) \mathrm{d{\bf r}}=\int_{\Omega} h(\phi^n) \mathrm{d{\bf r}}.
\eena
\end{thm}
\noindent {\bf Proof:} The proof is similar to that of theorem 4.5 and is thus omitted.

We define the discrete energy under the volume conservation condition as follows
\bena
F^n=\int_{\Omega} [\gamma_1(\nabla\phi^{n})^2+\gamma_2{(\phi^{n})}^2]\mathrm{d{\bf r}}+({r}^{n})^2-C_0.
\eena
\begin{thm}
The semi-discrete  scheme obeys the following energy dissipation law
\bena
F^{n+1}-F^n=-\Delta t\int_{\Omega}[\hs {\tilde \mu}\ohs M\hs {\tilde \mu}]\mathrm{d{\bf r}}.
\eena
So it is unconditional stable.
\end{thm}
\noindent {\bf Proof:} The proof is similar to that of theorem 4.3 and is thus omitted.

This scheme can be recast into
\ben
A\phi^{n+1}+(\phi^{n+1}, c_1)d_1+(\phi^{n+1},c_2)d_2+(c_3,(\phi^{n+1},c_1))d_2=b^n,
\een
So we have
\bena
(\phi^{n+1},c_1)+(\phi^{n+1},c_1)(A^{-1} d_1, c_1)+(\phi^{n+1},c_2)(A^{-1} d_2, c_1)+(c_3,(\phi^{n+1},c_1))(A^{-1}d_2,c_1)=(A^{-1} b^n,c_1),\\
(\phi^{n+1},c_2)+(\phi^{n+1},c_1)(A^{-1} d_1, c_2)+(\phi^{n+1},c_2)(A^{-1} d_2, c_2)+(c_3,(\phi^{n+1},c_1))(A^{-1}d_2,c_2)=(A^{-1} b^n,c_2).\\
\eena
We solve for $(\phi^{n+1},c_1)$ and $(\phi^{n+1},c_2)$ from the above equation after we obtain $(x,y,z)$ from
\bena
A[x,y,z]=[d_1,d_2, b^n],
\eena
where
\bena
A=\bI+\Delta t \ohs M(-\gamma_1\Delta +\gamma_2\bI),\\
c_1=\ohs g,
d_1=\Delta t \ohs M\frac{\ohs g}{4},\\c_2=\ohs {h'}\ohs M(-\gamma_1 \Delta +\gamma_2I),\\
d_2=-\Delta t \ohs M\frac{\ohs{h'}}{(\ohs{h'},\ohs{M}\ohs{h'})},\\
c_3=\ohs {h'}\ohs M \frac{\ohs g}{4},\\
b^n=\phi^n-\Delta t\ohs M(-\gamma_1\Delta \phi^n +\gamma_2\phi^n+r^n\ohs{g}-\ohs{g}\int_\Omega\frac{\ohs{g}}{4}\phi^n\mathrm{d{\bf r}})\\
-(\ohs{h'}\ohs M,-\gamma_1\Delta \phi^n+\gamma_2\phi^n+r^n\ohs{g}-\ohs{g}\int_\Omega\frac{\ohs{g}}{4}\phi^n\mathrm{d{\bf r}})d_2.\\
\eena
Finally, we have
\bena
\phi^{n+1}=z-(\phi^{n+1}, c_1)x-(\phi^{n+1},c_2)y-(c_3,(\phi^{n+1},c_1))y,\\
r^{n+1}=r^n+(\frac{\ohs g}{2},(\phi^{n+1}-\phi^n)).\\
\eena

In order to make a comparison with the volume preserving Cahn-Hilliard equation for binary material systems, we present two energy stable schemes for the Cahn-Hilliard equation next.

\subsubsection{Numerical methods for the Cahn-Hilliard model using EQ}

 \noindent \indent We use $
q=\sqrt{f(\phi)-\gamma_2\phi^2+C_0} $ to
  recast the free energy into
\bena
F=\int_{\Omega} [\gamma_1(\nabla\phi)^2+\gamma_2{\phi}^2+q^2-C_0]\mathrm{d{\bf r}}.
\eena
Then,
 \bena
\mu=\frac{\delta F}{\delta \phi}=-2\gamma_1 \nabla ^2\phi+2\gamma_2\phi +2qq'\label {model23}, \quad q'=\parl{q}{\phi}.
 \eena
We rewrite (\ref{BCH}) as
\bena
\parl{\phi}{t}= \nabla \cdot (M \nabla\mu),\quad
\parl{q}{t}= q' \phi_t.\\
\eena
We discretize it using the linear Crank-Nicolson method in time  to arrive at a second order semi-discrete  scheme.

\begin{sch}
 Given initial conditions $\phi^0,q^0$, we compute $\phi^1, q^1$ by a first order scheme. Having computed $\phi^{n-1},q^{n-1}$, and $\phi^n,q^n$, we compute $\phi^{n+1},q^{n+1}$ as follows.
\bena
\phi^{n+1}-\phi^n= \Delta t \nabla \cdot (\overline{M}^{n+1/2} \nabla \hs  \mu),\label{EQCH}\\
\hs \mu=(\hs{-2\gamma_1 \nabla ^2\phi+2\gamma_2\phi)} +2\hs{q} \ohs{q'},\\
q^{n+1}-q^n= \ohs{q'}(\phi^{n+1}-\phi^n).
\eena
\end{sch}
We define the discrete energy as follows
\bena
F^n=\int_{\Omega} [\gamma_1(\nabla\phi^{n})^2+\gamma_2{(\phi^{n})}^2+({q}^{n})^2-C_0]\mathrm{d{\bf r}}.
\eena
Then, we can obtain the following theorem.
\begin{thm}
The semi-discrete  scheme obeys the following energy dissipation law
\bena
F^{n+1}-F^n=-\Delta t\int_{\Omega} [\nabla\hs{\mu}]\ohs{M}[\nabla\hs{\mu}]\mathrm{d{\bf r}}.\\
\eena
So, it is unconditionally stable.
\end{thm}
\noindent {\bf Proof:} Multiplying the above three equations with $\hs \mu$, $\phi^{n+1}-\phi^n$, $q^{n+1}+q^n$, and taking the sum, we obtain the desired results.

The implementation of the scheme is as follows
\bena
A\phi^{n+1}=b^n,\\
A=\bI-\Delta t \nabla \cdot (\overline{M}^{n+1/2} \nabla [-\gamma_1 \nabla ^2+\gamma_2\bI + (\ohs{q'})^2]),\\
b^n=\phi^n+\Delta t \nabla \cdot (\overline{M}^{n+1/2} \nabla[-\gamma_1 \nabla ^2\phi^n+\gamma_2\phi^n +2{q}^n \ohs{q'}-(\ohs{q'})^2\phi^n]),\\
q^{n+1}= q^n+\ohs{q'}(\phi^{n+1}-\phi^n).
\eena

\subsubsection{Numerical methods for the Cahn-Hilliard model using SAV}

\noindent \indent
We rewrite the energy functional as follows
\bena
F=\int_{\Omega} [\gamma_1(\nabla\phi)^2+\gamma_2\phi^2]\mathrm{d{\bf r}}+\int_{\Omega} [f(\phi)-\gamma_2\phi^2]\mathrm{d{\bf r}}.
\eena
We define $E_1(\phi)=\int_{\Omega} [f(\phi)-\gamma_2\phi^2]\mathrm{d{\bf r}}$ and choose a constant $C_0$ such that $E_1(\phi)\geq -C_0$. Setting $U=\vparl{E_1}{\phi}$, and introducing $r=\sqrt{E_1+C_0}$ as the scalar auxiliary variable,   we arrive at a new scheme as follows.
\begin{sch}
 Given initial conditions $\phi^0$ and $r^0$, we  compute $\phi^1$ and $r^1$ by a first order scheme. Having computed $\phi^{n-1}$, $r^{n-1}$, $\phi^n$ and $r^n$, we compute $\phi^{n+1}$ as follows.
\bena
\phi^{n+1}-\phi^n= \Delta t \nabla \cdot(\overline{M}^{n+1/2} \nabla \hs{\mu}),\\ \label{SAVAC}
\hs \mu=\hs{(-2\gamma_1 \nabla ^2\phi+2\gamma_2\phi)}+\hs r \ohs g,\\
r^{n+1}-r^n= \int_{\Omega} \frac{\ohs g}{2}(\phi^{n+1}-\phi^n)\mathrm{d{\bf r}},
\eena
where $\ohs g=\ohs{(\frac{ U( \phi)}{\sqrt {E_1 ( \phi )+C_0}})}$.
\end{sch}
We define the discrete energy as follows
\bena
F^n=\int_{\Omega} [\gamma_1(\nabla\phi^{n})^2+\gamma_2(\phi^n)^2]\mathrm{d{\bf r}}+({r}^{n})^2-C_0.
\eena

\begin{thm}
Scheme (\ref{SAVAC}) obeys the following energy dissipation law
\bena
F^{n+1}-F^n=-\Delta t\int_{\Omega} [\nabla\hs{\mu}]\ohs{M}[\nabla\hs{\mu}]\mathrm{d{\bf r}}.\\
\eena
So, it is unconditionally stable.
\end{thm}
\noindent {\bf Proof:} Multiplying the above three equations with $\hs \mu$, $\phi^{n+1}-\phi^n$, $r^{n+1}+r^n$, and taking the sum, we  get the desired results \cite{shen2018scalar}.

We can solve it effectively in the following steps
\bena
A[x,y]=[ d, b^n],\\
(\phi^{n+1},c)=\frac{(y,c)}{1+(x,c)},\\
\phi^{n+1}=-(\phi^{n+1},c) x+y,\\
r^{n+1}= r^n+\int_{\Omega} \frac{\ohs g}{2}(\phi^{n+1}-\phi^n)\mathrm{d{\bf r}},
\eena
where
\bena
A=\bI-\Delta t(\nabla \ohs M \nabla+\ohs M \Delta)(-\gamma_1 \Delta +\gamma_2 \bI),\\
c=\frac{\ohs g}{4},
d=-\Delta t(\nabla \ohs M \nabla+\ohs M \Delta){\ohs g},\\
b^n=\phi^n+\Delta t(\nabla \ohs M \nabla+\ohs M \Delta)((-\gamma_1 \Delta +\gamma_2 \bI)\phi^n+r^n{\ohs g}-\ohs{g}\int_\Omega\frac{\ohs{g}}{4}\phi^n\mathrm{d{\bf r}}).
\eena
%For one of the linear schemes, we prove that the linear system resulting from that scheme is uniquely solvable. This result applies to all the linear schemes presented in this section besides the ones for the Cahn-Hilliard equation. The proofs of EQ and SAV scheme for the Cahn-Hilliard model are similar. We summarize the results into a theorem.
%\begin{thm}
%The linear system resulting from the linear schemes admits a unique solution.
%\end{thm}

\subsubsection{Solvability of the linear systems resulting from the schemes}

\noindent \indent We summarize the unique solvability  for all linear systems resulting from the schemes presented in this section   for the Allen-Cahn  and nonlocal Allen-Cahn models into a theorem. A similar theorem for the Cahn-Hilliard model can be found in \cite{chen2018regularized}.
\begin{thm}
The linear system resulting from any scheme for the Allen-Cahn and nonlocal Allen-Cahn models  admits a unique weak solution.
\end{thm}
 The proof of the uniqueness of the solution for the EQ or SAV scheme can be summarized into two classes. The first class is for the penalizing potential model, and the second one is for the   model with a Lagrange multiplier.

Li et. al proved the uniqueness of the weak solution of  the EQ scheme with periodic boundary conditions for the nonlocal Allen-Cahn model with a penalizing potential \cite{li2018unconditionally}. Here we expand it into both EQ and SAV schemes for the Allen-Cahn model with a penalizing potential subject to a Neumann boundary condition. For simplicity, we assume $M$ is a positive constant mobility coefficient. The  schemes for the nonlocal models can be generically written into:
\ben
A\phi^{n+1}+\sum_{i=1}^{N}(\phi^{n+1},c_i)d_i=r^n,
\een
 where $A=\bI+\Delta t \ohs M(-\gamma_1 \Delta +\gamma_2 \bI+(\ohs{q'})^2)$, $r^n$ is the given right hand side term, and  $N=1$ for the EQ scheme, while $A=\bI+\Delta t \ohs M(-\gamma_1\Delta +\gamma_2 \bI)$ and $N=2$ for the SAV scheme.
Both  EQ and SAV schemes guarantee the value of $\frac{d_i}{c_i}$ is a positive constant, where $d_i=\lambda_i c_i$ ($\lambda_i>0$). Note that $c_i=d_i=0$ for the Allen-Cahn schemes.
We now express the linear system as
\ben
\mathscr{A}\phi^{n+1}=b.
\een
 To prove the uniqueness of the solution, we need to prove that the linear spatial operator $\mathscr{A}$ is symmetric positive definite.

\noindent {\bf Proof}:
We define the bilinear form:
\bena
(\mathscr{A}\phi,\psi)=(A\phi,\psi)+\sum_{i=1}^{N}(\phi,c_i)(d_i,\psi)=(\psi,\mathscr{A}\phi).
\eena
Then, $\exists C>0$ such that
\bena
(\mathscr{A}\phi,\phi)=(A\phi,\phi)+\sum_{i=1}^{N}(\phi,c_i)(d_i,\phi)=(A\phi,\phi)+\sum_{i=1}^{N}(\phi,c_i)(\lambda_i c_i,\phi)\geqslant (A\phi,\phi) > C \left\|\phi\right\|_{L^2}^2.
\eena
We can easily show that the operator $\mathscr{A}$ is bounded above by $\bar{C} \|\phi\|^2$.

Now we define $\left\|\phi\right\|_{\mathscr{A}}=\sqrt{(\mathscr{A}\phi,\phi)}$ for any $\phi \in L_{Neu}^2(\Omega)=\left\{\phi \in L^2(\Omega): \left\|\phi\right\|_{\mathscr{A}}<\infty, \frac{\partial \phi}{\partial n}|_{\partial \Omega}=0 \right\}$, where $L_{Neu}^2(\Omega)$ is a subset of $L^2(\Omega)$.  It is obviously that $\left\|\phi\right\|_{\mathscr{A}}$ is a norm for $L_{Neu}^2(\Omega)$ and $L_{Neu}^2(\Omega)$ is a Hilbert subspace. Applying the Lax-Milgram theorem, the uniqueness of the solution of the linear systems in $\Phi$ is established.

 For the nonlocal Allen-Cahn model with a Lagrangian multiplier,  we only present the detailed proof for the EQ scheme since the proof for the SAV schemes is similar.
From scheme (\ref{EQLAC}) we understand that in order to prove the uniqueness of the solution of linear system $A\phi^{n+1}+(c,\phi^{n+1})d=b^n$ with the Neumann boundary,  we only need to prove that
\bena
A\phi^{n+1}+(c,\phi^{n+1})d=0
\eena
has the zero solution only, where
\bena
A=\bI+\Delta t \ohs M(-\gamma_1 \Delta +\gamma_2 \bI+(\ohs{q'})^2), \quad c=\ohs M(-\gamma_1 \Delta +\gamma_2 \bI+(\ohs{q'})^2),\\
 d=-\frac{\Delta t \ohs M}{\int_\Omega \ohs M \mathrm {d \bf r}},\quad
  \quad b^n=\phi^n-\delta t \ohs M [-\gamma_1\Delta \phi^n+\gamma_2\phi^n-{(\ohs{q'})}^2\phi^n+2\ohs{q'}q^n-\\
  \frac{\int_\Omega \ohs{M}(-\gamma_1\Delta \phi^n+\gamma_2\phi^n-{(\ohs{q'})}^2\phi^n+2\ohs{q'}q^n)\mathrm{d \bf r}}{\int_\Omega \ohs M \mathrm{d \bf r}}].
 \eena

Notice that S Boussaid et al \cite{boussaid2015convergence} proved the uniqueness of the solution for
\bena
(A-\bI) \phi+(c,\phi)d=0.
\eena
 So it is easy to arrive at the uniqueness of the solution for equation $A\phi+(c,\phi)d=0$ since $\bI$ is a positive definite operator.

Note that we can only prove the uniqueness of the solution in the case of $h'(\phi)=1$ here. For a general $h(\phi)$, the uniqueness proof still eludes us.

\subsection{Spatial discretization}

\noindent \indent We use  the finite difference method to discretize the semidiscrete schemes presented above in space. The Neumann boundary condition is adopted in the discretization.
 We divide the 2D domain $\Omega=[0,L_x] \times [0, L_y]$ into uniform rectangular meshes with mesh sizes $h_x=L_x/N_x$ and $h_y=L_y/N_y$, where $L_x$, $L_y$ are two positive real numbers and $N_x$, $N_y$ are the number of meshes in the x and y direction, respectively. The sets of the cell center points $C_x$ and $C_y$ in the  uniform partition are defined as follows
\bena
C_x=\left\{x_{i}|i=0, 1, \cdots, N_x \right\},\\
C_y=\left\{y_{j}|j=0, 1, \cdots, N_y \right\},
\eena
where $x_i=(i-\frac{1}{2})h_x$ and $y_j=(j-\frac{1}{2})h_y$. The phase field variable is discretized at the cell center points $C_x\times C_y$.

We  define the east-west-edge-to-center and center-to-east-west-edge difference operators $d_x$ and $D_x$, respectively,
\bena
d_x\phi_{ij}=\frac{\phi_{i+\frac{1}{2},j}-\phi_{i-\frac{1}{2},j}}{h_x}, \qquad D_x\phi_{i+\frac{1}{2,}j}=\frac{\phi_{i+1,j}-\phi_{i-1,j}}{h_x}.
\eena
Similarly, we can get the north-south-edge-to-center and center-to-north-south-edge difference operators $d_y$ and $D_y$, respectively,
\bena
d_y\phi_{ij}=\frac{\phi_{i,j+\frac{1}{2}}-\phi_{i,j-\frac{1}{2}}}{h_y}, \qquad D_y\phi_{i,j+\frac{1}{2}}=\frac{\phi_{i,j+1}-\phi_{i,j-1}}{h_y}.
\eena
The fully discrete Laplacian operator is defined by
\bena
\Delta_h=d_x(D_x\phi)+d_y(D_y\phi).
\eena
The  inner product at cell center is denoted by
\bena
\langle f,g \rangle =h_xh_y\sum_{i,j}f_{i,j}g_{i,j}.\label{DiscreteInnerproduct}
\eena
Replacing the differential operators in the semidiscrete schemes by the discrete operators properly, we obtain the fully discrete schemes. The energy stability  and  volume-preserving property are retained in the fully discrete schemes as well. For more details on the spatial discretization, we  refer readers to the papers \cite{Gong2018Fully,Gong2018Linear}.
%We now also express the linear system of Lagrangian schemes as $\mathscr{A}\phi^{n+1}=b$. To prove the uniqueness of the solution, we also need to prove that the linear spatial operator $\mathscr{A}$ is symmetric positive definite firstly.\\
%\noindent {\bf Proof}:
%It is effortless to obtain
%\bena
%(\mathscr{A}\phi,\psi)=(A\phi,\psi)+(\phi,c_1)(d_1,\psi)=(\psi,\mathscr{A}\phi),
%\eena
%where $c_1=\ohs{h'}M\Delta t (-\gamma_1 \Delta +\gamma_2\bI+(\ohs {q'})^2)=\ohs{h'}(A-\bI)$ and $d_1=-\frac{M\ohs{h'}}{(\ohs{h'},\ohs {h'})}$. Substituting $c_1$,$d_1$ in to it, we can get
%\bena
%(\mathscr{A}\phi,\phi)=(A\phi,\phi)-\frac{(\ohs{h'},\phi) (\ohs{h'},(A-\bI)\phi)}{(\ohs{h'},\ohs {h'})}.
%\eena
%If $h'(\phi)=1$, then
%\bena
%(\mathscr{A}\phi,\phi)=(A\phi,\phi)-\frac{(1,\phi) (1,(A-\bI)\phi)}{(1,1)}.
%\eena
%The magnitude of $\frac{(1,\phi) (1,(A-\bI)\phi)}{(1,1)}$ is $O(\Delta t)$, if we choose $\Delta t$ small enough, we can guarantee $(\mathscr{A}\phi,\phi)>0$. \\
%If $h'(\phi)=\frac{(m+1)(m+2)}{m}\phi^m{(1-\phi)}^m$, then
%\bena
%(\mathscr{A}\phi,\phi)=(A\phi,\phi)-\frac{(\ohs{h'},\phi) (\ohs{h'},(A-\bI)\phi)}{(\ohs{h'},\ohs {h'})}.
%\eena
%If both the width of the interface of the phases $\varepsilon$ and $\Delta t$ are small enough, $\frac{(\ohs{h'},\phi)(\ohs{h'}(A-I),\phi)}{(\ohs{h'},\ohs {h'})}$can be bounded by a small positive constant with a fixed $m$ , so $(\mathscr{A}\phi,\phi)>0$. In this sense, we can prove that $\mathscr{A}$ is symmetric positive definite.\\
%Similar as before, we can prove the uniqueness of the weak solution for Lagrangian model by using the Lax-Milgram theorem.

\section{Numerical Results and Discussions}

\noindent \indent In this section, we conduct mesh refinement tests to validate the accuracy of the proposed schemes and then present some numerical examples to assess the schemes for the nonlocal Allen-Cahn models against those for the Cahn-Hilliard model.  When considering the definition of the volume, we use two different choices of $h(\phi)$. If $h(\phi)=\phi$, we call the model Lagrangian model 1, otherwise we name it  Lagrangian model 2. For convenience, we refer the numerical schemes designed by EQ strategy for the Allen-Cahn model, the  Allen-Cahn model with the penalizing potential, Lagrangian model 1 and 2, and the Cahn-Hilliard model as AC-EQ, AC-P-EQ, AC-L1-EQ, AC-L2-EQ, CH-EQ, respectively. Similarly, we name the numerical schemes obtained using SAV approaches for the models as AC-SAV, AC-P-SAV, AC-L1-SAV, AC-L2-SAV, CH-SAV, respectively.

 In the following, we set $\eta=1\times 10^5$ in the EQ and SAV schemes for the nonlocal Allen-Cahn model with a penalizing potential unless noted otherwise. In addition, we set the constant in the free energy at $C_0=1\times 10^ 5$ in all computations.

\subsection{Mesh refinement}
\noindent \indent
We refine the  mesh systematically to test the accuracy by setting $\varepsilon=0.1$ ($\gamma_1=5\times 10^{-2}$ and $\gamma_2=10$) with the initial condition given by
\bena
\phi(0,x,y)=\frac{1}{2}+\frac{1}{2}\cos(\pi x)\cos(\pi y). \label{Ini1}
\eena
The computational domain is set as $\Omega=[- 1,1]^2$.
We choose the solution obtained at $\Delta t=10^{-4}$ and $ \Delta x =\Delta y= \frac {1}{256}$ as the "exact solution". In table 1 and 2, we list the $L^2$ errors of the phase variable between the numerical solutions and the "exact solution" at $t=2$ with respect to different time steps.
In table 1 and table 2, we show the convergence rates  match second order accuracy with respect to different time steps, $\Delta t= 5\times 10^{-4}, 10^{-3},5\times 10^{-3},10^{-2},5\times 10^{-2}$, respectively.  We note that the schemes are also second order accurate in space and omit the details of the mesh refinement tests in space.

\vskip 10 pt
{
\centerline{{\bf Table 1} Numerical errors and convergence rates  of the EQ schemes in time.}
\vskip 10 pt
\centerline{
\begin {tabular} [!htbp] {|c|c|c|c|c|c|c|c|c|c|c|} %{lccr}
\hline
%\cline{2-4}
Scheme & \multicolumn {2}{|c|} {AC-EQ} & \multicolumn {2}{|c|} {AC-P-EQ} & \multicolumn {2}{|c|} {AC-L1-EQ}& \multicolumn {2}{|c|} {AC-L2-EQ}& \multicolumn {2}{|c|} {CH-EQ}\\
\cline{0-10}
\hline
$\delta t$ & $L^2$ error & order& $L^2$ error & order & $L^2$ error & order& $L^2$ error & order & $L^2$ error & order\\
\cline{0-10}
5.00E-02&2.16E-03&-&7.68E-07&-&7.80E-07 &- &4.89E-07&- &1.22E-06&-\\
1.00E-02&3.63E-04&1.11&7.14E-08&1.48&7.17E-08&1.48 &4.46E-08&1.49&4.88E-08&2.00\\
5.00E-03&9.71E-05&1.90&2.05E-08&1.80&2.05E-08&1.81 &1.29E-08&1.79&1.22E-08&2.00\\
1.00E-03&3.93E-06&1.99&9.14E-10&1.93&8.86E-10&1.95 &5.78E-10&1.93&4.83E-10&2.01\\
5.00E-04&9.52E-07&2.05&2.25E-10&2.02&1.99E-10&2.15 &1.37E-10&2.08&1.18E-10&2.03\\
\hline
\end {tabular}
}}
\vskip 10 pt
\centerline{{\bf Table 2} Numerical errors  and convergence rates  of the SAV schemes in time.}
\vskip 10 pt
\centerline{
\begin {tabular} [!htbp] {|c|c|c|c|c|c|c|c|c|c|c|} %{lccr}
\hline
%\cline{2-4}
Scheme & \multicolumn {2}{|c|} {AC-SAV} & \multicolumn {2}{|c|} {AC-P-SAV} & \multicolumn {2}{|c|} {AC-L1-SAV}& \multicolumn {2}{|c|} {AC-L2-SAV}& \multicolumn {2}{|c|} {CH-SAV}\\
\cline{0-10}
\hline
$\delta t$ & $L^2$ error & order& $L^2$ error & order & $L^2$ error & order& $L^2$ error & order & $L^2$ error & order\\
\cline{0-10}
5.00E-02&9.77E-04&-&1.41E-07&-&1.40E-07&- &2.86E-07&- &1.22E-06&-\\
1.00E-02&1.26E-04&1.27&9.80E-09&1.66&9.73E-09&1.66&1.83E-08&1.71&4.18E-08&2.00\\
5.00E-03&3.08E-05&2.03&2.72E-09&1.85&2.70E-09&1.85&5.02E-09&1.87&1.22E-08&2.00\\
1.00E-03&1.12E-06&2.06&1.18E-10&1.95&1.17E-10&1.95&2.14E-10&1.96&4.83E-10&2.01\\
5.00E-04&2.66E-07&2.07&2.90E-11&2.02&2.89E-11&2.02&5.10E-11&2.07&1.17E-10&2.05\\
\hline				
\end {tabular}
}
\vskip 10 pt
\centerline{{\bf Table 3} Computational efficiency of the models in 1000 time steps with $M=1\times 10^{-4}$.}
\vskip 10 pt
\centerline{
\begin {tabular} [!htbp] {|c|c|c|c|c|c|c|c|c|c|c|} %{lccr}
\hline
%\cline{2-4}
Scheme & \multicolumn {2}{|c|} {AC-EQ/SAV} & \multicolumn {2}{|c|} {AC-P-EQ/SAV} & \multicolumn {2}{|c|} {AC-L1-EQ/SAV}& \multicolumn {2}{|c|} {AC-L2-EQ/SAV}& \multicolumn {2}{|c|} {CH-EQ/SAV}\\
\cline{0-10}
\hline
%$\delta t$ & $L^2$ error & order& $L^2$ error & order & $L^2$ error & order& $L^2$ error & order & $L^2$ error & order\\
\cline{0-10}
Time for EQ (s)&\multicolumn {2}{|c|} {27}&\multicolumn {2}{|c|} {42}&\multicolumn {2}{|c|} {57}&\multicolumn {2}{|c|} {58}&\multicolumn {2}{|c|} {37}\\
Time for SAV (s)&\multicolumn {2}{|c|} {38}&\multicolumn {2}{|c|} {41}&\multicolumn {2}{|c|} {58}&\multicolumn {2}{|c|} {50}&\multicolumn {2}{|c|} {53}\\
\hline				
\end {tabular}
}
\vskip 10 pt
\centerline{{\bf Table 4} Computational efficiency of the models in 1000 time steps with $M=1\times 10^{-2}$.}
\vskip 10 pt
\centerline{
\begin {tabular} [!htbp] {|c|c|c|c|c|c|c|c|c|c|c|} %{lccr}
\hline
%\cline{2-4}
Scheme & \multicolumn {2}{|c|} {AC-EQ/SAV} & \multicolumn {2}{|c|} {AC-P-EQ/SAV} & \multicolumn {2}{|c|} {AC-L1-EQ/SAV}& \multicolumn {2}{|c|} {AC-L2-EQ/SAV}& \multicolumn {2}{|c|} {CH-EQ/SAV}\\
\cline{0-10}
\hline
%$\delta t$ & $L^2$ error & order& $L^2$ error & order & $L^2$ error & order& $L^2$ error & order & $L^2$ error & order\\
\cline{0-10}
Time for EQ (s)&\multicolumn {2}{|c|} {27}&\multicolumn {2}{|c|} {41}&\multicolumn {2}{|c|} {60}&\multicolumn {2}{|c|} {60}&\multicolumn {2}{|c|} {48}\\
Time for SAV (s)&\multicolumn {2}{|c|} {37}&\multicolumn {2}{|c|} {41}&\multicolumn {2}{|c|} {57}&\multicolumn {2}{|c|} {53}&\multicolumn {2}{|c|} {60}\\
\hline				
\end {tabular}
}
\vskip 12 pt

 We also test the performance of the schemes  when implemented with respect to the initial condition in terms of the computational efficiency. In table 3, we list the performance of the ten schemes computed in 1000 time steps. The performance depends on the mobility, the other model parameters as well as the initial condition. In the tests conducted with the chosen parameter values, we observe that the best performance, among the volume preserving schemes, is given by the  scheme for the Allen-Cahn model with a penalizing potential at $M=1\times 10^{-4}$ for the initial condition given in equation(\ref{Ini1}). The performance of the schemes computed using $M=1\times10^{-2}$ is given  in table 4. The scheme for the Allen-Cahn equation with a penalizing potential  once again performs better. The pre-factor of the penalizing term for the Allen-Cahn model with penalizing potential is set at $\eta=1\times 10^{5}$ in the computations. In the case where the mobility is larger, the performance in the SAV schemes improve, in which some even surpass the corresponding EQ schemes. Overall, the performance of the schemes between the EQ and the SAV class is comparable.

% We next assess the numerical schemes on a couple of benchmark examples.

\subsection{Assessment of the numerical schemes}

\noindent \indent In this section, we will  assess the numerical schemes derived by using EQ and SAV methods on two benchmark problems.
Firstly, we study merging of  two drops  using the numerical schemes to examine the volume preserving property of the nonlocal models as well as energy dissipation.

We put two drops, next to each other, in  the computational domain.  The drops and the ambient are represented  by $\phi=1$ and $\phi=0$, respectively.  The parameter values of the models are chosen as $\gamma_1=5\times 10^{-3}$, $\gamma_2=100$. The initial condition is given by
\begin {equation}
\left\{
\begin{array}{lr}
1,& r_1\leq0.2-\delta \quad \text{or}\quad r2\leq0.2-\delta,\\
tanh(\frac{0.2-r_1}{\delta}),& 0.2-\delta<r_1\leq0.2,\\
tanh(\frac{0.2-r_2}{\delta}),& 0.2-\delta<r_2\leq0.2,\\
0,& \text {other},
\end{array}
\right.
\end{equation}
where $r_1=\sqrt{(x-0.3)^2+(y-0.5)^2}$, $r_2=\sqrt{(x-0.7)^2+(y-0.5)^2}$ and $\delta=0.01.$

\begin{figure*}
\centering
\subfigure[]{
\begin{minipage}[b]{0.49\linewidth}
\includegraphics[width=1\linewidth]{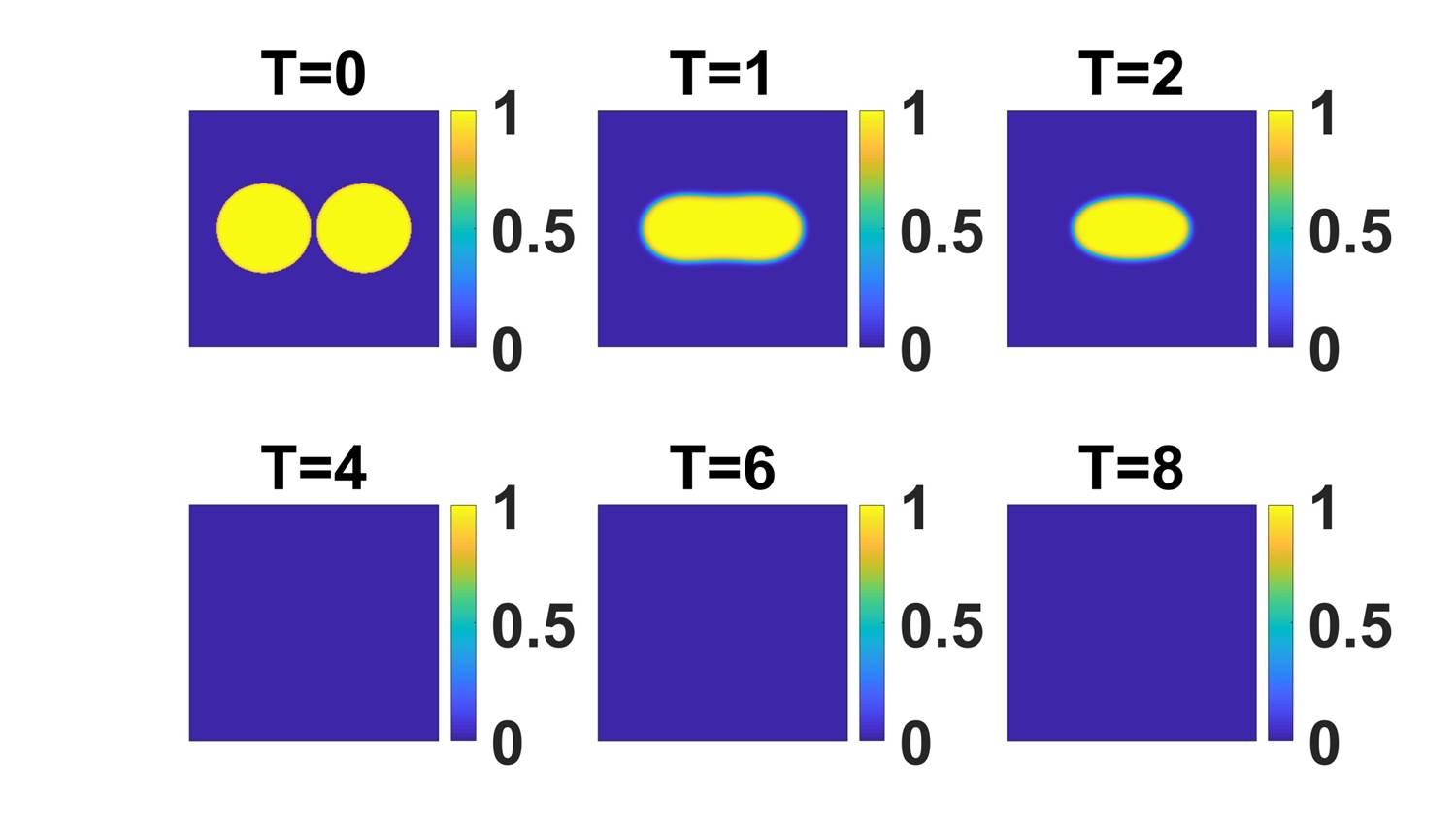}
\end{minipage}}
\subfigure[]{
\begin{minipage}[b]{0.49\linewidth}
\includegraphics[width=1\linewidth]{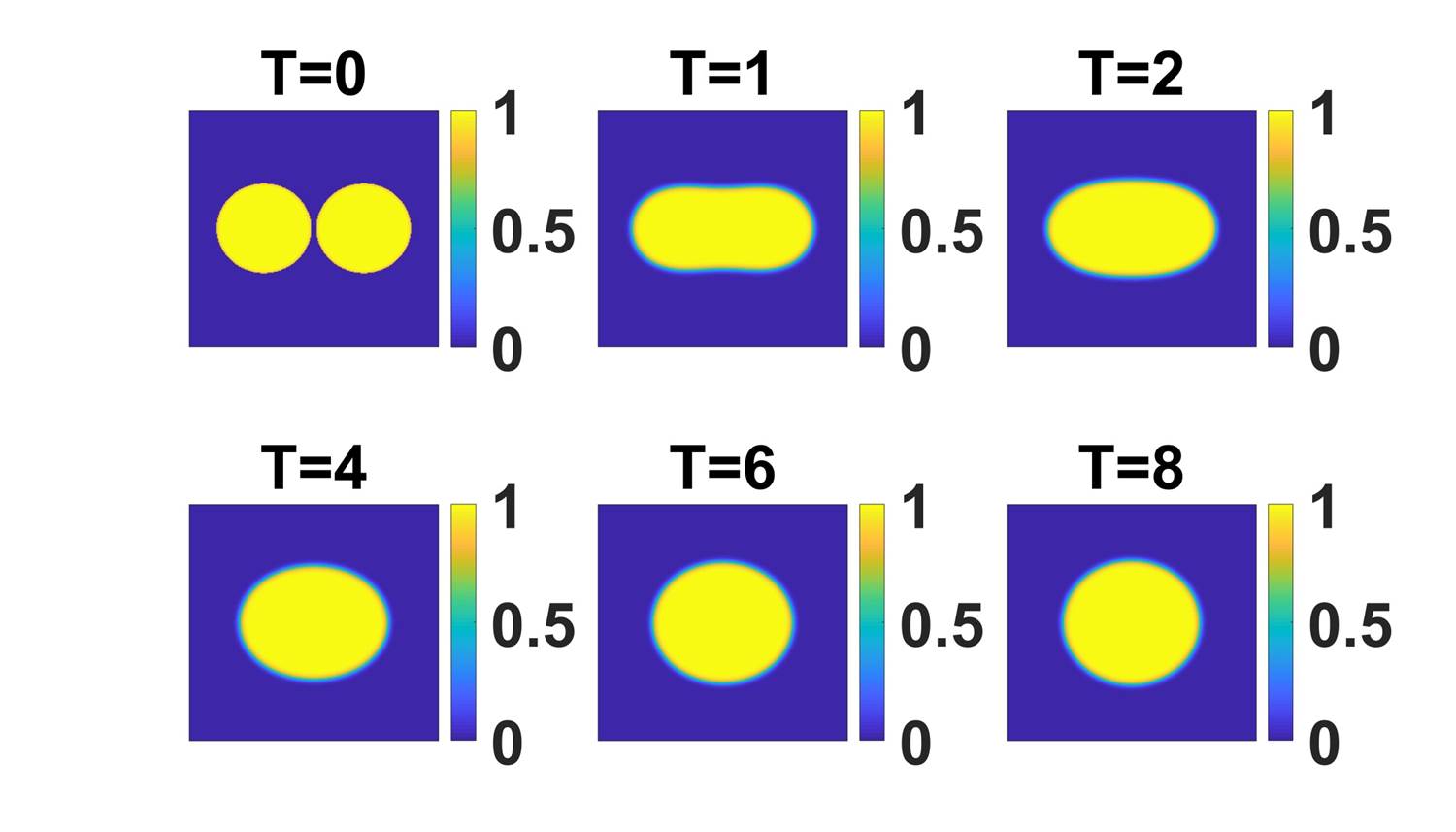}
\end{minipage}}
\subfigure[]{
\begin{minipage}[b]{0.49\linewidth}
\includegraphics[width=1\linewidth]{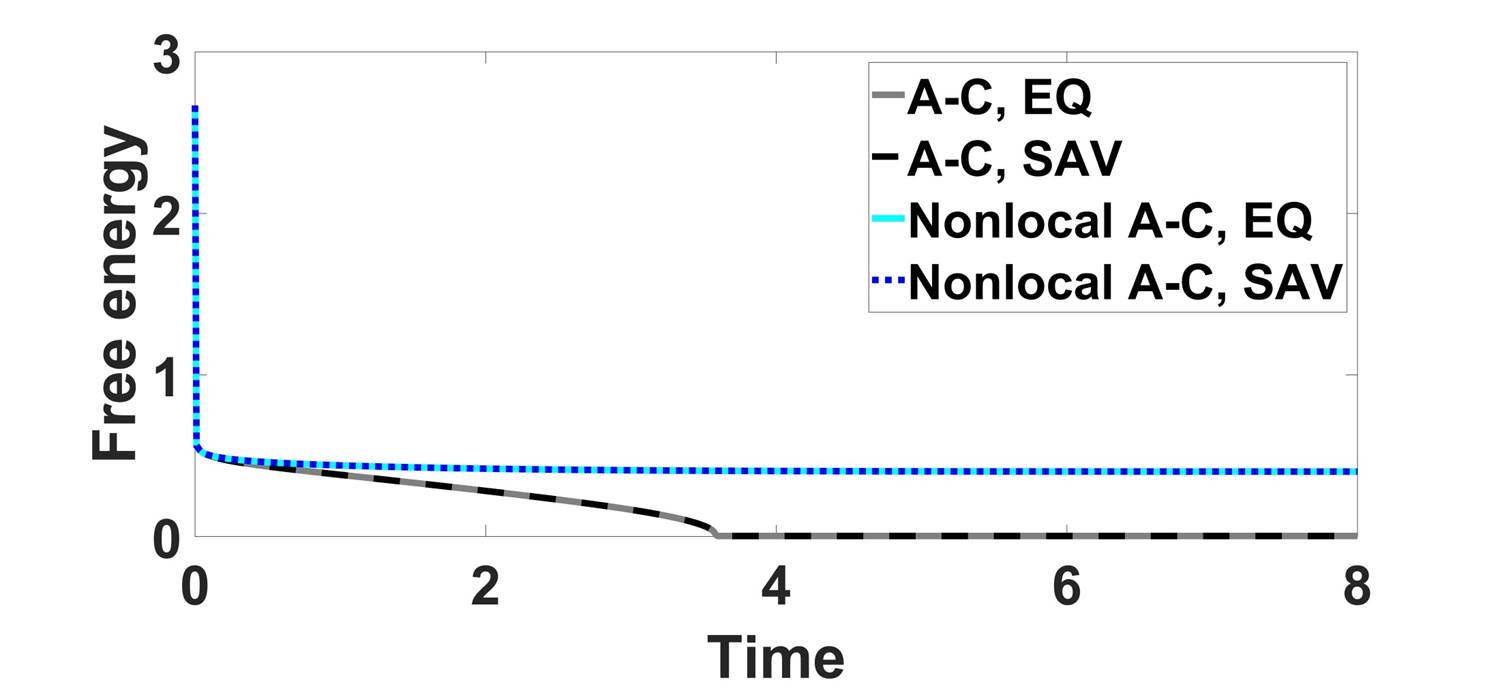}
\end{minipage}}
\subfigure[]{
\begin{minipage}[b]{0.49\linewidth}
\includegraphics[width=1\linewidth]{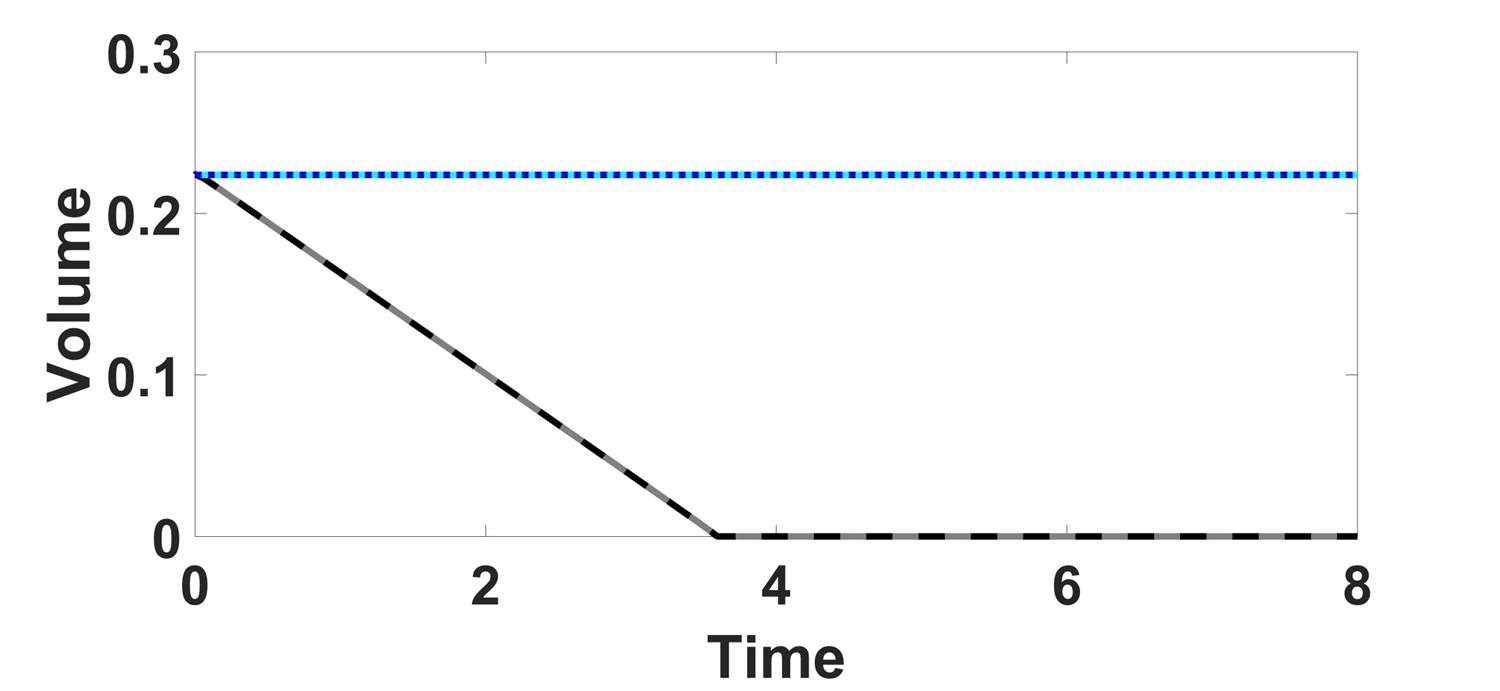}
\end{minipage}}
\caption{Merging of two droplets simulated using the Allen-Cahn and the Allen-Cahn models with nonlocal constraints at $M=1$. The dynamical behavior of the droplets of Allen-Cahn model and the Allen-Cahn model with nonlocal constraints is shown in (a) and (b) computed using AC-EQ and AC-L1-SAV scheme, respectively. Snapshots of the numerical approximation of volume fraction $\phi$ are taken at $\mathrm T=0, 0.8, 1.6, 2.4, 3.2, 4$, respectively,  in both cases. The time evolution of the free energy and phase volume are shown in (c) and (d), respectively.  Since all other models except for the Allen-Cahn model predict the similar dynamical behavior, we only show the phase transition dynamics computed using the AC-L1-SAV scheme in (b). We compare the time evolutions of free energy and volume computed by the EQ and SAV schemes in (c) and (d). The volume computed using the Allen model vanishes before $\mathrm T=4$ and in the meantime the energy hits zero as well.   All nonlocal Allen-Cahn models preserve the  volume and dissipate energies in time.  $\gamma_1=5\times10^{-3}$, $\gamma_2=100$ are used in the computations. The temporal and spatial step are set as $\Delta t=1\times 10^{-5}$ and $h_x=h_y=1/256$, respectively.}\label{Fig1}
\end{figure*}

We first simulate merging of the two drops using the Allen-Cahn and the  Allen-Cahn model with nonlocal constraints, respectively, with $M=1$.  The results computed from the EQ and SAV schemes for the Allen-Cahn model are identical, likewise the results computed using the EQ and SAV schemes for the nonlocal Allen-Cahn models are identical. We don't see any differences between the results of the Allen-Cahn model with a penalizing term and those of the Allen-Cahn model of a Lagrange multiplier. Figure \ref{Fig1}-(a) depicts the results computed from AC-EQ scheme and Figure \ref{Fig1}-(b) shows the results computed from AC-L1-SAV. We choose these two simulation results to show as representative examples.

From the simulations, we observe that the drops computed using the Allen-Cahn model first merge into a single drop and then the drop dissipates until eventually vanishes at the end of the simulation; while the drops computed using the Allen-Cahn models with nonlocal constraints merge into a single drop and eventually rounded up at the end of the simulation.  Figure \ref{Fig1}-(c) and (d) depict the computed free energy and volume of drops using the two numerical schemes. Obviously, the volume decays in the Allen-Cahn model while conserved in the simulation of the other models. The free energy decays in the Allen-Cahn model with respect to time and vanishes as the drop disappears. In contrast, the free energy of the Allen-Cahn model with a nonlocal constraint saturates at a nonzero value at the end of the simulation.

%  \begin{figure}[htbp]
%  \centering
%  \includegraphics[width=1 \textwidth]{Figure1.jpg}
%\caption{Merging of droplets simulated by the Allen-Cahn and the Allen-Cahn models with nonlocal constraints at $M=1$. The dynamical behaviors of the droplets of Allen-Cahn model and the Allen-Cahn model with nonlocal constraints are shown in A and B, respectively. Snapshots of the numerical approximation of $\phi$ are taken at t=0, 0.8, 1.6, 2.4, 3.2, 4 in both cases. The time evolution of the free energy and volume are shown in (C) and (D), respectively. We show the phase transition of Allen-Cahn model simulated by EQ scheme in A. Since all other models except for the Allen-Cahn model predict the similar dynamical behavior, we only show the phase transition dynamics of the Lagrangian model with $h'(\phi)=1$ simulated by SAV scheme in B. We compared the time evolutions of free energy and volume simulated by EQ and SAV schemes for the models in C and D. All nonlocal Allen-Cahn models preserve the same volume and dissipate in time. We set $\gamma_1=5\times10^{-3}$, $\gamma_2=100$. The pre-factor of the penalizing term in Penalizing model is set as $1\times 10^5$. The time step and space step are set as $\Delta t=1\times 10^{-5}$ and $\Delta x=\Delta y=1/256$, respectively.}\label{Fig1}
%\end{figure}
\begin{figure*}
\centering
\subfigure[]{
\begin{minipage}[b]{0.49\linewidth}
\includegraphics[width=1\linewidth]{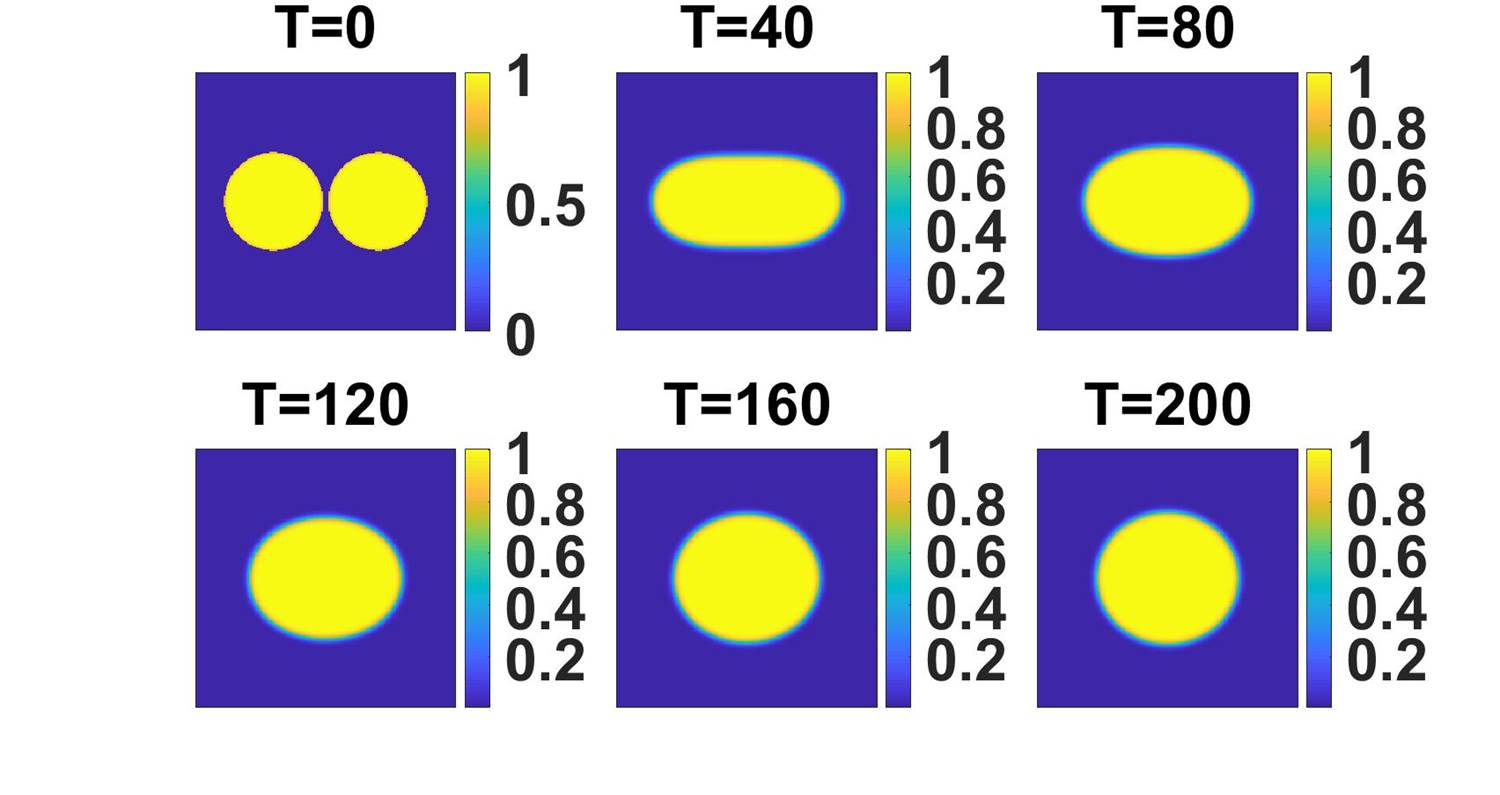}
\end{minipage}}
\subfigure[]{
\begin{minipage}[b]{0.49\linewidth}
\includegraphics[width=1\linewidth]{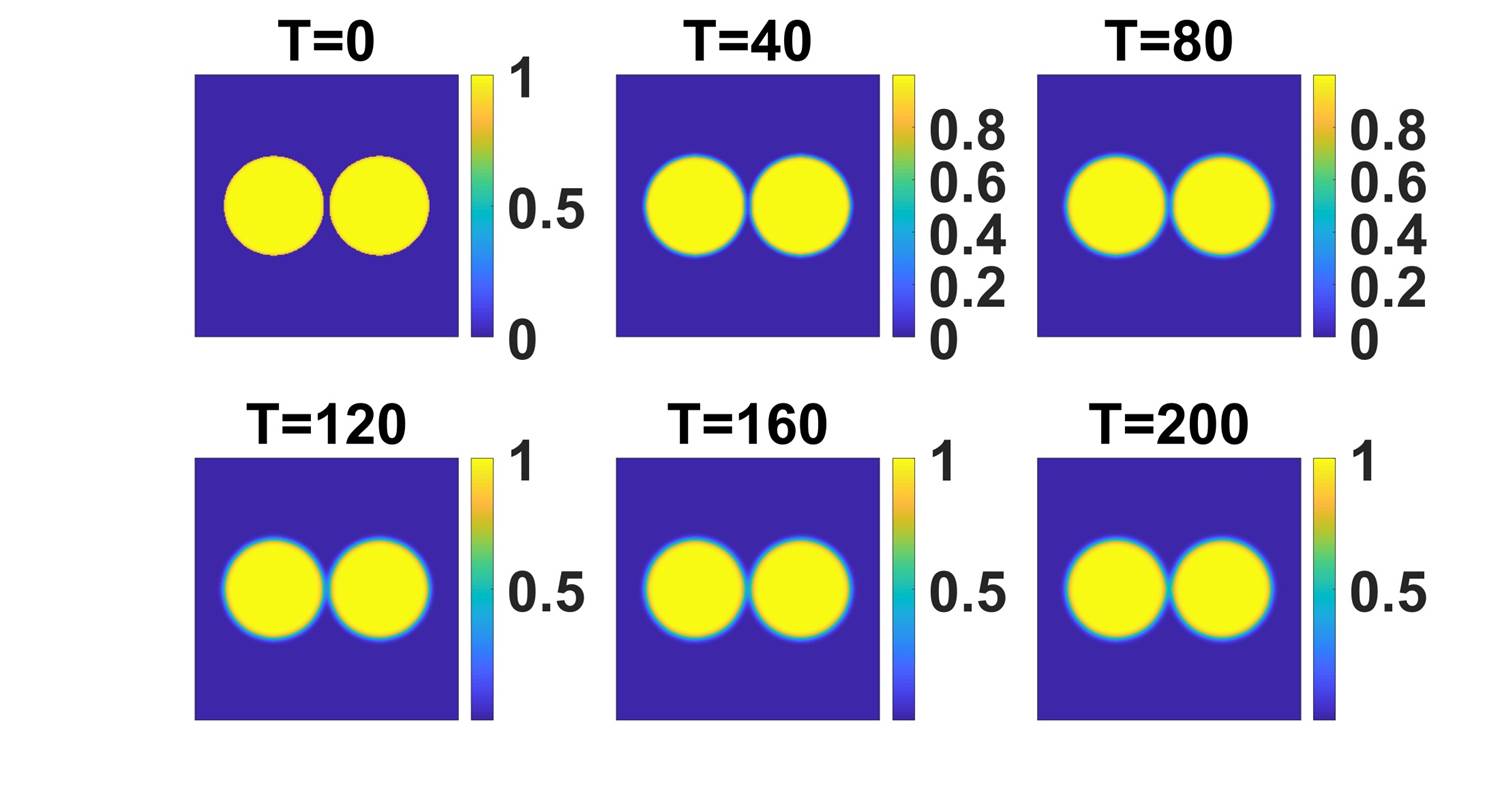}
\end{minipage}}
\subfigure[]{
\begin{minipage}[b]{0.49\linewidth}
\includegraphics[width=1\linewidth]{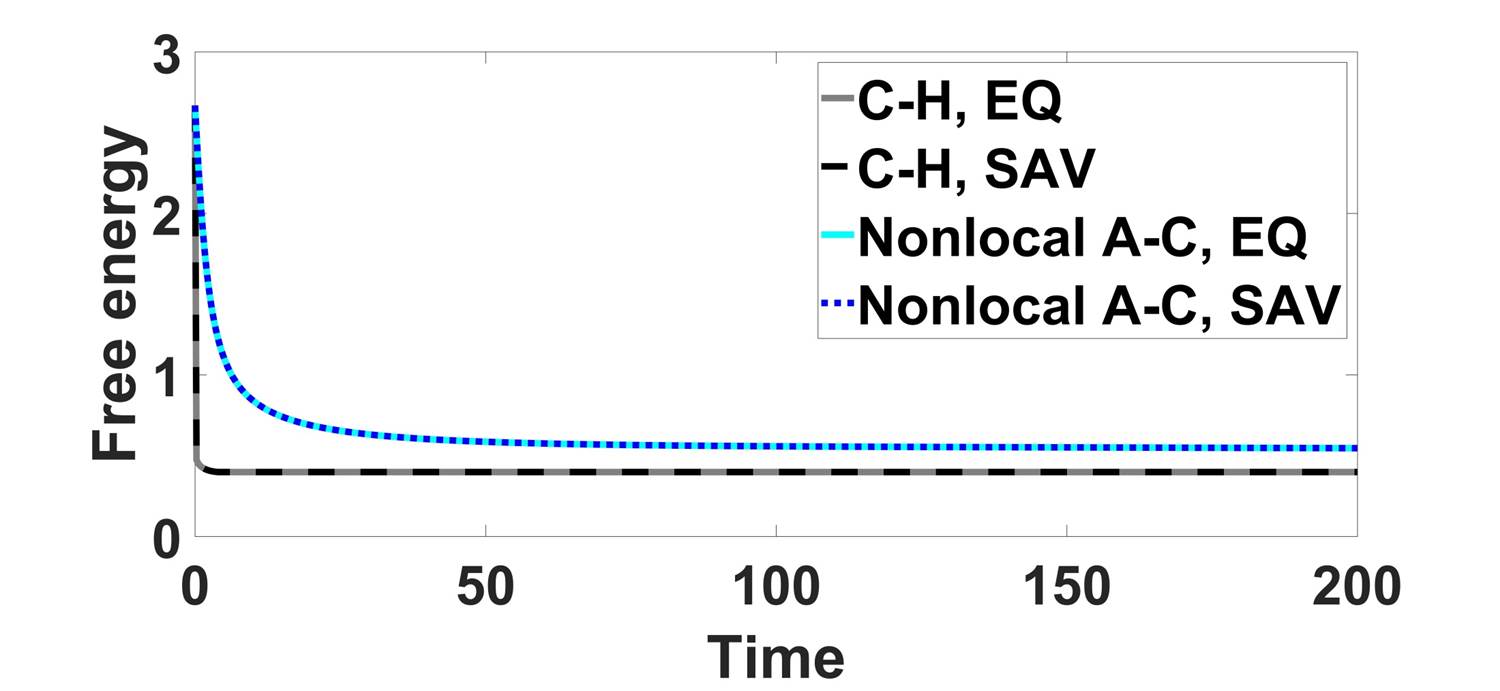}
\end{minipage}}
\subfigure[]{
\begin{minipage}[b]{0.49\linewidth}
\includegraphics[width=1\linewidth]{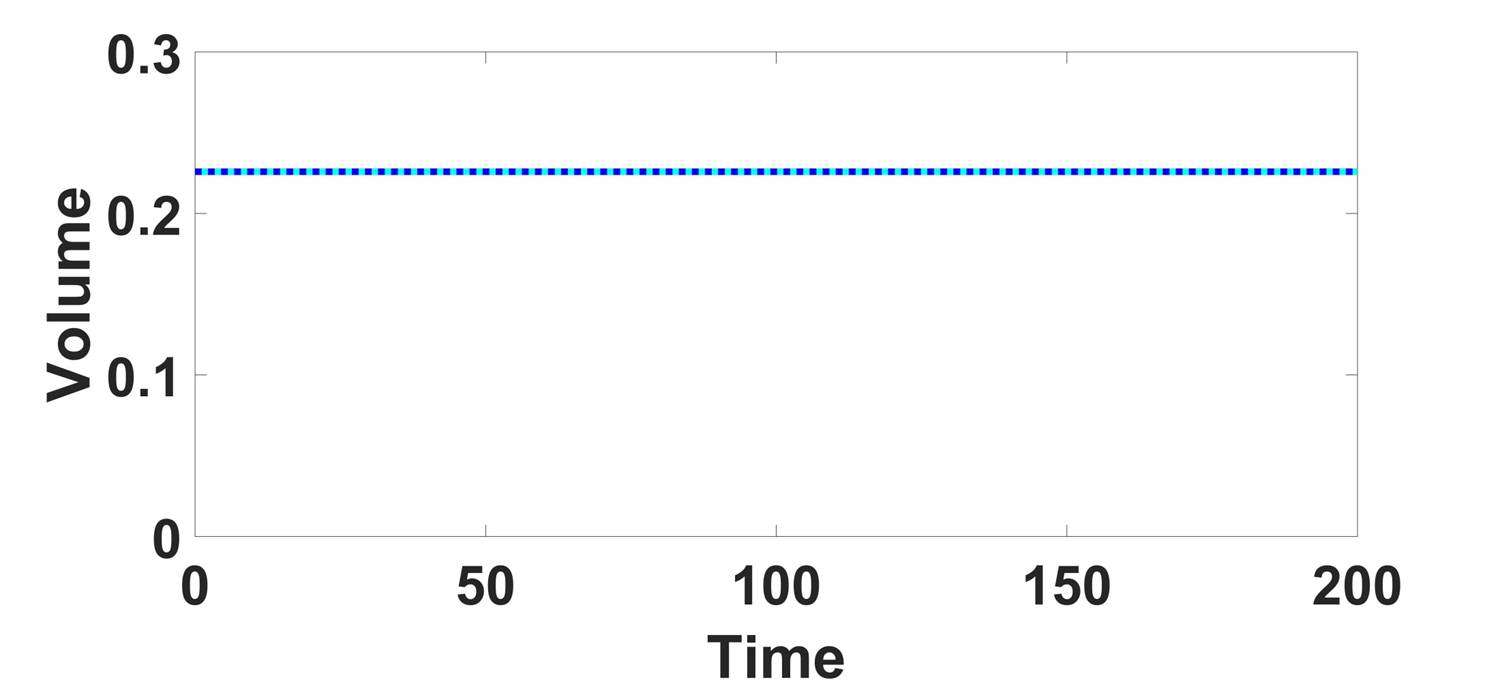}
\end{minipage}}
\caption{Merging of two droplets simulated using the Cahn-Hilliard and the Allen-Cahn models with nonlocal constraints at $M=1\times10^{-4}$. The phase evolution of the droplets of Cahn-Hilliard model and the Allen-Cahn model with nonlocal constraints is shown in (a) and (b) using the CH-EQ and AC-L1-SAV scheme, respectively. Snapshots of the numerical solution of $\phi$ are taken at $\mathrm T=0, 40, 80, 120, 160, 200$, respectively,  in both cases. The time evolution of the free energy and the volume are shown in (c) and (d), respectively. The Cahn-Hilliard model and the nonlocal Allen-Cahn models preserve the same volume and dissipate energies in time. Since the energy dissipates is faster in the Cahn-Hilliard model than in the nonlocal Allen-Cahn models initially, the energy profile predicted by the former is lower than the one predicted by the latter in the simulations. We set $\gamma_1=5\times10^{-3}$, $\gamma_2=100$ in the simulations.   The temporal and spatial steps are set as $\Delta t=1\times 10^{-5}$ and $h_x=h_y=1/256$, respectively.}\label{Fig2}
\end{figure*}

Then, we repeat the simulations using the  Allen-Cahn model with nonlocal constraints and the Cahn-Hilliard model with mobility coefficient  $M=1\times10^{-4}$. Since the coarsening rate in the Cahn-Hilliard model is much faster than that in the nonlocal Allen-Cahn models. At $\mathrm T=200$, the drops described by the Cahn-Hilliard model have merged into a single rounded drop, while the drops described by the nonlocal Allen-Cahn model just begin fusing. Figure \ref{Fig2}-(a) and (b) depict a representative drop merging simulation using the CH-EQ scheme for the Cahn-Hilliard model and one using the AC-L1-SAV scheme for the Allen-Cahn with a nonlocal constraint, respectively. The time evolution of the free energy computed from the Cahn-Hilliard model and that from the nonlocal Allen-Cahn model is depicted in Figure \ref{Fig2}-(c) and (d), respectively, in which the Cahn-Hilliard  yields a smaller free energy than the nonlocal Allen-Model does. This is because at $\mathrm T=200$ the Cahn-Hilliard dynamics has come into the steady state comparing with the dynamics of the Allen-Cahn model with nonlocal constraints.

%  \begin{figure}[htbp]
%  \centering
%  \includegraphics[width=1 \textwidth]{Figure2.jpg}
%  \caption{Merging of droplets simulated by the Cahn-Hilliard and the Allen-Cahn models with nonlocal constraints at $M=1\times10^{-4}$. The dynamical behaviors of the droplets of Cahn-Hilliard model and the Allen-Cahn model with nonlocal constraints are shown in A and B. Especially, A  and B are simulated by SAV and EQ scheme for the Cahn-Hilliard model and Lagrangian model with $h'(\phi)=1$. Snapshots of the numerical approximation of $\phi$ are taken at t=0, 40, 80, 120, 160, 200 in both cases. The time evolution of the free energy and volume are shown in (C) and (D), respectively.  We compared the EQ and SAV schemes for the models in C and D respectively. All Cahn-Hilliard model and nonlocal Allen-Cahn models preserve the same volume and dissipate in time. We set $\gamma_1=5\times10^{-3}$, $\gamma_2=100$. The pre-factor of the penalizing term in Penalizing model is set as $1\times 10^5$. The time step and space step are set as $\Delta t=1\times 10^{-5}$ and $\Delta x=\Delta y=1/256$, respectively.}\label{Fig2}
%\end{figure}

\begin{figure*}
\centering
\subfigure[]{
\begin{minipage}[b]{0.49\linewidth}
\includegraphics[width=1\linewidth]{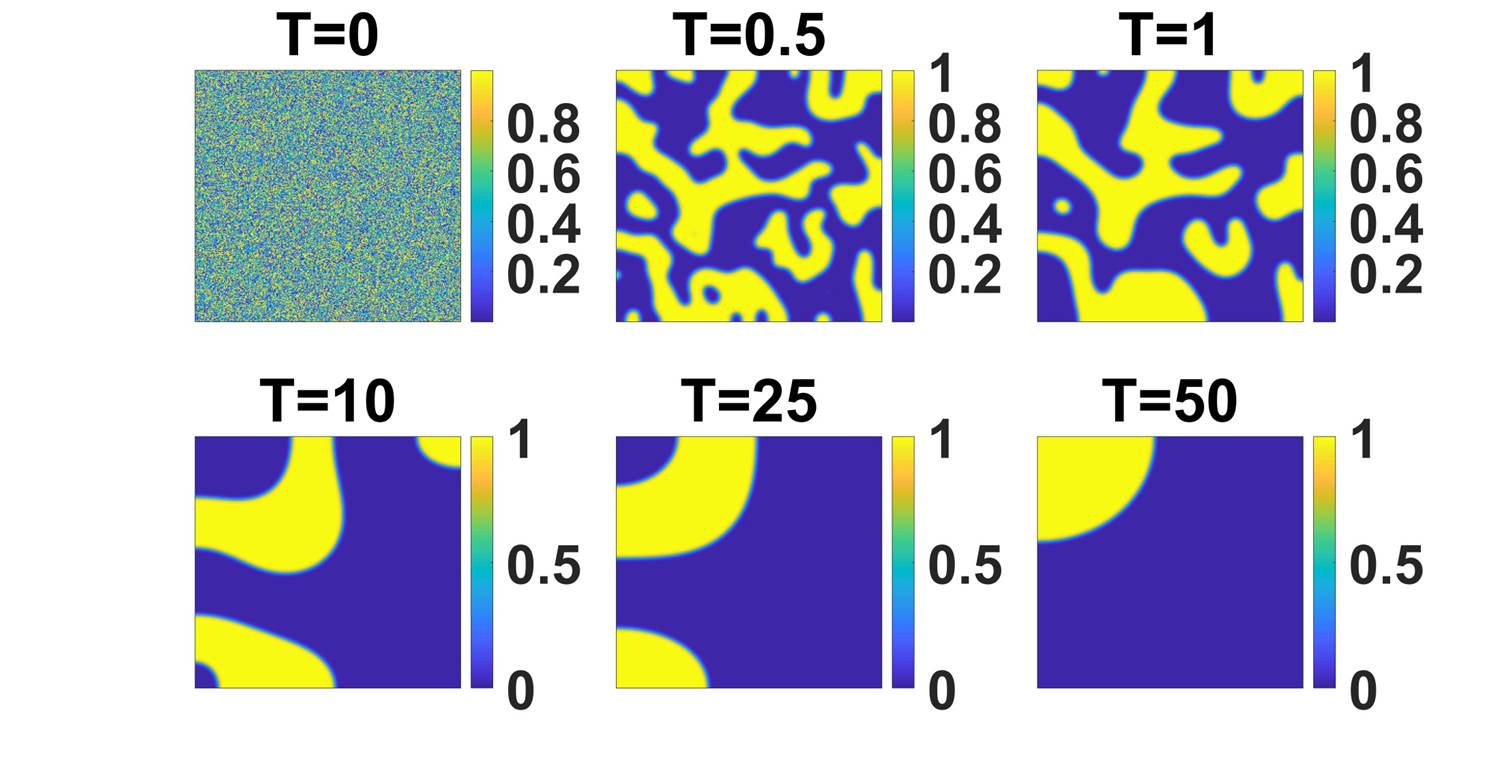}
\end{minipage}}
\subfigure[]{
\begin{minipage}[b]{0.49\linewidth}
\includegraphics[width=1\linewidth]{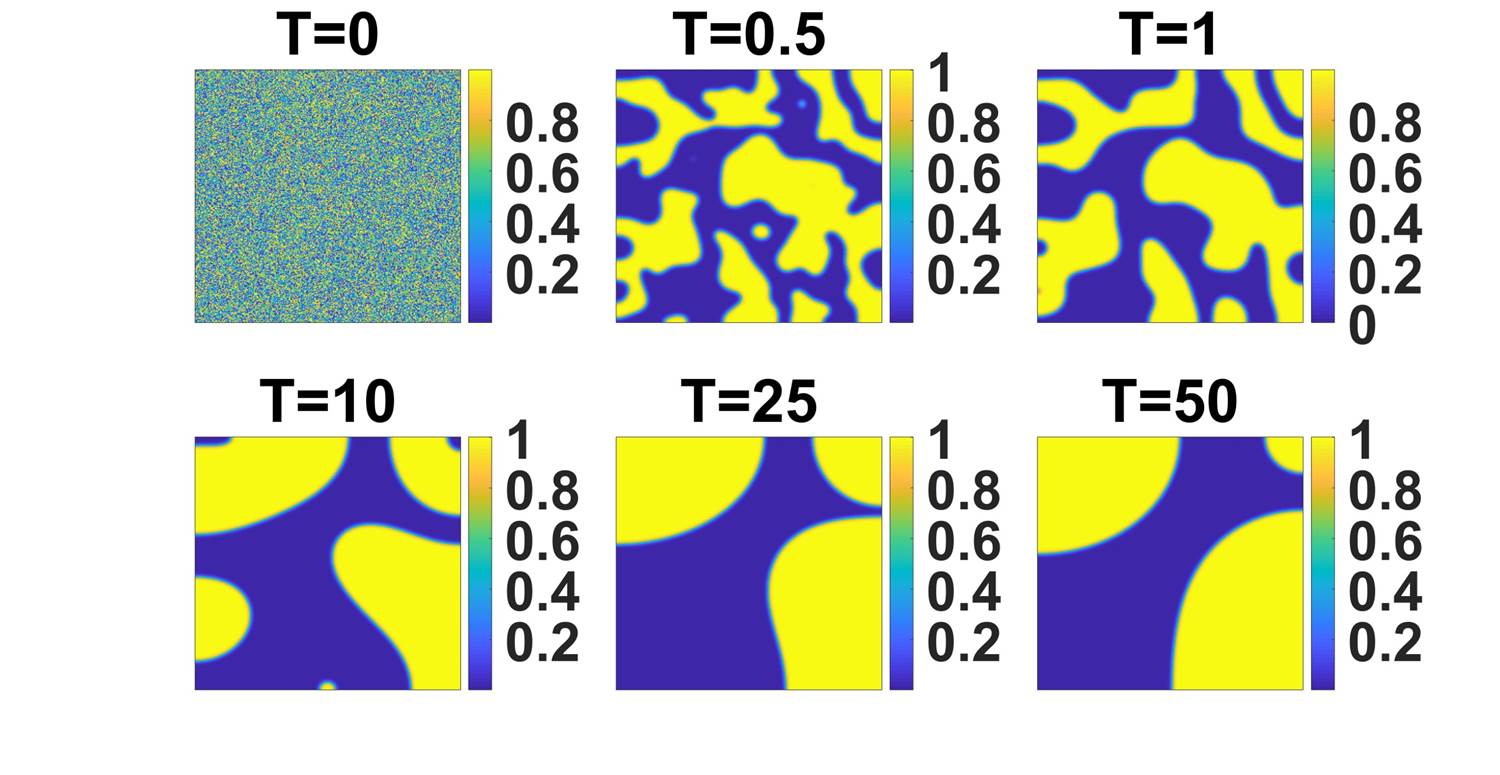}
\end{minipage}}
\subfigure[]{
\begin{minipage}[b]{0.49\linewidth}
\includegraphics[width=1\linewidth]{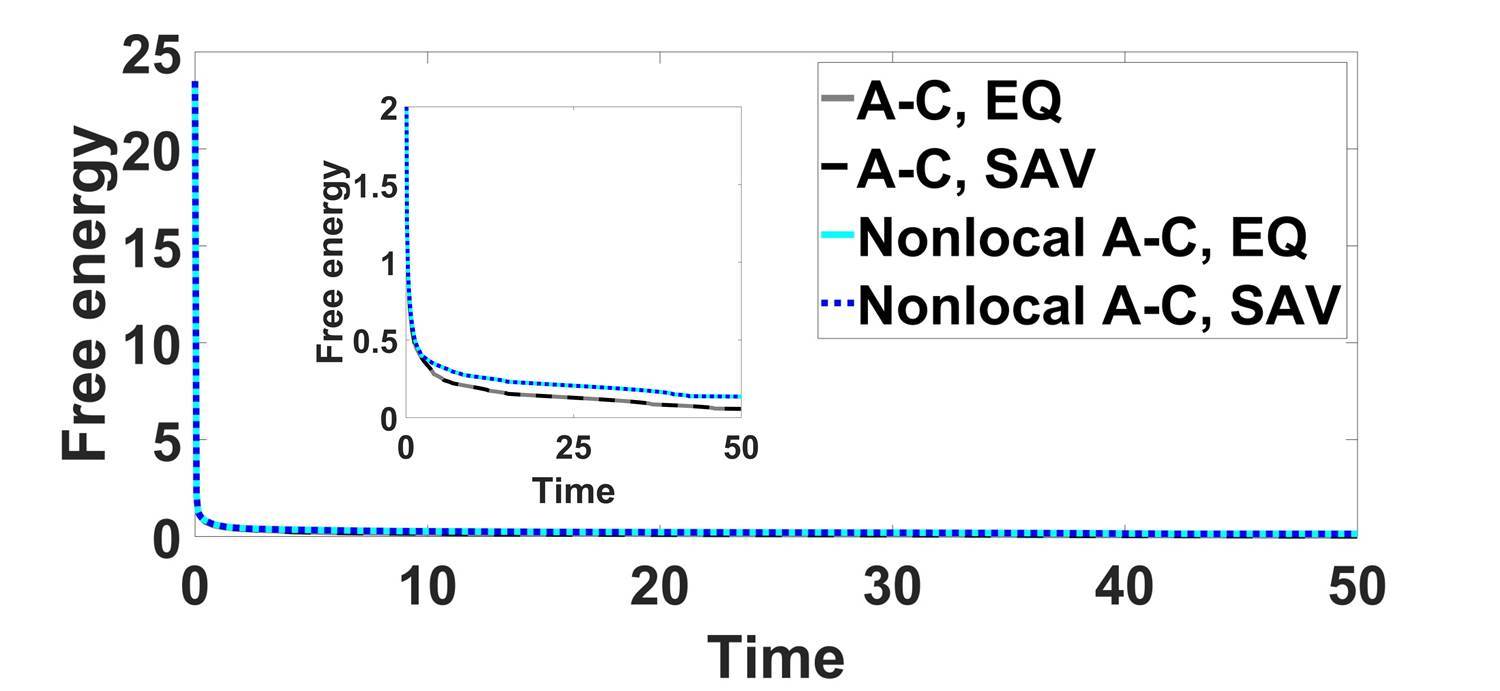}
\end{minipage}}
\subfigure[]{
\begin{minipage}[b]{0.49\linewidth}
\includegraphics[width=1\linewidth]{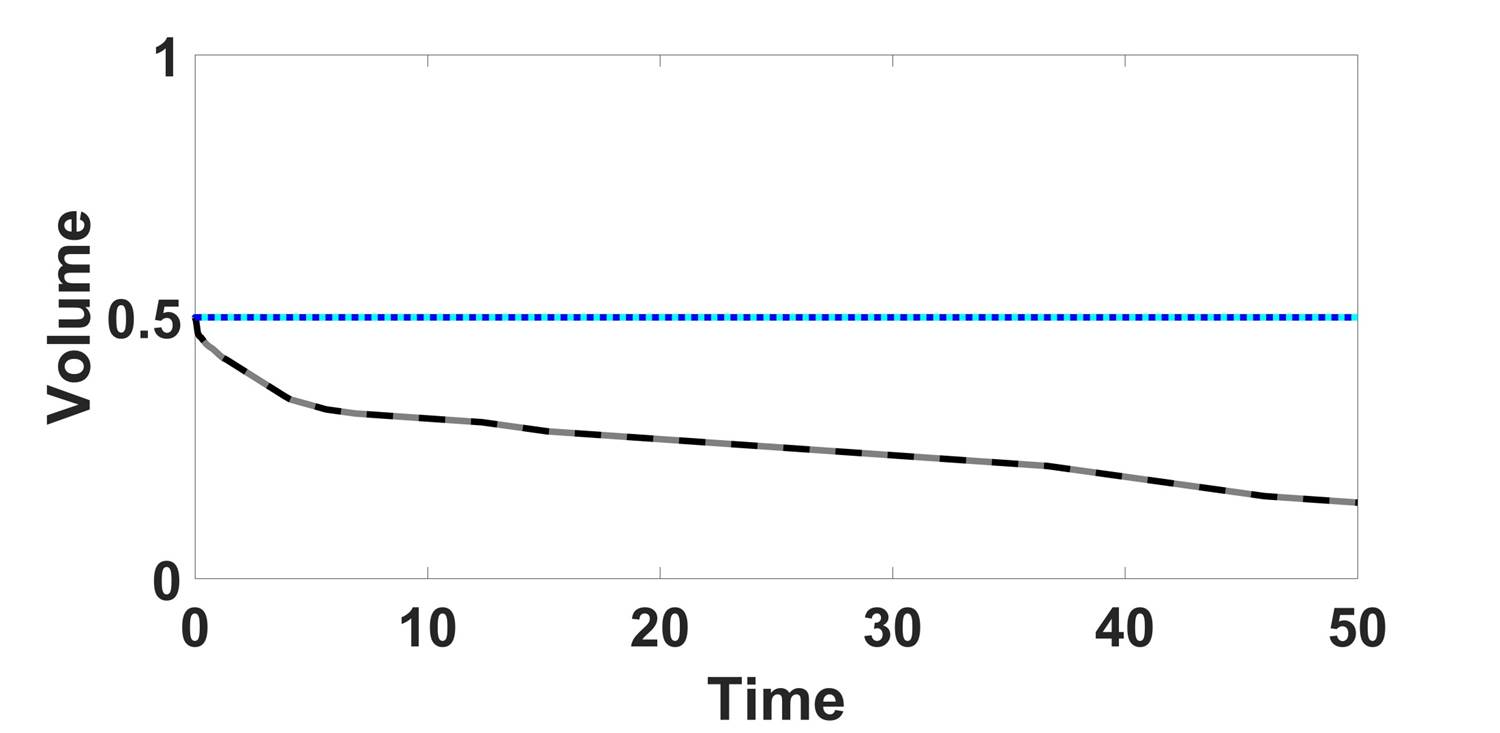}
\end{minipage}}
\caption{Coarsening dynamics  simulated using the Allen-Cahn and the Allen-Cahn models with nonlocal constraints at  $M=1$. The coarsening dynamics of Allen-Cahn model and the Allen-Cahn model with nonlocal constraints are shown in (a) and (b) using the AC-EQ and AC-L1-EQ scheme, respectively. Snapshots of the numerical approximation of $\phi$ are taken at $\mathrm T=0, 0.5, 1, 10, 25, 50$, respectively,  in both cases. The time evolution of the free energy and the volume are shown in (c) and (d), respectively.   The Allen-Cahn model does not conserve the volume, while the nonlocal Allen-Cahn models preserve the  volume and dissipate energies in time. The energy obtained using the Allen-Cahn model is lower than that using the nonlocal Allen-Cahn models. We set $\gamma_1=1\times10^{-3}$, $\gamma_2=50$ in the computations.  The temporal and spatial  steps are set as $\Delta t=1\times 10^{-5}$ and $h_x=h_y=1/256$, respectively.}\label{Fig3}
\end{figure*}

Secondly we use a random initial condition to assess the property of volume preserving nonlocal Allen-Cahn models in phase coarsening dynamics. Once again, The Allen-Cahn model gives one phase diagram at $t=50$, while the  Allen-Cahn models with nonlocal constraints yield another at the same parameter values and initial conditions. There is simply no comparison between these two model predictions in the terminal phase diagram. Figure \ref{Fig3}-(a) and (b) depict a typical simulation using the AC-EQ scheme for the Allen-Cahn model and the AC-L1-EQ scheme for the Allen-Cahn model with a Lagrangian multiplier, respectively. The time evolution of the free energy and the volume computed from the numerical schemes for the Allen-Cahn and the nonlocal Allen-Cahn models are shown in Figure \ref{Fig3}- (c) and (d), respectively. Between the EQ and SAV schemes of the same model, we have yet seen any visible differences between the numerical results.
%  \begin{figure}[htbp]
%  \centering
%  \includegraphics[width=1 \textwidth]{Figure3.jpg}
%   \caption{Coarsening dynamics of the system simulated by the Allen-Cahn and the Allen-Cahn models with nonlocal constraints at $M=1$. The dynamical behaviors of the droplets of Allen-Cahn model and the Allen-Cahn model with nonlocal constraints are shown in A and B, simulated by corresponding EQ schemes. Snapshots of the numerical approximation of $\phi$ are taken at t=0, 0.5, 1, 10, 25, 50 in both cases. The time evolution of the free energy and volume are shown in (C) and (D), respectively. Especially, since all other models except for the Allen-Cahn model predict the similar dynamical behavior, we only show the phase transition dynamics of the Penalizing model in B. We compared the EQ and SAV schemes for the models respectively. All nonlocal Allen-Cahn models preserve the same volume and dissipate in time. We set $\gamma_1=1\times10^{-3}$, $\gamma_2=50$. The pre-factor of the penalizing term in Penalizing model is set as $1\times 10^5$. The time step and space step are set as $\Delta t=1\times 10^{-5}$ and $\Delta x=\Delta y=1/256$, respectively.}\label{Fig3}
%\end{figure}

\begin{figure*}
\centering
\subfigure[]{
\begin{minipage}[b]{0.49\linewidth}
\includegraphics[width=1\linewidth]{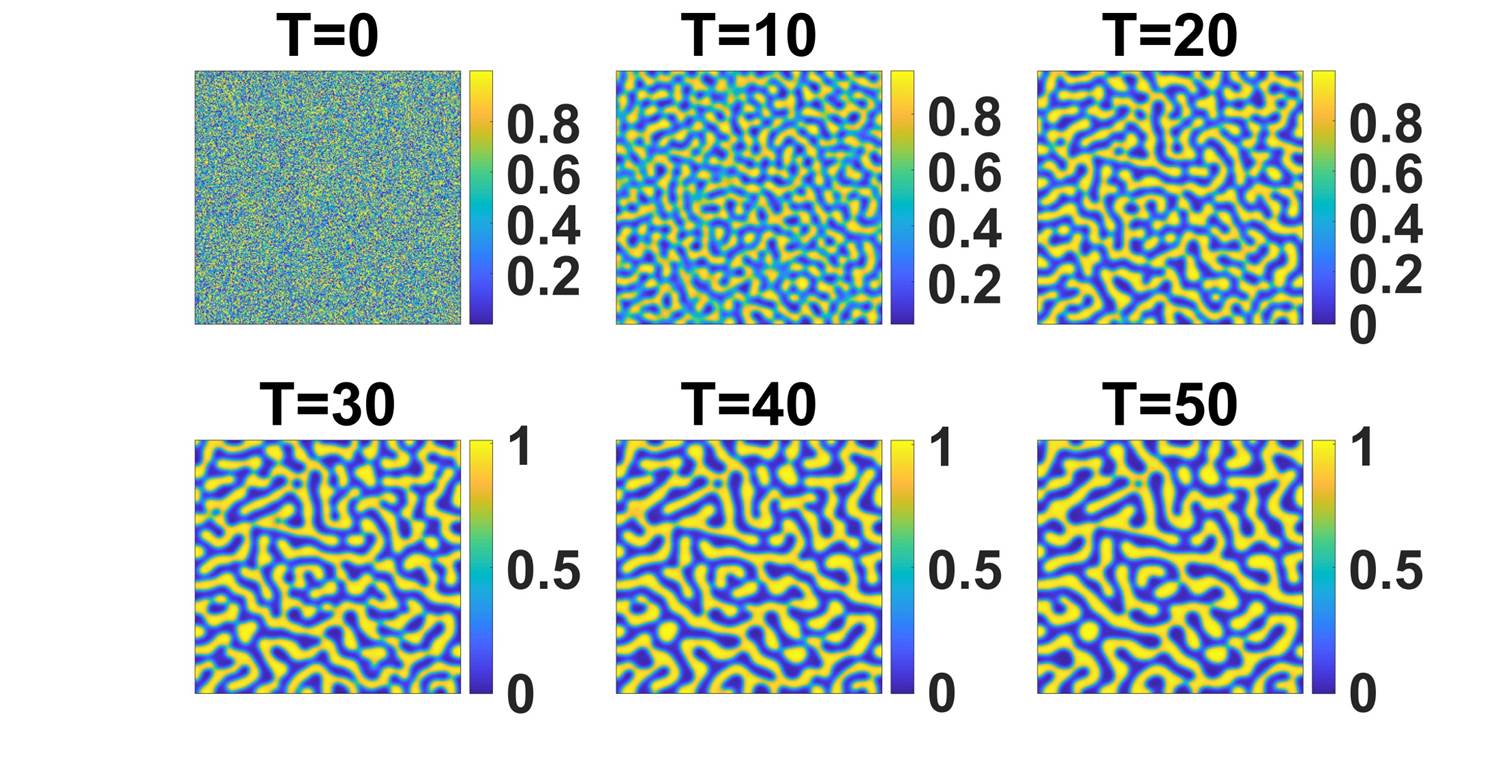}
\end{minipage}}
\subfigure[]{
\begin{minipage}[b]{0.49\linewidth}
\includegraphics[width=1\linewidth]{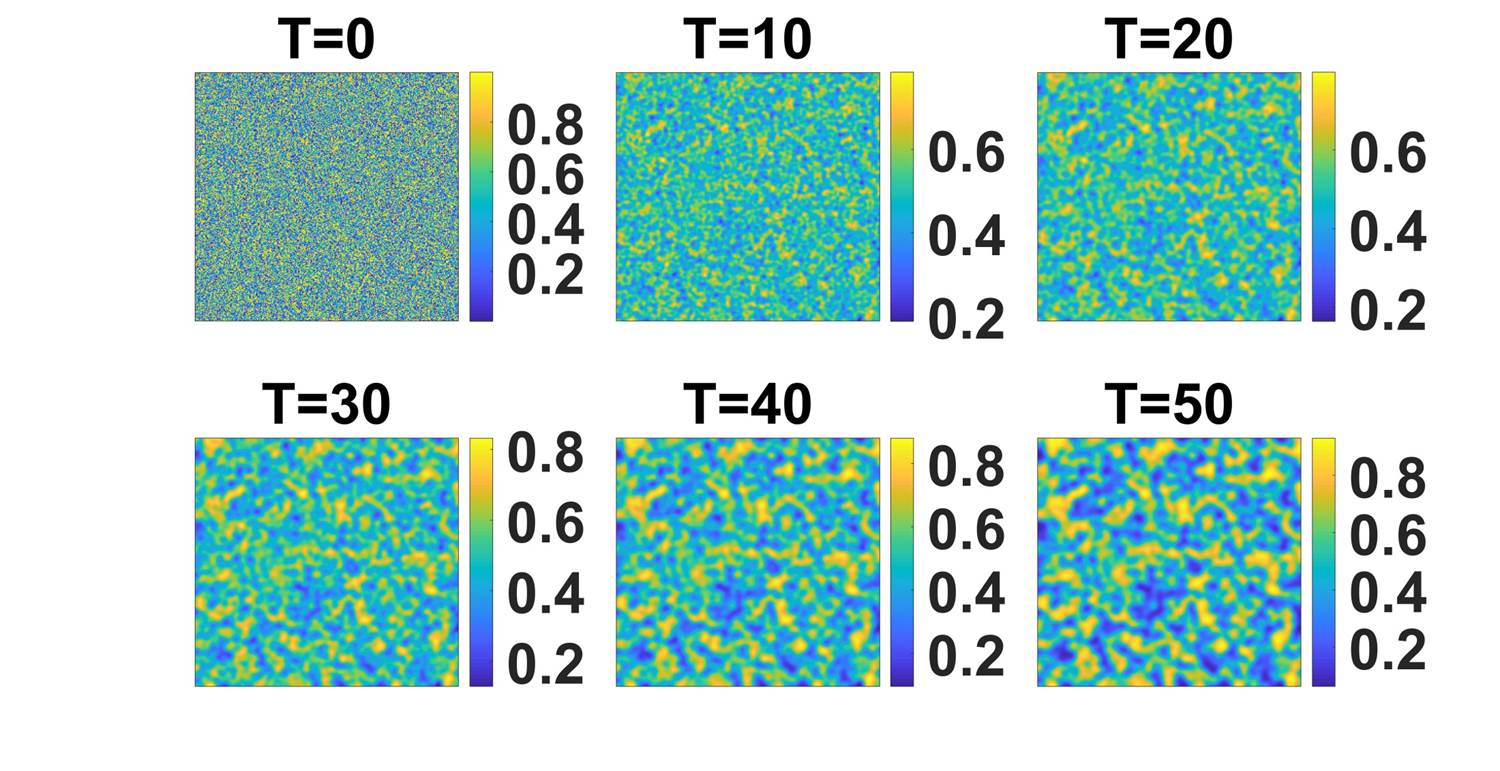}
\end{minipage}}
\subfigure[]{
\begin{minipage}[b]{0.49\linewidth}
\includegraphics[width=1\linewidth]{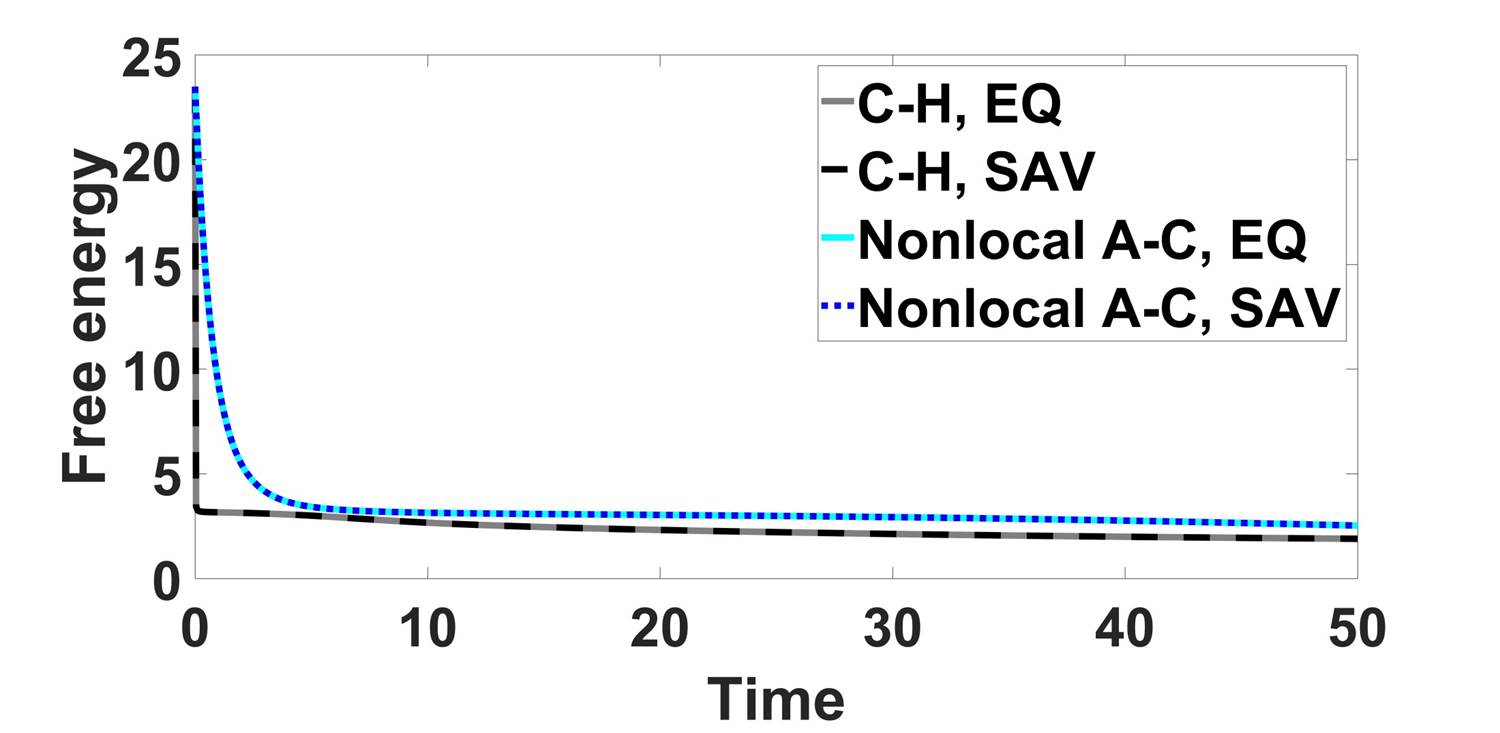}
\end{minipage}}
\subfigure[]{
\begin{minipage}[b]{0.49\linewidth}
\includegraphics[width=1\linewidth]{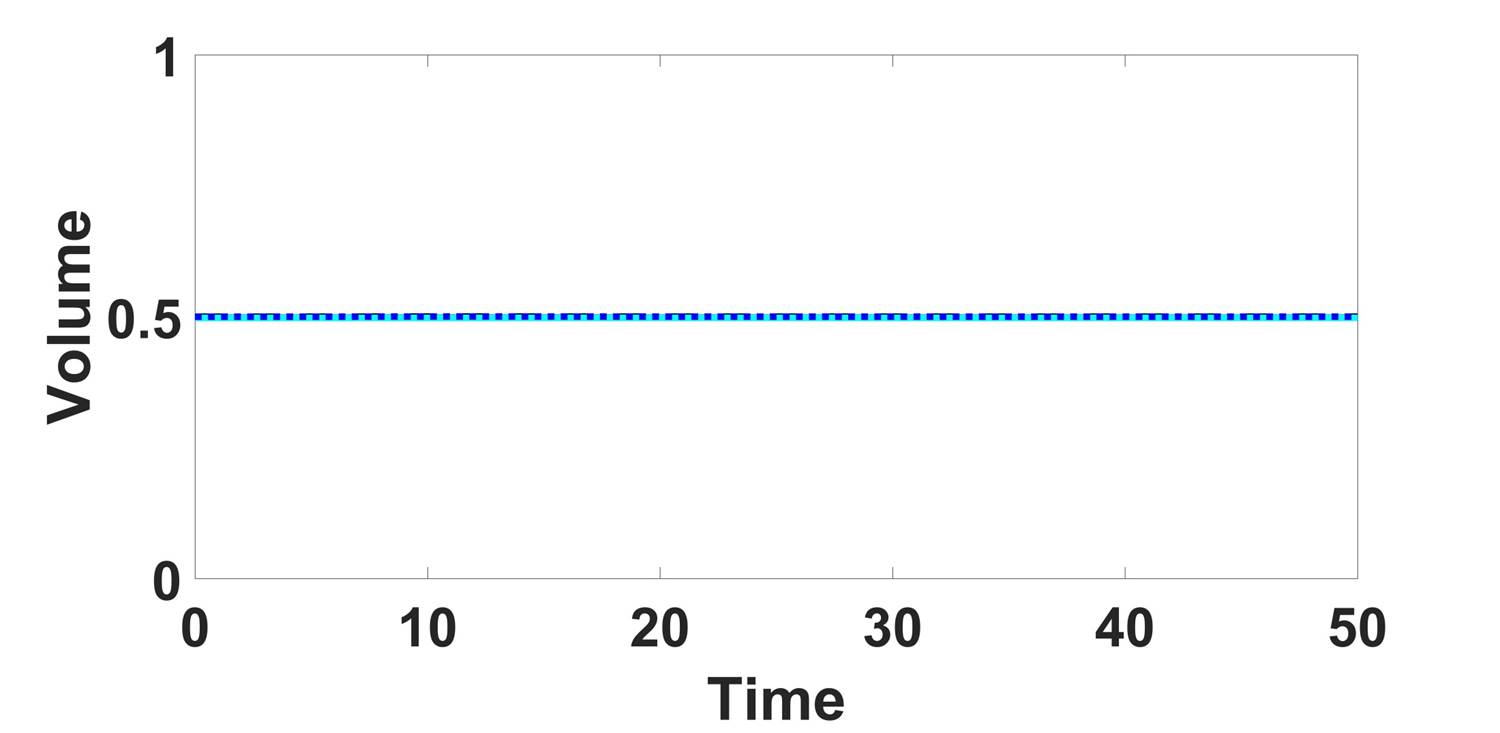}
\end{minipage}}
\caption{Coarsening dynamics simulated using the Cahn-Hilliard model at $M=1\times 10^{-6}$ and the Allen-Cahn models with nonlocal constraints at $M=1\times 10^{-3}$. The coarsening dynamics of the Cahn-Hilliard model and the Allen-Cahn model with nonlocal constraints are shown in (a) and (b), simulated by corresponding SAV schemes (AC-SAV and AC-L1-SAV), respectively. Snapshots of the numerical solutions of $\phi$ are taken at $\mathrm T=0, 10, 20, 30, 40, 50$ in both cases. The time evolution of the free energy and volume are shown in (c) and (d), respectively.  All  models preserve the  volume and dissipate energies in time. The faster coarsening dynamics in the Cahn-Hilliard model makes its energy slightly slower than the one predicted using the nonlocal Allen-Cahn models. We set $\gamma_1=1\times10^{-3}$, $\gamma_2=50$. The temporal and spatial step sizes are set as $\Delta t=1\times 10^{-5}$ and $h_x=h_y=1/256$, respectively.}\label{Fig4}
\end{figure*}

The  behavior of such coarsening dynamics is also compared between the Cahn-Hilliard model and the nonlocal Allen-Cahn models. Since the coarsening rate in the Cahn-Hilliard model is faster than that in the nonlocal Allen-Cahn models, we increase the magnitude of the mobility coefficient in the Allen-Cahn model with nonlocal constraints by 1000 folds and then repeat the simulations. This speeds up the coarsening dynamics of the nonlocal Allen-Cahn system significantly, although the result from the Cahn-Hilliard model still reaches a coarser grain than that of the Allen-Cahn models with nonlocal constraints. Figure \ref{Fig4}-(a) and (b) depict two representative examples on phase coarsening dynamics using CH-SAV and AC-L1-SAV schemes, respectively. Figure \ref{Fig4}-(c) and (d) show the decaying free energy  and the volume preserving results for the two selected simulations.\\
% \begin{figure}[htbp]
%  \centering
%  \includegraphics[width=1 \textwidth]{Figure4.jpg}
%  \caption{Coarsening dynamics of the system simulated by the Cahn-Hilliard at $M=1\times 10^{-6}$ and the Allen-Cahn models with nonlocal constraints at $M=1\times 10^{-3}$. The dynamical behaviors of the droplets of Cahn-Hilliard model and the Allen-Cahn model with nonlocal constraints are shown in A and B, simulated by corresponding SAV schemes. Snapshots of the numerical approximation of $\phi$ are taken at t=0, 10, 20, 30, 40, 50 in both cases. The time evolution of the free energy and volume are shown in (C) and (D), respectively. We compared the EQ and SAV schemes for the models respectively. All nonlocal Allen-Cahn models preserve the same volume and dissipate in time. We set $\gamma_1=1\times10^{-3}$, $\gamma_2=50$. The pre-factor of the penalizing term in Penalizing model is set as $1\times 10^5$. The time step and space step are set as $\Delta t=1\times 10^{-5}$ and $\Delta x=\Delta y=1/256$, respectively.}\label{Fig4}
%\end{figure}

\begin{figure}[htbp]
  \centering
  \includegraphics[width=0.65 \textwidth]{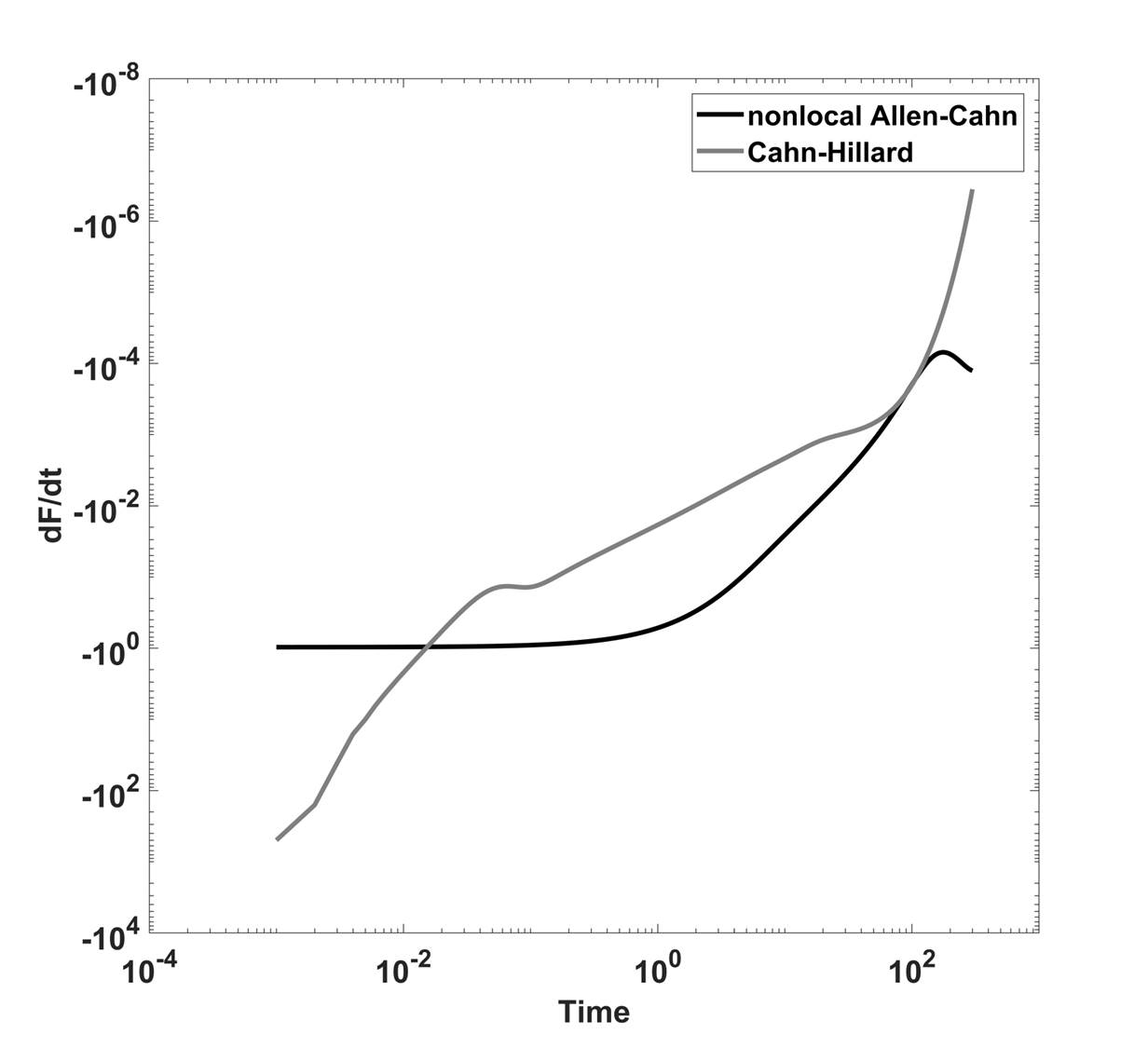}
   \caption{Time evolutions of the dissipation rate for the Allen-Cahn model with nonlocal constraints and the Cahn-Hilliard model (Simulated using AC-L1-SAV and CH-EQ, respectively). The dissipation rate of the energy is calculated from Fig \ref{Fig2}. At $t\leq 0.02$, the dissipation rate in the Cahn-Hilliard equation is larger than that in the Allen-Cahn equations with nonlocal constraints. This explains why the Cahn-Hilliard dynamics merges drops faster than that of the Cahn-Hilliard with nonlocal constraints does. }\label{Fig5}
\end{figure}

The results show the non-volume-conserving Allen-Cahn model can't be used to simulate the merging of droplets and the coarsening dynamics accurately, whereas the Allen-Cahn model with a penalizing potential and a Lagrangian multiplier can conserve the volume as the Cahn-Hilliard model does (see Fig \ref{Fig1}-\ref{Fig4}).  Compared with the Allen-Cahn model and the Cahn-Hilliard model, the Allen-Cahn models with nonlocal constraints not only conserves the volume of the density field in domain, but also shows a slower dissipation rate than both of them. One can enlarge the mobility of the  Allen-Cahn model with nonlocal constraints to accelerate the merging dynamics. We also compared the time evolution of the dissipation rates for the nonlocal Allen-Cahn models and the Cahn-Hilliard model with the same mobility coefficient $M$. The result indicates the dissipation rate of the free energy of nonlocal Allen-Cahn models is slower than that of the Cahn-Hilliard (see Fig \ref{Fig5}).

%%%%%%%%%%%%%%%%%%%%%%%%%% %%%%%%%%%%%%%%%%%%%%%%%%%%%%%%%%%%%%%%%%%%%%%%%%%%%%%%%%%%%%%%%%%%%%%%%%%%

\subsubsection{Practical implementation of the schemes}

\begin{figure*}
\centering
\subfigure[]{
\begin{minipage}[b]{1\linewidth}
\centering
\includegraphics[width=0.75\linewidth]{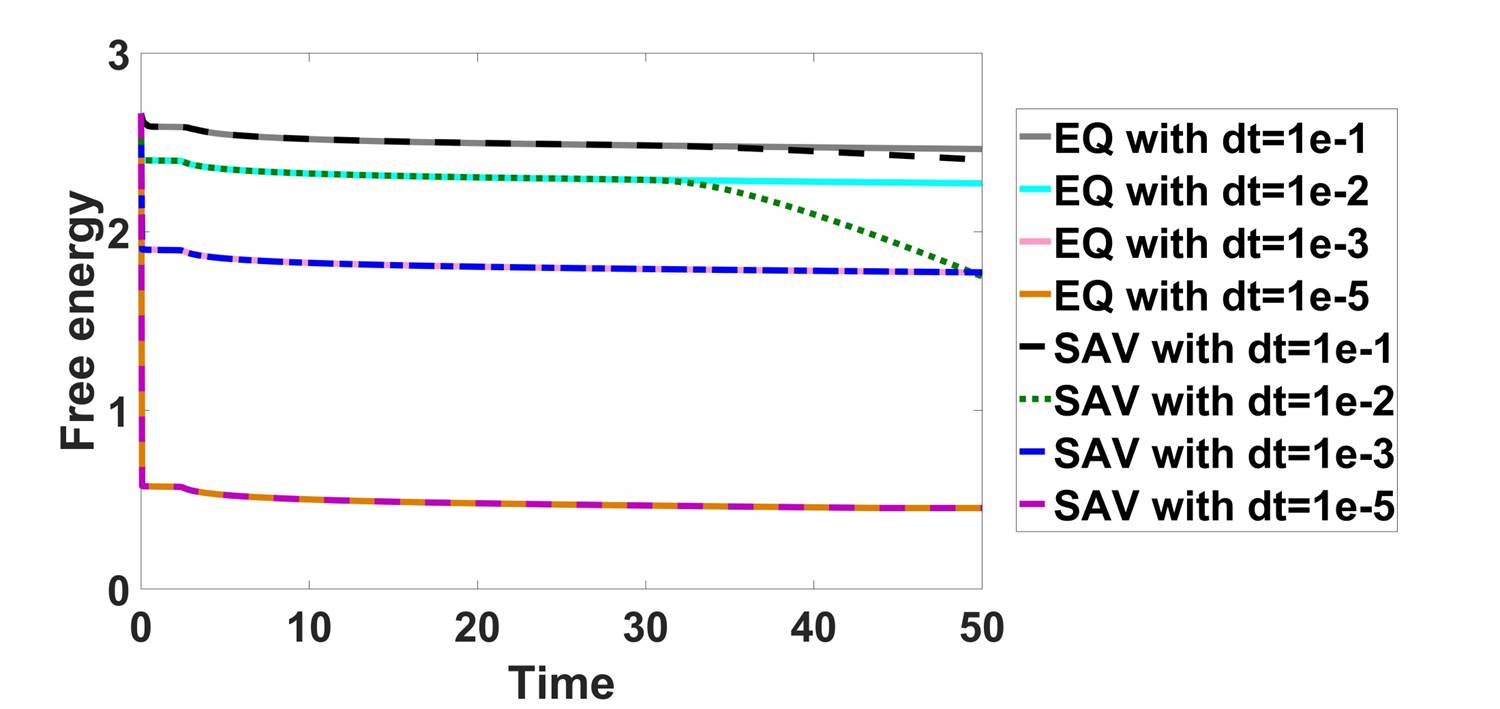}
\end{minipage}}
\subfigure[]{
\begin{minipage}[b]{1\linewidth}
\centering
\includegraphics[width=0.75\linewidth]{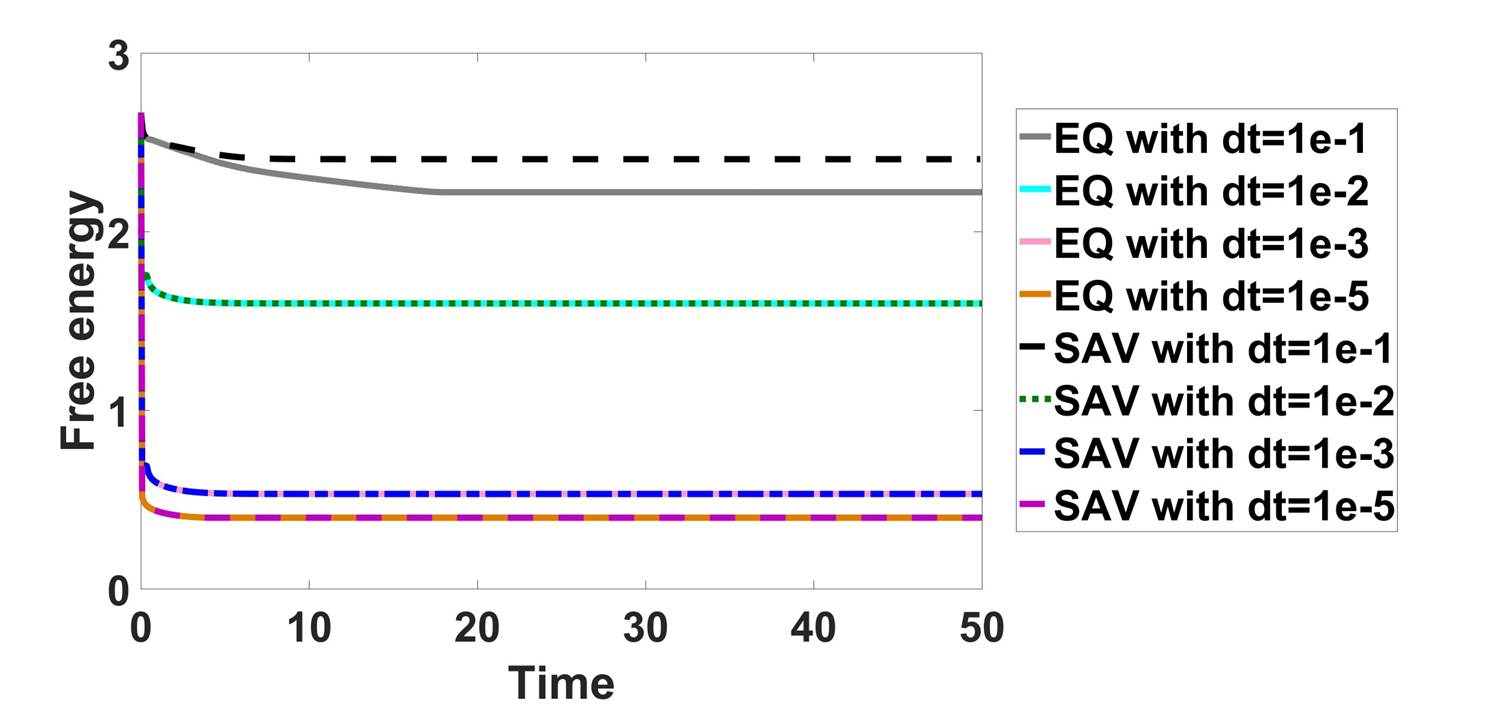}
\end{minipage}}
\caption{Accuracy of numerical solutions of the  Cahn-Hilliard model and the Allen-Cahn model with a Lagrange multiplier. (a) The free energy computed using the Cahn-Hilliard model at four selected time steps, respectively. (b). The free energy computed using AC-L1-EQ and AC-L1-SAV schemes at four selected time steps, respectively. The time step sizes used are  $1\times 10^{-1}$, $1\times 10^{-2}$, $1\times 10^{-3}$, $1\times 10^{-5}$. We use $M=1\times 10^{-4}$ for the Cahn-Hilliard model and $M=1$ for the nonlocal Allen-Cahn model. The initial condition and other model parameters are the same as those in \ref{Fig2}. In the computation using the Cahn-Hilliard model, the result converges at $\Delta t=1\times 10^{-5}$. The outcome is slightly better for the Allen-Cahn model with a lagrange multiplier. }\label{Fig6}
\end{figure*}

\noindent \indent  Although the schemes  are  shown  unconditionally energy stable, the numerical results are not guaranteed to be always accurate if the time step size is large due to the essentially sequential decoupling of the schemes in Fig \ref{Fig6}. In all these schemes based on energy quadratization strategies, the equation for the auxiliary variables or the intermediate variables are ordinary differential equations, which are derived from original algebraic equations (which define the intermediate variables) by taking time derivatives. The numerical schemes devised to solve these ordinary differential equations may not be accurate enough to warrant the accumulation of the energy accurately.

\begin{figure*}
\centering
\subfigure[]{
\begin{minipage}[b]{0.49\linewidth}
\centering
\includegraphics[width=0.75\linewidth]{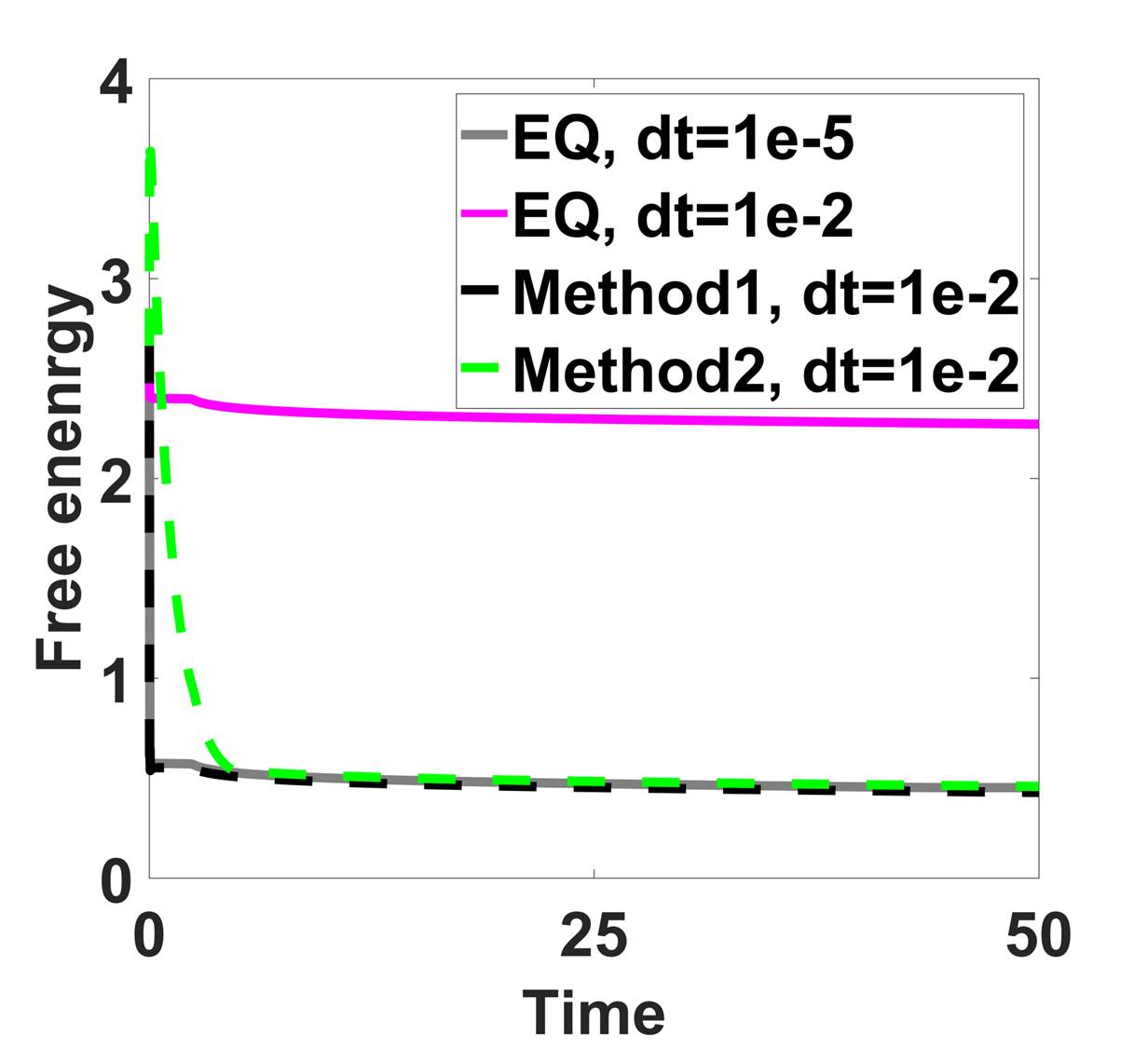}
\end{minipage}}
\subfigure[]{
\begin{minipage}[b]{0.49\linewidth}
\centering
\includegraphics[width=0.75\linewidth]{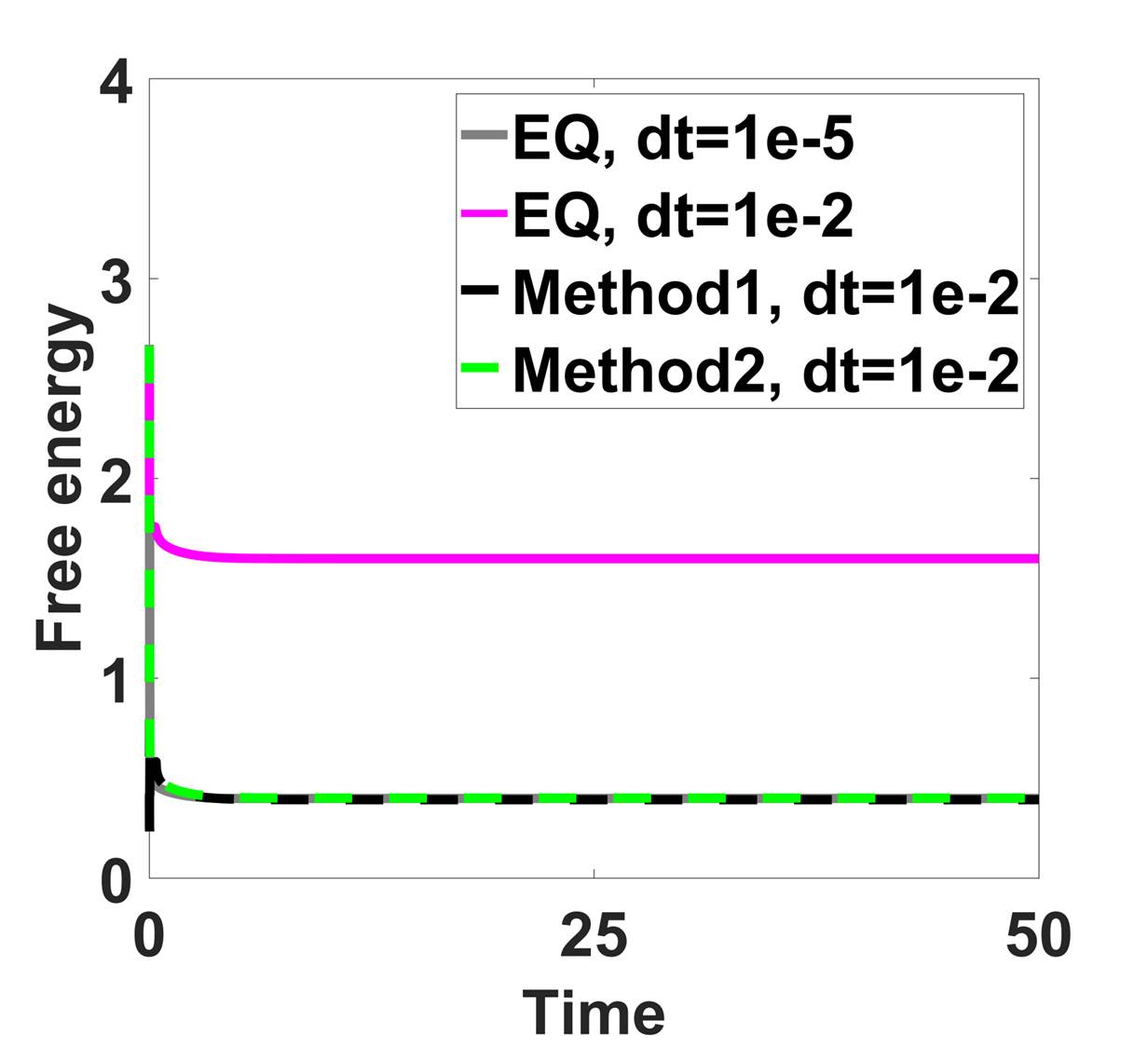}
\end{minipage}}
\caption{Accuracy of EQ schemes when implemented using the tricks. (a). Results computed using CH-EQ scheme. (b). Results computed using the  AC-L1-EQ scheme. We use $M=1\times 10^{-4}$, $\alpha=1, C=5\times 10^{-5}$ for the Cahn-Hilliard model and $M=1$, $\alpha=1, C=1.5\times 10^{-4}$  for the Allen-Cahn model with a Lagrange multiplier in the simulations.   The second numerical trick works better than the first one in the case of the nonlocal Allen-Cahn model, but the first method works better for the Cahn-Hilliard model.  The initial conditions and the parameter values are the same as those used in \ref{Fig2}.}\label{Fig7}
\end{figure*}

To remedy the inherent deficiency, we propose   two   methods to modify the schemes in practical implementations to improve their numerical  accuracy  with a large time step. For simplicity, we define $f_1(\phi^n)$ as the non-quadratic, nonlinear term in the bulk potential, which is $f_1(\phi^n)=f(\phi^n)-\gamma_2{(\phi^n)}^2$ in this paper. The two methods are given  below.
\begin{enumerate}
\item After obtaining $\phi^{n+1}$, we update $q^{n+1}$ using  $q^{n+1}=\sqrt{f_1(\phi^{n+1})+C_0}$.\\
\item  After obtaining $\phi^{n+1}$, if $\int_\Omega (q^{n+1}-\sqrt{f_1(\phi^{n+1})+C_0}\mathrm) d{\bf r}\leqslant \varepsilon$,
we update $q^{n+1}$ using  $q^{n+1}=q^n+\ohs{q'}(\phi^{n+1}-\phi^n)$, otherwise, using $q^{n+1}=\sqrt{f_1(\phi^{n+1})+C_0}-\frac{2-\alpha\Delta t}{2+\alpha\Delta t}(q^n-\sqrt{f_1(\phi^n)+C_0})$.
\end{enumerate}
We note that, in method 2, the result is insensitive to the choice of $\alpha>0$. Figure \ref{Fig6} depicts a pair of comparative studies on drop merging simulations using the Cahn-Hilliard model and the Allen-Cahn model with a Lagrange multiplier discretized using both EQ and SAV methods. The benchmark or good results are obtained using small time steps. Using these two tricks, we can alleviate the constraints imposed on the time step size for the EQ schemes considerably.  In Fig \ref{Fig7}, we show results of the numerical methods at a relatively large step size. The improvement is  significant.  We also test our tricks for the SAV schemes, but we observe  that only method 1 performs well.
%\begin{figure}[htbp]
%  \centering
%  \includegraphics[width=1\textwidth]{Figure6.jpg}
%  \caption{The effectiveness and accuracy of EQ approaches with a large time step for  Cahn-Hilliard model (A) and nonlocal Allen-Cahn model (B). We shown the time evolutions of free energy for the Cahn-Hilliard and nonlocal Allen-Cahn models. The free energy shown in B is the simulation result of the Lagrangian model. The time step shown in this figure is chose as $1\times 10^{-1}$, $1\times 10^{-2}$, $1\times 10^{-3}$, $1\times 10^{-5}$. We set $M=1\times 10^{-4}$ for the Cahn-Hilliard model and $M=1$ for the nonlocal Allen-Cahn model. The initiation conditions and other parameters are the same as those in \ref{Fig2}. In the Cahn-Hilliard model computations, the result converges at $\Delta t=1\times 10^{-5}$. The outcome is slightly better for the Allen-Cahn model with a lagrange multiplier. }\label{Fig6}
%\end{figure}

%\begin{figure}[htbp]
%  \centering
%  \includegraphics[width=1 \textwidth]{Figure7.jpg}
%   \caption{The modified numerical implementations of EQ schemes for Cahn-Hilliard model(A) and nonlocal Allen-Cahn model(B), respectively. We set $M=1\times 10^{-4}$ for the Cahn-Hilliard model and $M=1$ for the nonlocal Allen-Cahn model in the simulations. We use $\alpha=1, C=5\times 10^{-5}$ for the Cahn-Hilliard model and $\alpha=1, C=1.5\times 10^{-4}$ for the nonlocal Allen-Cahn model . The second numerical trick works better than the first one in the nonlocal Allen-Cahn cases, but the first method works better in the Cahn-Hilliard case.  The initiation conditions and parameters are same as that in \ref{Fig2}.}\label{Fig7}
%\end{figure}
In the numerical experiments, we observe that the Allen-Cahn model with a penalizing potential and the Allen-Cahn model with a Lagrange multiplier seem to render comparable numerical results.  In this study, we have used two different definitions of the phase volume $ V(t)=\int _\Omega h(\phi(t)) \mathrm {d {\bf r}}$ with two choices of $h(\phi)$:
\bena
h(\phi)=\phi,\\
h'(\phi)=\frac{(m+1)(2m+1)}{m}[\phi(1-\phi)]^m, m \hbox{ is a  positive integer}.
\eena
Our numerical experiments do not seem to be able to differentiate between the nonlocal Allen-Cahn models using either volume definitions. Thus both can be used at the discretion of the user in practice. Physically, the second  definition seems to be more sound because the compensation to the time change of the volume fraction  primarily takes place around the interface while the first definition seems to compensate the phase variable globally \cite{Rubinstein1992Nonlocal}.

\begin{figure*}
\centering
\subfigure[]{
\begin{minipage}[b]{0.49\linewidth}
\includegraphics[width=0.75\linewidth]{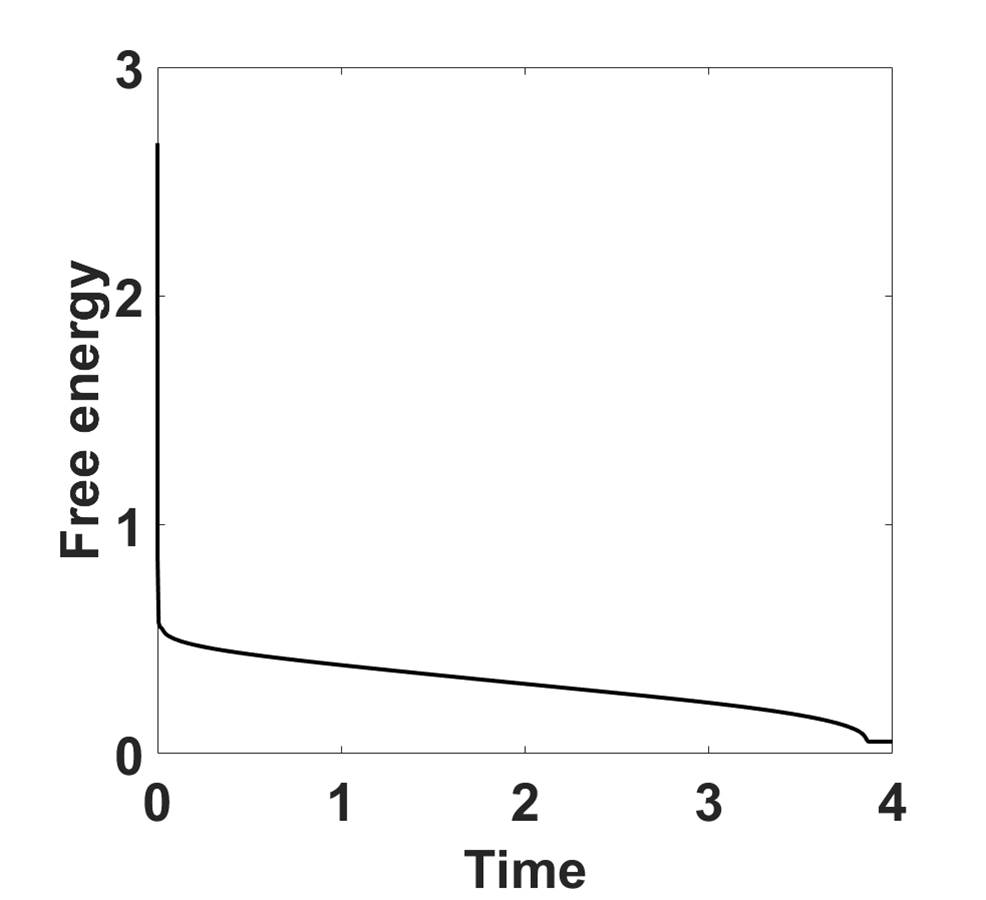}
\end{minipage}}
\subfigure[]{
\begin{minipage}[b]{0.49\linewidth}
\includegraphics[width=0.75\linewidth]{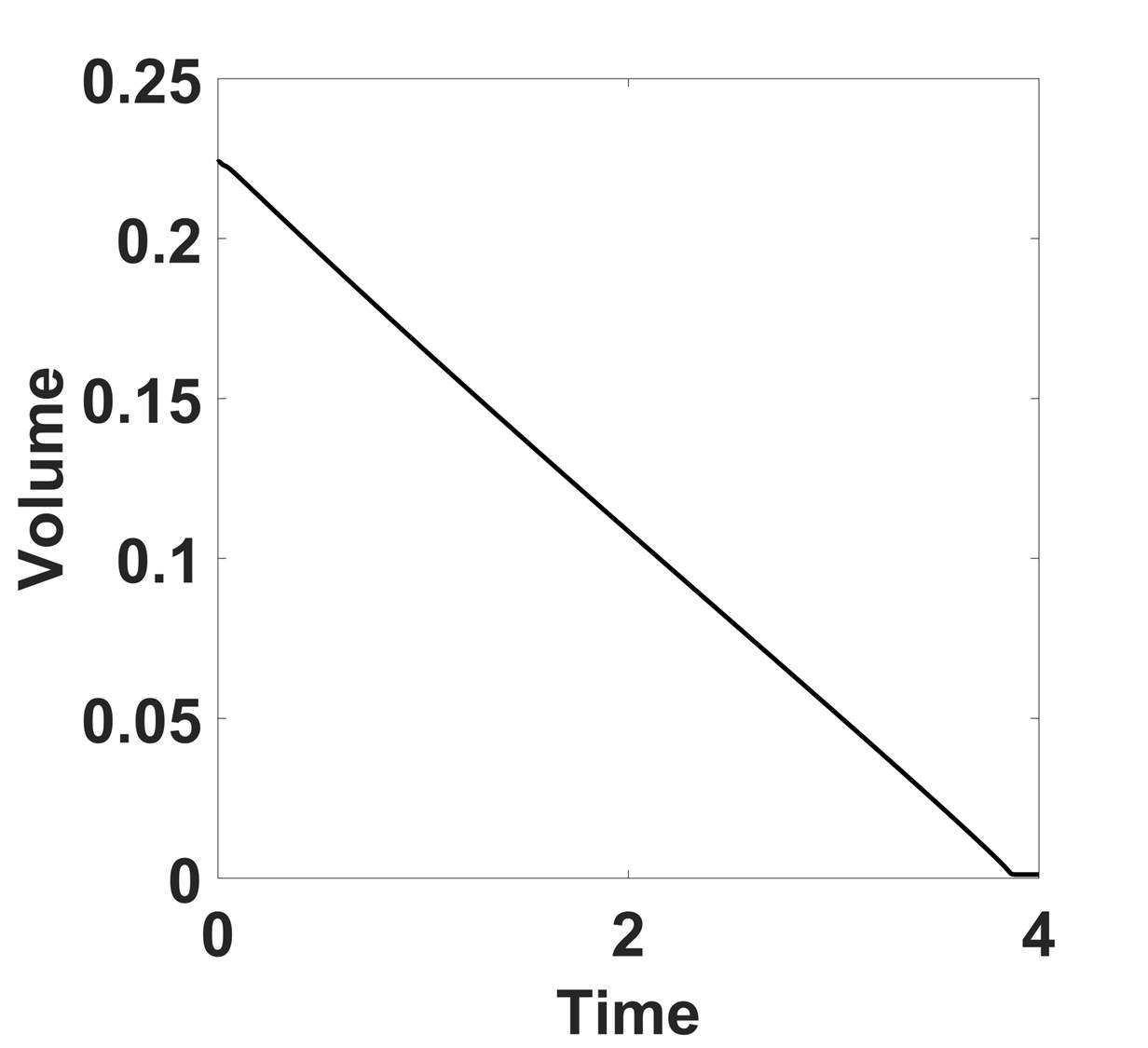}
\end{minipage}}
\subfigure[]{
\begin{minipage}[b]{0.49\linewidth}
\includegraphics[width=0.75\linewidth]{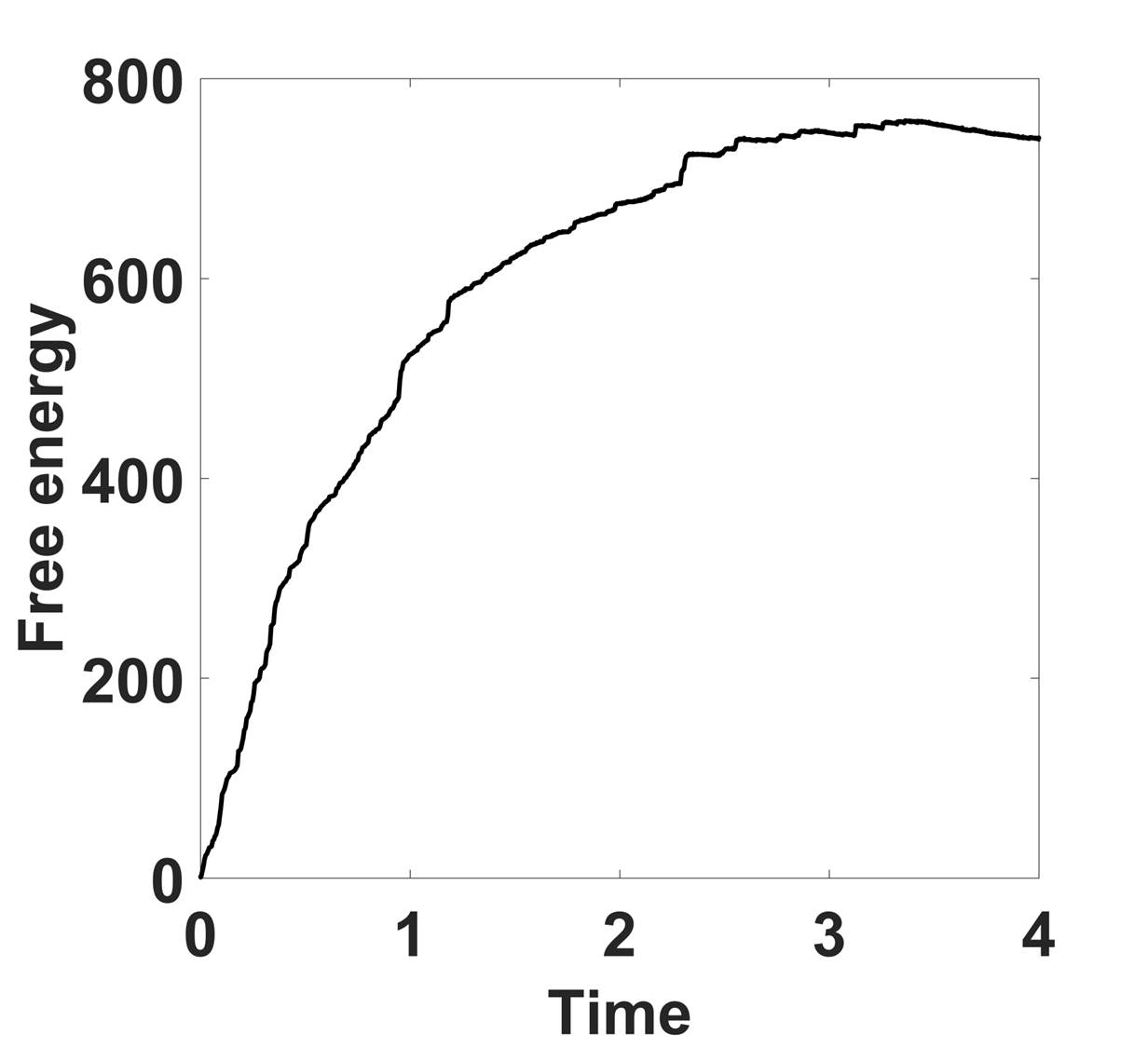}
\end{minipage}}
\subfigure[]{
\begin{minipage}[b]{0.49\linewidth}
\includegraphics[width=0.75\linewidth]{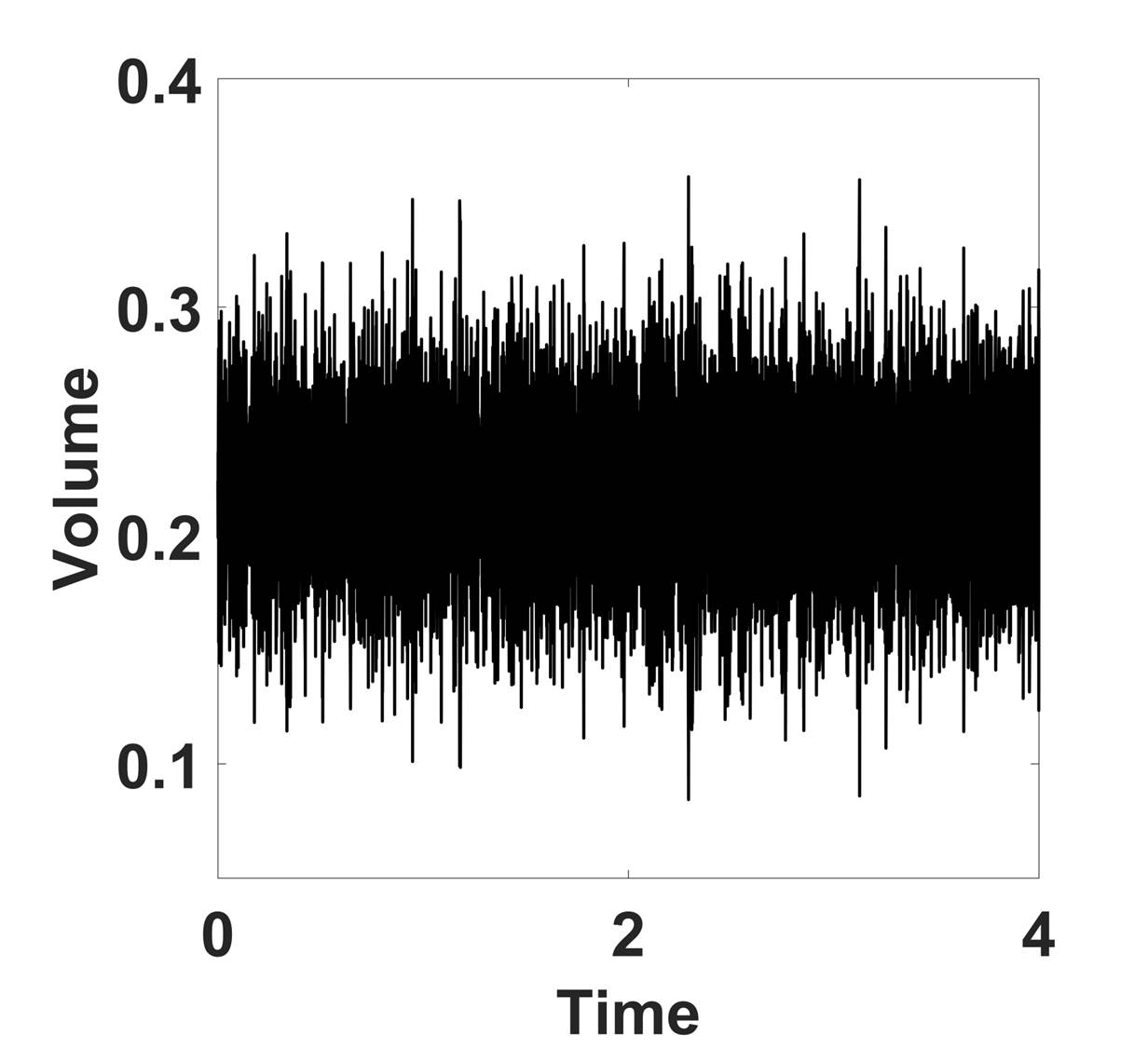}
\end{minipage}}
\caption{Two examples of failure for the Allen-Cahn model with a penalizing potential at $M=1$. $\eta=1$ and $\eta=1\times 10^{18}$ are chosen in the simulations. We plotted time evolution of the free energy and the volume, respectively. Results in (a) and (b) are obtained using the AC-P-EQ scheme with $\eta=1$.  Results in (c) and (d) are obtained using  the AC-P-SAV scheme with $\eta=1\times 10^{18}$. If $\eta$ is not "properly chosen", the results are not volume-conserving nor energy stable. The time step size and space time size are chosen as $\Delta t=1\times 10^{-5}$, $h_x=h_y=1/256$ in the simulations, respectively. The initial conditions and other parameters are chosen the same as those in \ref{Fig2}.} \label{Fig8}
\end{figure*}

 The phase evolution, time evolution of the volume and the free energy of the two Allen-Cahn models with nonlocal constraints simulated by EQ or SAV are essentially the same, if the penalizing parameter $\eta$  is chosen appropriately. The choice of $\eta$ can be fairly arbitrary  \cite{li2018unconditionally}. For very large  $\eta$, however,  the code does not perform well in that it slows down significantly. We thus believe there exists an "optimal" $\eta$ that renders the best result, which ought to be  determined empirically in numerical implementations. We show two examples of failure  at either large $\eta$ or a small $\eta$ in Fig \ref{Fig8}, in which the energy   is not dissipative nor the volume   conserved.
% \begin{figure}[htbp]
%  \centering
%  \includegraphics[width=1 \textwidth]{Figure8.jpg}
%  \caption{Two examples of failure for the nonlocal Allen-Cahn model with a penalizing term at $M=1$. $\eta=1$ and $\eta=1\times 10^{18}$ are chosen. We plotted the time evolution of free energy and volume, respectively. A and C are simulated by the EQ schemes for the penalizing model with $\eta=1$.  C and D are simulated by the SAV schemes for the penalizing model with $\eta=1\times 10^{18}$. If $\eta$ is not suitable, the results are not volume conserved and even not energy stable. The time step size and space time size are $\Delta t=1\times 10^{-5}$, $\Delta x=1/256$ and $\Delta y=1/256$. The initiation conditions and other parameters are same as that in \ref{Fig2}.} \label{Fig8}
%  \end{figure}

\begin{figure*}
\centering
\subfigure[]{
\begin{minipage}[b]{0.49\linewidth}
\centering
\includegraphics[width=0.75\linewidth]{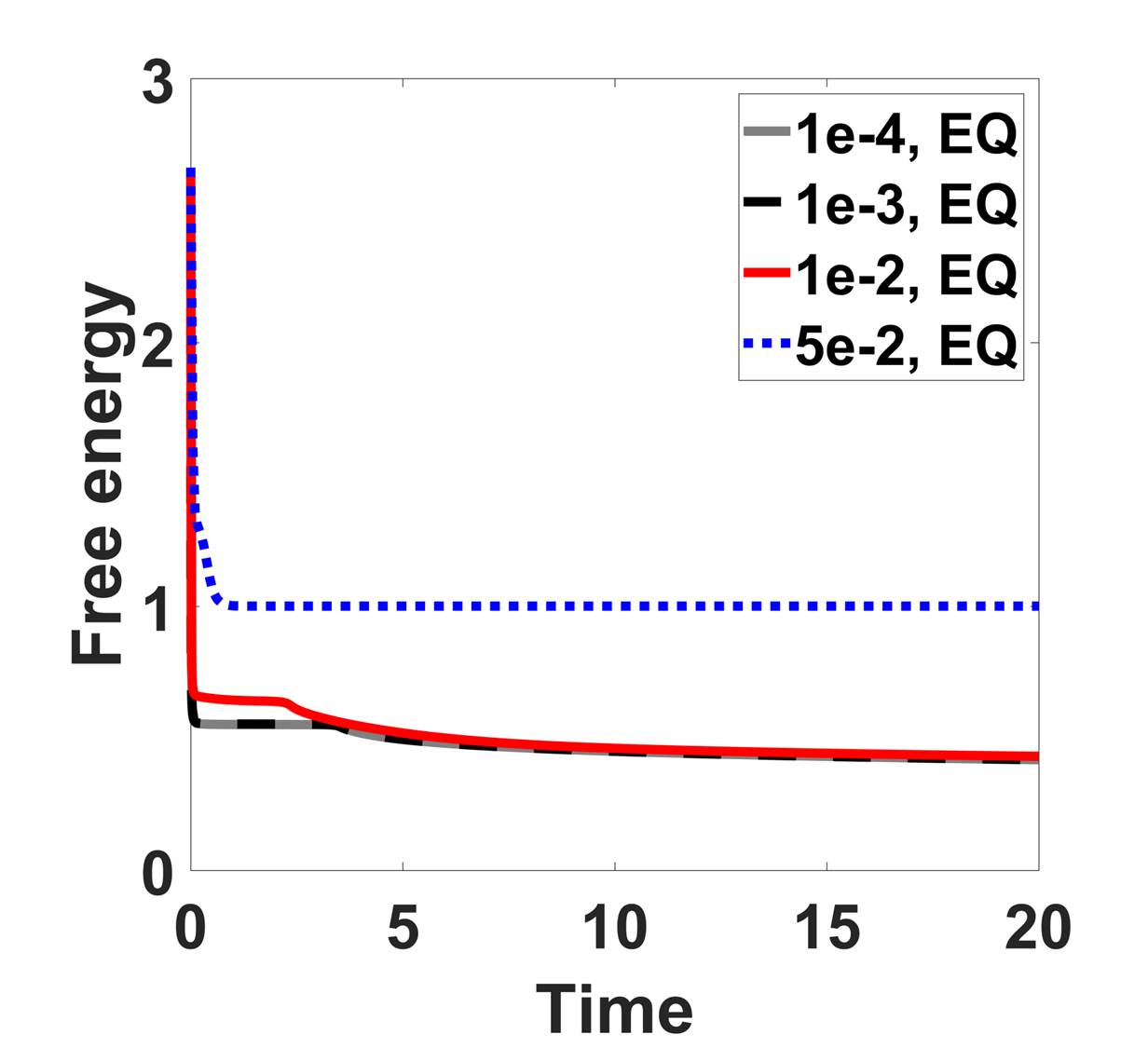}
\end{minipage}}
\subfigure[]{
\begin{minipage}[b]{0.49\linewidth}
\centering
\includegraphics[width=0.75\linewidth]{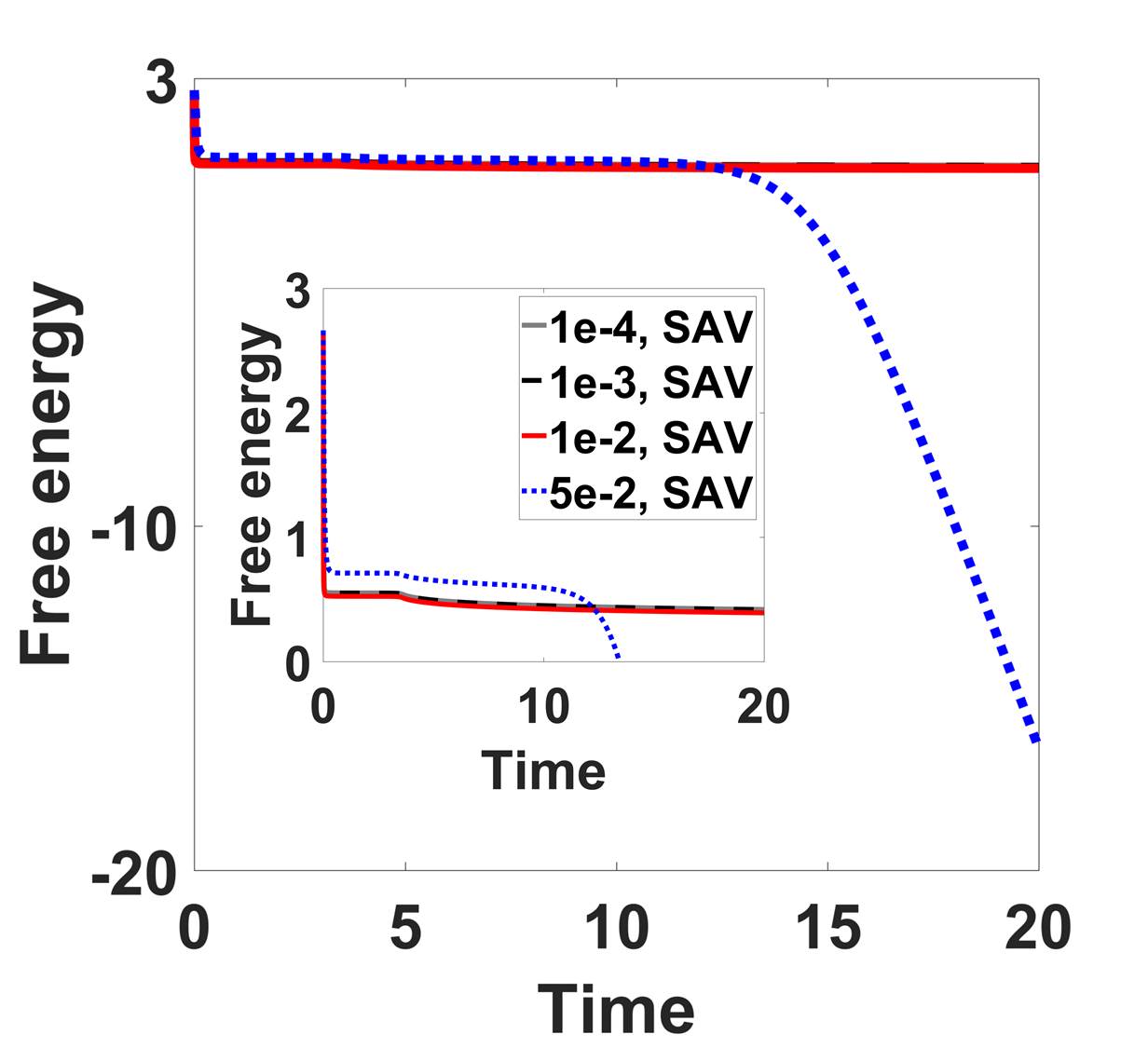}
\end{minipage}}
\caption{Merging of two droplets simulated by the AC-L1-EQ scheme in (a) and the AC-L1-SAV scheme in (b). The time step sizes used are $1\times 10^{-4}$, $1\times 10^{-3}$, $1\times 10^{-2}$, $5\times 10^{-2}$, respectively. The spatial mesh is $256 \times 256$ in 2D. The initial conditions and other parameter values are the same as those in \ref{Fig2}. At $\Delta t =0.05$, the scheme AC-L1-SAV would yield erroneous results eventually.  } \label{Fig9}
\end{figure*}

In the case with a double well potential, if we define $q=\sqrt{\gamma_2} \phi (1-\phi)$ rather than the one in  section 4.1 in the schemes, we observe a significant improvement in convergence at large time steps.
 This is because $q'=\vparl{q}{\phi}=\sqrt{\gamma_2}(1-2\phi)$ is a linear function, the extrapolation of $\ohs{q'}=\frac{3}{2}(q')^{n}-\frac{1}{2}(q')^{n-1}$ is exact. Compared with the results computed by previous schemes with large time step for nonlocal model in Fig \ref{Fig6}, larger time steps can perform well in the newly developed  EQ or SAV schemes in Fig \ref{Fig9}. Notice that EQ schemes outperforms   SAV schemes  at $\Delta t=5\times 10^{-2}$ but under-performs SAV schemes at $\Delta t=1\times 10^{-2}$. Their their performance is comparable.
%\begin{figure}[htbp]
%  \centering
%  \includegraphics[width=1 \textwidth]{Figure9.jpg}
%  \caption{Merging of two droplets performed by the new schemes for the Lagrangian model designed by EQ (A) and SAV(B) approaches. The time step sizes are $1\times 10^{-4}$, $1\times 10^{-3}$, $1\times 10^{-2}$, $5\times 10^{-2}$. The mesh is 256 $\times$ 256. The initiation conditions and other parameters are same as that in \ref{Fig2}.} \label{Fig9}
%  \end{figure}

%The two implementation tricks  are also applied to the CH-SAV and AC-L1-SAV scheme with large time step size, however, only method 1 can work well (See Fig \ref{FigS1}).\\
%\begin{figure*}
%\centering
%\subfigure[]{
%\begin{minipage}[b]{0.49\linewidth}
%\includegraphics[width=1\linewidth]{FigureS1A.jpg}
%\end{minipage}}
%\subfigure[]{
%\begin{minipage}[b]{0.49\linewidth}
%\includegraphics[width=1\linewidth]{FigureS1B.jpg}
%\end{minipage}}
%\caption{The modified numerical implementations of SAV schemes for the Cahn-Hilliard model (a) and nonlocal Allen-Cahn model (b), respectively. We set $M=1$ for the  nonlocal %Allen-Cahn model and $M=1\times 10^{-4}$ for the Cahn-Hilliard model in the simulations. (a) and (b) are simulated by CH-SAV and AC-L1-SAV, respectively. (Method1: After %obtaining $\phi^{n+1}$, we update   $q^{n+1}=\sqrt{\int_\Omega f_1(\phi^{n+1}) \mathrm {d {\bf r}}+C_0}$ for the SAV schemes). The initiation conditions and parameters are same %as that in \ref{Fig2}.} \label{FigS1}
%\end{figure*}
\section{Conclusions}

\noindent \indent We have developed an exhaustive set of linear, second order, energy stable schemes for the  Allen-Cahn equation with nonlocal constraints that preserve the phase volume and compared them with the energy stable, linear schemes for the Allen-Cahn and the Cahn-Hilliard models. These schemes are devised based on the energy quadratization strategy in the form of EQ and SAV format. We show that they are unconditionally energy stable and uniquely solvable. All schemes can be solved using efficient methods, making the models   alternatives to the Cahn-Hilliard model to describe interface dynamics of immiscible materials while preserving the volume. The nonlocal Allen-Cahn models shows a slower coarsening rate than that of Cahn-Hilliard at the same mobility, but one can increase the mobility coefficient of the nonlocal Allen-Cahn model to accelerate their dynamics. Two implementation tricks are introduced to enhance the stability property as well as the accuracy of the numerical schemes at a large step size, but the second method doesn't work well on the SAV schemes. In addition, we have compared the two  Allen-Cahn models with nonlocal constraints numerically. The computational efficiency  of the Allen-Cahn model with a penalizing potential is slightly than the one with a Lagrange multiplier, but the accuracy of the former depends on a suitable choice of model parameter $\eta$. Through numerical experiments, we show that the practical implementation using the defining algebraic functions for the auxiliary variable makes the  EQ scheme superior to the SAV scheme in their computational efficiency and accuracy. When the equation of the auxiliary variable can be solved more accurately, large time step size can be applied. In the end, we note that the size of mobility and the time step size are dominating factors that determine ultimately the efficiency and accuracy of the schemes.

\section*{Acknowledgements}
 %Jun Li's work is supported by the National Natural Science Foundation of China (Grant No. 11301287).
 Qi Wang's research is partially supported by NSFC awards \#11571032, \#91630207 and NSAF-U1530401.
\section*{Appendix}
\begin{appendices}
\section{Sherman-Morrison formula and its application to solving the integro-differential equation}
\noindent \indent Here we give a brief review on the Sherman-Morrison formula \cite{bodewig1959matrix} and explain its applications in the practical implementation of our various relevant schemes.

Suppose $A$ is an invertible square matrix, and $u$,$v$ are column vectors. Then $A+uv^T$ is invertible iff $1+v^TA^{-1}u \neq 0$. If $A+uv^T$ is invertible, then its inverse is given by
\bena
(A+uv^T)^{-1}=A^{-1}-\frac{A^{-1}uv^TA^{-1}}{1+v^TA^{-1}u}.
\eena
So
if $Ay=b$ and $Az=u$, $(A+uv^T)x=b$ has the solution given by
\bena
x=y-\frac{v^Ty}{1+v^Tz}z.
\eena

For the integral term(s) in the  semi-discrete schemes in this study such as (\ref{Sherman-M}), we need to discretize it properly. $\forall f$, we discretize $\int_\Omega f \mathrm{d{\bf r}}$ using the composite trapezoidal rule and adding all the elements of the new matrix $w_1 w_2^T f$, where $w_1=\frac{h_x}{2}S$, $w_2=\frac{h_y}{2}S$, $h_x$, $h_y$ are the spatial step sizes and $S={[1,2,2,...,2,2,1]}^T$. For convenience, we use $w_1w_2^T f$ to represent the integral discretized by the composite trapezoidal rule.

To solve equation (\ref{Sherman-M}), we  discretize the integral
 or the inner product of functions $(c,\phi^{n+1})d $ as $ u { v}^T \phi^{n+1}$. The scheme is recast to $A\phi^{n+1}+u { v}^T \phi^{n+1}=b^n$.  After using the Sherman-Morrison formula, we get
\bena
\phi ^{n+1}= A^{-1} {b}^n-\frac{{v}^T{ A}^{-1} { {b}}^n}{{1+{ v}^T  A}^{-1}  u}{ A}^{-1}  u,
\eena
In the  inner product of  vectors,
(\ref{Sherman-M}) can be rewritten into
\bena
{\phi} ^{n+1}= A^{-1}{b}^n-\frac{\langle c,A^{-1}b^n\rangle }{1+\langle c,A^{-1}d\rangle }{A}^{-1}d.
\eena
So, indeed the approach we take in the study using discrete inner product is essentially equivalent to applying the Sherman-Morrison formula.
%\section{The implementation of the SAV schemes}
% \begin{figure}[htbp]
%  \centering
%  \includegraphics[width=0.8 \textwidth]{FigureS1.jpg}
%  \caption{The modified numerical implementations of EQ schemes for the Cahn-Hilliard model(A) and nonlocal Allen-Cahn model(B), respectively. We set $M=1$ for the  nonlocal Allen-Cahn model and $M=1\times 10^{-4}$ for the Cahn-Hilliard model in the simulations. B is simulated by the SAV scheme for Lagrangian model. (Method1: After obtaining $\phi^{n+1}$, we update   $q^{n+1}=\sqrt{\int_\Omega f_1(\phi^{n+1}) \mathrm {d {\bf r}}+C_0}$ for the SAV schemes). The initiation conditions and parameters are same as that in \ref{Fig2}.} \label{FigS1}
%\end{figure}
%\noindent \indent
%\clearpage
%\section{}
%\subsection{}
%\section{}
%\subsection{}

\end{appendices}

\bibliographystyle{plain}
\bibliography{mybibtex}
\end{document}